\documentclass[10pt,a4paper]{amsproc}
\usepackage{amssymb,amsthm,amsmath}
\usepackage[english]{babel}
\usepackage{fullpage}
\numberwithin{equation}{section}

\title{Laplacians and gauged Laplacians\\~ on a quantum Hopf bundle}

\date{21 March  2010}

\author{~\\~\\~\\~
\large{Alessandro Zampini}
\\ [10pt]
~\\~\\
\normalsize{Max Planck Institut f\"ur  Mathematik - Bonn,} \\
\normalsize{Vivatsgasse 7, D-53111 Bonn, Germany.}\\
\normalsize{ {\tt zampini@mpim-bonn.mpg.de}}
\\~\\~\\~
\normalsize{Hausdorff Zentrum f\"{u}r Mathematik}\\
\normalsize{der Universit\"{a}t Bonn,}\\
\normalsize{Endenicher Allee 62, D-53115 Bonn, Germany.}
}

\newtheorem{theo}{Theorem}[section]
\newtheorem{lemm}[theo]{Lemma}
\newtheorem{prop}[theo]{Proposition}

\newtheorem{rema}[theo]{Remark}

\newcommand{\nn}{\nonumber}
\newcommand{\ce}{\mathcal{E}}
\newcommand{\dd}{{\rm d}}
\newcommand{\ca}{\mathcal{A}}

\newcommand{\ch}{\mathcal{H}}
\newcommand{\cl}{\mathcal{L}}
\newcommand{\cn}{\mathcal{N}}
\newcommand{\cp}{\mathcal{P}}
\newcommand{\cb}{\mathcal{B}}
\newcommand{\cq}{\mathcal{Q}}
\newcommand{\oca}[1]{\Omega^{#1}(\ca)}
\newcommand{\och}[1]{\Omega^{#1}(\ch)}
\newcommand{\cu}{\mathcal{U}}        
\newcommand{\SU}{\mathrm{SU}_q(2)}  
\newcommand{\ASU}{\ca(\mathrm{SU}_q(2))}  
\newcommand{\sq}{\mathrm{S}^2_{q}}  
\newcommand{\Asq}{\ca(\mathrm{S}^2_{q})}  
\newcommand{\su}{\cu_q(\mathfrak{su}(2))}  

\newcommand{\eps}{\varepsilon}      
\newcommand{\cop}{\Delta}           
\newcommand{\co}[2]{#1_{(#2)}}      
\newcommand{\hs}[2]{\left\langle #1,#2\right\rangle}  

\newcommand{\ket}[1]{\left | #1 \right\rangle }
\newcommand{\bra}[1]{\left\langle #1 \right |}
\newcommand{\oh}{{\tfrac{1}{2}}}
\newcommand{\shalf}{{\scriptstyle\frac{1}{2}}} 
\newcommand{\half}{{\mathchoice{\oh}{\oh}{\shalf}{\shalf}}} 
\newcommand{\lt}{{\triangleright}}    
\newcommand{\rt}{{\triangleleft}}
\newcommand{\IC}{{\mathbb C}} 
\newcommand{\IR}{{\mathbb R}} 
\newcommand{\IN}{{\mathbb N}} 
\newcommand{\IZ}{{\mathbb Z}} 
\DeclareMathOperator{\Ad}{Ad}       
\DeclareMathOperator{\id}{id}       
\DeclareMathOperator{\U}{U}       

\newcommand{\abs}[1]{\left|#1\right|}

\newcommand{\figureheight}{8cm}
\newcommand{\putfig}[2]{\begin{figure}[htp]
        \special{isoscale c:/itex/texfig/#1.wmf, \the\hsize \figureheight}
        \vspace{\figureheight}
        \caption{#2}\label{fig:#1}
        \end{figure}}
\newcommand{\pictureheight}{4cm}
\newcommand{\putpicture}[2]{\begin{figure}[htp]
        \special{isoscale c:/itex/texfig/#1.wmf, \the\hsize \pictureheight}
        \vspace{\pictureheight}
        \caption{#2}\label{fig:#1}
        \end{figure}}

\newcommand{\beqa}{\begin{eqnarray}}
\newcommand{\eeqa}{\end{eqnarray}}
\newcommand{\beq}{\begin{equation}}
\newcommand{\eeq}{\end{equation}}

\newcommand{\del}{\partial}

\newcommand{\zed}{{\bb Z}} 

\newcommand{\idop}{{\bf 1}} 

 \font\mybb=msbm10 at 12pt
\def\bb#1{\hbox{\mybb#1}}

\newcommand{\Psin}{\Psi^{\left(n\right)}}
\newcommand{\Psint}{\tilde{\Psi}^{\left(n\right)}}

\newcommand{\mn}{\abs{n}}
\newcommand{\mj}{\abs{j}}

\newcommand{\qpp}{\mathfrak{p}^{\left(n\right)}}
\newcommand{\qpn}{\check{\mathfrak{p}}^{\left(n\right)}}
\newcommand{\qpt}{\tilde{\mathfrak{p}}^{\left(n\right)}}
\newcommand{\bz}{B_{0}}
\newcommand{\bp}{B_{+}}
\newcommand{\bm}{B_{-}}
\newcommand{\delb}{\bar{\del}}

\newcommand{\e}{\mathrm{e}}
\newcommand{\f}{\mathrm{f}}
\newcommand{\h}{\mathrm{h}}

\newcommand{\Ks}{\mathrm{K}}
\newcommand{\A}{\mathrm{a}}
\newcommand{\An}{\mathrm{A^{(n)}}}
\newcommand{\Apn}{\mathrm{A^{\prime(n)}}}
\newcommand{\mcl}{\mathfrak{L}}
\newcommand{\zce}{\mathcal{E}^{(0)}}
\newcommand{\bsigma}{\bra{\sigma}}

\begin{document}
\maketitle

\thispagestyle{empty}

\maketitle

\begin{abstract}
\noindent
This paper presents an analysis of the set of connections and covariant derivatives on a $\U(1)$ quantum Hopf bundle on the standard quantum sphere $\sq$, whose total space algebra $\SU$ is equipped with the 3d left covariant differential calculus by Woronowicz. The introduction of a Hodge duality on both $\Omega(\SU)$ and on $\Omega(\sq)$ allows for the study of Laplacians and of gauged Laplacians.
\end{abstract}
~\\~
\begin{center}
{\sl This paper is dedicated to Sergio Albeverio, on the occasion of his 70th birthday.} 
\end{center}
\newpage

\section{Introduction}

This paper is focussed on the analysis of  a class of Hall Hamiltonians in the noncommutative set up. 
It is intended as a survey of the general formulation of quantum principal bundles, and as a description of a specific procedure to introduce, on both the total space and the base space of a quantum Hopf bundle, a set of Laplacian operators and to couple them with gauge connections.  It also presents a detailed formulation of the classical Hopf bundle. The emphasis in the presentation of structures from   differential geometry will be given to their algebraic aspects  extended to the noncommutative setting.


Classical Hall Hamiltonians are  gauged Laplace operators acting on the space of sections of the vector bundles associated to the principal bundles $\pi:G\to G/\Ks$ over homogeneous spaces (with $G$ semisimple and $\Ks$ compact) and can be constructed in terms of the Casimir operators of $G$ and $\Ks$. With $(\rho,V)$ a representation of $\Ks$, one has the identification of sections of the associated vector bundle $E=G\times_{\rho(\Ks)}V$ with equivariant maps from $G$ to $V$, $\Gamma(G/\Ks,E)\simeq C^{\infty}(G,V)_{\rho(\Ks)}\subset C^{\infty}(G)\otimes V$. Given a connection on $G$ one has a covariant derivative $\nabla$ on $\Gamma(G/\Ks,E)$, so that the gauged Laplacian operator is $\Delta^{E}=(\nabla\nabla^{*}+\nabla^{*}\nabla)=\star\nabla\star\nabla$, where the dual $\nabla^{*}$ is defined from the metric induced on the homogeneous space basis $G/\Ks$ by the Cartan-Killing  metric on $G$, or equivalently the Hodge duality comes from the induced metric on $G/\Ks$. If the connection is the canonical one, given by the orthogonal splitting of the Lie algebra $\mathfrak{g}$ of $G$ in terms of the Lie algebra $\mathfrak{k}$ of the gauge group and of its orthogonal complement, then the gauged Laplacian operator can be cast in terms of the quadratic Casimirs of $\mathfrak{g}$ and $\mathfrak{k}$:
\beq
\Delta^{E}=\left.(\Delta^{G}\otimes 1-1\otimes C_{\mathfrak{k}})\right|_{C^{\infty}(G,V)_{\rho(\Ks)}}=\left.(C_{\mathfrak{g}}\otimes 1-1\otimes C_{\mathfrak{k}})\right|_{C^{\infty}(G,V)_{\rho(\Ks)}}
\label{uno}
\eeq
The  above formula \cite{bgv}  simplifies the diagonalisation of the gauged Laplacian, and has  important applications in the study of the heat kernel expansion and index theorems on principal bundles.

The natural evolution is to develop  models of the Hall effect on noncommutative spaces whose symmetries are described in terms of quantum groups.  In \cite{lrz} the first model of 'excitations moving on a quantum 2-sphere' in the field of a magnetic monopole has been studied. 
It is described by quantum principal $\U(1)$-bundle over a quantum sphere $\sq$ having as a total space the manifold of the quantum group $\SU$ \cite{bm93}. The natural associated line bundles are classified by the winding number $n\in\zed$: equipped $\SU$ with the three dimensional left covariant calculus from Woronowicz \cite{wor89}, the gauge monopole connection is studied and a gauged Laplacian acting on sections of the associated bundle is completely diagonalised. That paper presents a first generalisation of the relation \eqref{uno}. Its most interesting aspect is that the corresponding energies are not invariant under the exchange monopole/antimonopole, namely the spectrum of the gauged Laplacian  is not invariant under the inversion of the direction of the magnetic field, a manifestation of the phenomenon usually referred to as  'quantisation removes degeneracy'. A parallel study of the relation \eqref{uno} is presented in \cite{dala}, where Laplacians on a quantum projective plane are gauged via the monopole connection. 
 
The analysis in \cite{lrz} embodies two specific starting points. The first one is that the quantum Casimir $C_{q}$ for the universal envelopping algebra $\mathcal{U}_{q}(\mathfrak{su(2)})$ dual to $\SU$ -- thus playing the quantum role of the classical envelopping algebra dual to the classical Lie group --  \emph{is}  a quadratic operator in the generators of $\mathcal{U}_{q}(\mathfrak{su(2)})$ acting on $\SU$, but \emph{can not} be cast in the form of a whatever rank polynomial in the left invariant generators of the left invariant three dimensional differential calculus by Woronowicz, so to say in the basis of  natural left invariant derivations associated to this differential calculus. The second starting point is given by the   studies performed in \cite{ma05}. In that paper a $\star$-Hodge operator on the exterior algebra on the Podle\'s sphere $\sq$ -- coming from  the differential two dimensional calculus induced on $\sq$  by the three dimensional calculus on $\SU$ -- had been introduced, so to make it possible the  definition of  a Laplacian operator on $\sq$.

This paper develops the analysis started in \cite{lrz}, and describes another generalisation of the relation \eqref{uno} to the setting of the same  quantum Hopf bundle. A family of compatible $\star$-Hodge structures on the exterior algebras $\Omega(\SU)$ and  $\Omega(\sq)$, depending on a set of real parameters, are introduced, giving  the corresponding Laplacians $\Box_{\SU}=\star\dd\star\dd:\ASU\mapsto\ASU$, and $\Box_{\sq}=\star\dd\star\dd:\Asq\mapsto\Asq$. The connections on the principal bundle allows for a gauging of the Laplacian $\Box_{\sq}$ on each associated line bundle. When $\Box_{\sq}$ is gauged into $\Box_{D_{0}}$ via the monopole connection, one finds
\beq
q^{2n}\Box_{D_{0}}=(\Box_{\SU}+\gamma X_{z}X_{z}),
\label{unob}
\eeq
where the integer $n\in\,\IZ$ specifies the value of the monopole charge. This is the relation generalising  the first equality in \eqref{uno}: the role of the quadratic Casimir of the gauge group algebra is played by $\gamma X_{z}X_{z}\lt$, with $X_{z}$ the vertical derivation of the fibration, and $\gamma\in\,\IR_{+}$ appears in this formulation as a parametrisation for a set of compatible $\star$-Hodge structures  giving Laplacians satisfying the same relation~\eqref{unob}.

This paper begins with an  exposition of the classical  Hopf bundle $\pi:S^{3}\to S^2$. Section \ref{se:ic} presents a global -- i.e. charts independent -- description of the differential calculi on both the Lie group manifold $SU(2)\simeq S^{3}$ and on the homogeneous space $S^{2}=S^{3}/\U(1)$, and introduces on the exterior algebras $\Omega(S^{3})$ and $\Omega(S^{2})$ the Hodge duality structures coming from a Cartan-Killing type metric on the Lie algebra $\mathfrak{su(2)}$, in order to define Laplacian operators. The principal bundle structure is described in terms of a well known principal bundle atlas.
The aim of the section is to explicitly compute for such a specific  Hopf bundle, following the classical approach from differential geometry, the main structures which will be generalised to the quantum setting. A more general and complete analysis of a noncommutative geometry approach to the differential geometry of principal and quantum bundles is in \cite{qpol}. 

Section \ref{se:qphb} describes  the quantum formulation \cite{bm93} of the principal bundle having $\ASU$ as total space algebra, $\Asq$ as base manifold algebra and $\ca(\U(1))$ as gauge group algebra, with the differential calculus on $\SU$ given the 3d left-covariant calculus introduced by Woronowicz \cite{wor87,wor89}. 

Section \ref{se:hssl} presents a $\star$-Hodge duality on $\Omega(\SU)$, allowing for the definition of a Laplacian operator. The Hodge duality is introduced following \cite{kmt}; section \ref{se:hlss} describes an evolution of this approach, giving a $\star$-Hodge duality structure on $\Omega(\sq)$, and analysing its  compatibility with the one on $\Omega(\SU)$. 

Section \ref{se:conn} provides a complete explicit description of the set of connections on this specific realisation of the quantum Hopf bundle, and of the main properties of the covariant derivative operators on each associated line bundle. The emphasis is on the domain of the covariant derivative operators -- the set of horizontal coequivariant elements of the bundle -- which appears here as the quantum counterpart of the classical forms also called tensorial forms.
Section \ref{se:gL} studies the coupling of the Laplacian operator on $\Omega(\sq)$ to the gauge connections.

Section \ref{limcla}  applies to the commutative algebras $\{\ca(\SU),\ca(S^{2}),\ca(\U(1))\}$ the  formalism developed in the quantum setting, in order to recover the structure of the classical Hopf bundle from an algebraic perspective.
Section \ref{bproj} closes the paper with an evolution of section \ref{se:conn}, describing how a covariant derivative operator can be defined on $\Omega(\SU)$, the whole exterior  algebra on the total space $\SU$ of the quantum Hopf bundle, following the formalism developed in \cite{durI,durII}.  

\section{The  classical Hopf bundle}
\label{se:ic}

The first formulation of what are nowadays known as Hopf fibrations is contained in \cite{hopf31,hopf35} in terms of projecting spheres to spheres of lower dimensions: it came also as a geometric formulation of the Dirac's model of  magnetic monopole  \cite{dir31}. The following lines are intended as a concise introduction to the formalism of fiber -- and principal -- bundles, aimed to set the notations that will be used in this paper: excellent textbooks -- like for example \cite{hus,michor} -- deeply and extensively describe this subject.

\bigskip

With $\pi:\mathcal{P}\to\mathcal{M}$  a smooth surjective map from a manifold $\mathcal{P}$ to a manifold $\mathcal{M}$, $(\mathcal{P},\mathcal{M},\pi)$ is a fibre bundle with typical fibre $\mathcal{F}$ over $\mathcal{M}$ if there is a fibre bundle atlas with charts  $(U_{i},\lambda_{i})$, where  $U_{i}$ is an open covering of   $\mathcal{M}$  and  the diffeomorphisms $\lambda_{i}:\pi^{-1}(U_{i})\to U_{i}\times\mathcal{F}$   are such that 
$\pi:\pi^{-1}(U_{i})\to U_{i}$ is the composition of $\lambda_{i}$ with the projection onto the first factor in $U_{i}\times \mathcal{F}$. The manifold $\mathcal{P}$ is called the total space of the bundle, the manifold $\mathcal{M}$ is the base of the bundle. From the definition it follows that $\pi^{-1}(m)$ is diffeomorphic to $\mathcal{F}$ -- the fibre of the bundle -- for any $m\in\mathcal{M}$. For any $f\in\,\mathcal{F}$ one has   $\lambda_{i}\circ\lambda_{j}^{-1}(m,f)=(m,\lambda_{ij}(m,f))$ where $\lambda_{ij}:
(U_{i}\cap U_{j})\times \mathcal{F}\to\mathcal{F}$ is smooth and $\lambda_{ij}(m,\,\,)$ belongs to  the group Diff$(\mathcal{F})$ of diffeomorphisms of the fibre $\mathcal{F}$ for each $m\in\,U_{i}\cap U_{j}$. The mappings $\lambda_{ij}$ are called the transition functions of the bundle, and satisfy the cocycle condition $\lambda_{ij}(m,\,\,)\circ\lambda_{jk}(m,\,\,)=\lambda_{ik}(m,\,\,)$ for $m\in\,U_{i}\cap U_{j}\cap U_{k}$, with $\lambda_{ii}(m,\,\,)=id_{\mathcal{F}}$ for $m\in\, U_{i}$. 

A fibre bundle $(\mathcal{P},\mathcal{M},\pi)$ is called a vector bundle if its typical fibre $\mathcal{F}$ is a vector space and if the trivialisation diffeomorphisms $\lambda_{i}$  give transition functions   $\lambda_{ij}$  which are invertible linear maps, elements in $\mathrm{GL}(\mathcal{F})$ for any $m\in\,\mathcal{M}$.
A principal bundle $(\mathcal{P},\Ks,[\mathcal{M}],\pi)$ with structure group $\Ks$ is a fibre bundle $(\mathcal{P},\mathcal{M},\pi)$ with typical fibre $\Ks$ and  transition functions $\lambda_{ij}(m,\,\,)\in\,$ Aut$(\Ks)$ which give the left translation of the group $\Ks$ on itself.  On the total space of a principal bundle there is also a right action of the Lie group $\Ks$  -- that is $\mathrm{r}_{k^{\prime}}(\mathrm{r}_{k}(p))=\mathrm{r}_{kk^{\prime}}(p)$ for any $p\in\,\mathcal{P}$ and $k,k^{\prime}\in\,\Ks$ -- such that $\pi(\mathrm{r}_{k}(p))=\pi(p)$,  and such that the action is free and transitive. The base $\mathcal{M}$ of the bundle can be identified with the quotient $\mathcal{P}/\Ks$ with respect to such a right action.

Given $G$ a Lie group and $\Ks\subset G$ a closed Lie subgroup of it, the group manifold $G$ is the total space manifold of a principal bundle $(G, \Ks, G/\Ks,\pi)$ with base space $G/\Ks$ - the space of left cosets - and typical fiber given by the structure or gauge group $\Ks$, so that the bundle projection $\pi:G\to G/\Ks$ is the canonical projection. The right principal action of the gauge group $\Ks$ on $G$ is given as $\mathrm{r}_{k}(g)=gk$ for any $k\in\Ks$ and $g\in G$. This action trivially satisfies the requirements of being free and transitive.   If $\mathfrak{k}$ is the Lie algebra of the group $\Ks$, the fundamental vector field $X_{\tau}\in\mathfrak{X}(G)$ associated to $\tau\in\mathfrak{k}$ is defined as the infinitesimal generator of the right principal action $\mathrm{r}_{\exp s\tau}(g)=g\exp s\tau$ of the one parameter subgroup $\exp s\tau\subset \Ks$: the  mapping $\tau\in\mathfrak{k}\to\{X_{\tau}\}\in\mathfrak{X}(G)$ is a Lie algebra isomorphism between $\mathfrak{k}$ and the set of fundamental vector fields $\{X_{\tau}\}$.  A differential form $\phi\in\,\Omega(G)$ is called horizontal if $i_{X_{\tau}}\phi=0$ for any fundamental vector field $X_{\tau}$.

If $\rho:\Ks\to \mathrm{GL}(W)$ is a finite dimensional representation of $\Ks$ on the vector space $W$, the associated vector bundle to $G$ is  the vector bundle whose total space is  $\mathcal{E}=G\times_{\rho(\Ks)}W$, having typical fiber $W$. 
It is defined as the quotient of the product $G\times W$  by the equivalence relation $(\mathrm{r}_{k}(g)=gk;w)\sim(g;\rho(k)\cdot w)$ for any choice of $g\in\,G$, $k\in\,\Ks$ and $w\in\,W$: $(\mathcal{E},G/\Ks,\pi_{\mathcal{E}})$ is a fibre bundle with a projection $\pi_{\mathcal{E}}:\mathcal{E}\to G/\Ks$ which is consistently defined on the quotient as $\pi_{\mathcal{E}}[g,w]_{\rho(\Ks)}=\pi(g)$ from the principal bundle projection $\pi$. 

With $\mathrm{r}_{k}^{*}:\Omega(G)\to\Omega(G)$ the action of $\Ks$ on the exterior algebra $\Omega(G)$ induced as a pull-back of the right action $\mathrm{r}_{k}$ of $\mathrm{K}$ on $G$, the $\rho(\Ks)$-equivariant $r$-forms  of the  principal bundle are $W$-valued forms on $G$  defined as:
\beq
\Omega^{r}(G,W)_{\rho(\Ks)}=\{\phi\in \Omega^{r}(G,W)=\Omega^{r}(G)\otimes W: \mathrm{r}_{k}^{*}(\phi)=\rho^{-1}(k)\phi\}.
\label{Omrg}
\eeq
A section of the associated bundle $\mathcal{E}$ is an element in $\Gamma(G/\Ks,\mathcal{E})$, namely a  map $\sigma:G/\Ks\to\mathcal{E}$ such that $\pi_{\mathcal{E}}(\sigma(m))=m$ for any $m\in\,G/\Ks$. This definition is extended to $\Gamma^{(r)}(G/\Ks,\mathcal{E})$, the set of $r$-forms on the basis $G/\Ks$ of the principal bundle with values in $\mathcal{E}$. There is a canonical isomorphism 
\beq
\Gamma^{(r)}(G/\Ks,\mathcal{E})\simeq \Omega^{r}_{\mathrm{hor}}(G,W)_{\rho(\Ks)}
\label{isgec}
\eeq
from the space of $\mathcal{E}$-valued differential forms on $G/\Ks$ onto the space of horizontal $\rho(\Ks)$-equivariant $W$-valued differential forms on the principal bundle $(G,\Ks,\pi)$. For $r=0$ -- with $\Gamma(G/\Ks,\ce)\simeq\Gamma^{(0)}(G/\Ks,\ce)$ --  the isomorphism gives the well known equivalence between equivariant functions of a principal bundle and sections of its associated bundle.
In particular, for $W=\IR, \IC$ with trivial representation the isomorphism is 
\beq
\Omega(G/\Ks)\simeq\Omega_{\mathrm{hor}}(G)_{\rho(\Ks)=\Ks}=\{\phi\in\,\Omega(G):\,i_{X_{\tau}}\phi=0;\,\mathrm{r}_{k}^{*}\phi=\phi\},
\label{isotr}
\eeq
 giving a description of the exterior algebra on the basis of the principal bundle.

\bigskip

A connection on a principal bundle can be given via a connection 1-form.  A connection 1-form on $G$ is an element  $\omega\in\,\Omega(G,\mathfrak{k})$, taking values in $\mathfrak{k}$ and satisfying the  two local conditions:
$$
\omega(X_{\tau})=\tau,
$$
$$
\mathrm{r}_{k}^{*}(\omega)=\Ad_{k^{-1}}\omega,
$$ 
where the adjoint action of $\Ks$ is given by $(\Ad_{k^{-1}}\omega)(X)=k^{-1}\omega(X)k$ for any vector field $X\in\mathfrak{X}(G)$. At each point $g\in\,G$  there is on the tangent space $T_{g}G$  a natural notion of vertical subspace, whose basis is given by the vectors $X_{\tau}$ which are  tangent to the fiber group $\Ks$, while the connection 1-form selects the horizontal subspace $H^{(\omega)}_{g}(G)$ given by the kernel of $\omega$. Identifiying the element $\omega(X)\in\,\mathfrak{k}$ with the vertical vector field it generates, the expression  $X^{(\omega)}=X-\omega(X)$ denotes the horizontal projection of the vector field $X\in\mathfrak{X}(G)$.  

Given any $\rho(\Ks)$-equivariant form $\phi\in\Omega^{r}(G,W)_{\rho(\Ks)}$, the covariant derivative is defined as the map:
\beq
D:\Omega^{r}(G,W)_{\rho(\Ks)}\to\Omega^{r+1}_{\mathrm{hor}}(G,W)_{\rho(\Ks)},\qquad D\phi(X_{1},\ldots,X_{r+1})=\dd\phi(X_{1}^{(\omega)},\ldots,X_{r+1}^{(\omega)})
\label{2p1}
\eeq 
where $\dd$ is the exterior derivative on $G$. On a $\rho(\Ks)$-equivariant horizontal form $\phi\in\,\Omega_{\mathrm{hor}}(G,W)_{\rho(\Ks)}$ the action of the covariant derivative can be written in terms of the connection 1-form as:
\beq
D\phi=\dd\phi+\omega\wedge\phi.
\label{2p2}
\eeq

\bigskip

The following sections describe the Hopf fibration  $\pi:S^{3}\mapsto S^{2}$, with $G\simeq SU(2)$, $\Ks\simeq U(1)$ and $S^{2}$ the space of the orbits $SU(2)/\U(1)$, and the monopole connection.

\subsection{A differential calculus on the classical $SU(2)$ Lie group}\label{se:cdc}

The aim of this section is to describe the differential calculus on the total space of this bundle, in terms of a natural basis of global vector fields and 1-forms \cite{michor}. It is intended to give them an explicit expression in order to clarify  the classical limit of their quantum counterparts.

\bigskip

Recall that  a Lie group $G$ naturally acts on itself both from the right and from the left. The left action is the smooth map $\mathrm{l}:G\times G\to G$ defined via the left multiplication $\mathrm{l}(g^{\prime},g)=g^{\prime}g=\mathrm{l}_{g^\prime}(g)$: since $\mathrm{l}_{g^{\prime}g^{\prime\prime}}(g)=\mathrm{l}_{g^{\prime}}(\mathrm{l}_{g^{\prime\prime}} (g))$, the left action is a group homomorphism $\mathrm{l}_{g}:G\to \mathrm{Aut}(G)$. The right action is the smooth map $\mathrm{r}:G\times G\to G$ defined via the right multiplication $\mathrm{r}(g,g^{\prime})=gg^{\prime}=\mathrm{r}_{g^{\prime}}(g)$; it is then immediate to see that $\mathrm{r}_{g^{\prime}g^{\prime\prime}}(g)=gg^{\prime}g^{\prime\prime}=\mathrm{r}_{g^{\prime\prime}}(\mathrm{r}_{g^{\prime}}(g))$: the right action is a group anti-homomorphism $\mathrm{r}_{g}:G\to \mathrm{Aut}(G)$. For any $T\in\mathfrak{g}$, the Lie algebra of $G$, it is possible to define a vector field $R_{T}\in\mathfrak{X}(G)$. It acts  as  a derivation on a smooth complex valued function defined on $G$, and can be written in terms of  the pull-back $\mathrm{l}_{g}^{*}:C^{\infty}(G)\to C^{\infty}(G)$ induced by $\mathrm{l}_{g}$. On $\phi\in C^{\infty}(G)$:
\beq
R_{T}(\phi)=\frac{\dd}{\dd s}\left. (\mathrm{l}_{\exp sT}^{*}(\phi))\right|_{s=0}
\label{Rv}
\eeq
Although defined via the left action $\mathrm{l}_{g}$,  the vector field $R_{T}$ is called the right invariant vector field associated to $T\in\mathfrak{g}$; this set of fields owes its name to the fact that, given 
$\mathrm{r}_{g*}:\mathfrak{X}(G)\to\mathfrak{X}(G)$ the push-forward induced by the right action $\mathrm{r}_{g}$, they satisfy a property of right invariance as
$\mathrm{r}_{g*}(R_{T})=R_{T}$.
From the definition  of the pull-back map $\mathrm{l}_{g}^{*}:C^{\infty}(G)\to C^{\infty}(G)$ one has:  $$
\mathrm{l}_{g^{\prime}g^{\prime\prime}}^{*}(\phi)=\phi\circ\mathrm{l}_{g^{\prime}g^{\prime\prime}}=\phi\circ\mathrm{l}_{g^{\prime}}\circ\mathrm{l}_{g^{\prime\prime}}=\mathrm{l}_{g^{\prime\prime}}^{*}(\mathrm{l}_{g^{\prime}}^{*}(\phi))
$$ 
for any $\phi\in C^{\infty}(G)$. This  relation enables to prove that the map $\check{\mathrm{l}}:T\in\mathfrak{g}\to R_{T}\in\mathfrak{X}(G)$ is a Lie algebra anti-homomorphism, 
$[R_{T},R_{T^{\prime}}]=R_{[T^{\prime},T]}$.

The analogous definitions starting from the right action naturally hold. For any $T\in\mathfrak{g}$, the vector field $L_{T}\in\mathfrak{X}(G)$ is defined as a derivation on $C^{\infty}(G)$, namely as the infinitesimal generator of the  pull-back  $\mathrm{r}^{*}_{g}$ induced by the right action $\mathrm{r}_{g}$: 
\beq
L_{T}(\phi)=\frac{\dd}{\dd s}\left.(\mathrm{r}^{*}_{\exp sT}(\phi))\right|_{s=0}
\label{Lv}
\eeq
on any $\phi\in C^{\infty}(G)$. Left invariant vector fields satisfy a property of left invariance given as 
$\mathrm{l}_{g}^{*}(L_{T})=L_{T}$;
the map $\check{\mathrm{r}}:T\in\mathfrak{g}\to L_{T}\in\mathfrak{X}(G)$ is a Lie algebra homomorphism, with 
$[L_{T},L_{T^{\prime}}]=L_{[T,T^{\prime}]}$. 
The sets $\{L_{T}\},\,\{R_{T}\}$  are two  basis of the left free $C^{\infty}(G)$-module $\mathfrak{X}(G)$.

\bigskip

The total space of the classical Hopf bundle is the manifold $S^{3}$,
which represents the elements of the Lie group $SU\left(2\right)$. A 
point $g\in\,S^{3}$ can be then written  via a $2\times 2$ matrix
with complex entries and unit determinant: \beq
g\,=\,\left(\begin{array}{cc}
u & -\bar{v}\\
v & \bar{u}
\end{array}\right)\,\,\,:\,\,\,\bar{u}u+\bar{v}v=1;
\label{desu}
\eeq 
 the left invariant vector fields $\check{\mathrm{r}}(T)=L_{T}$ are given, following \eqref{Lv}, as the tangent vectors to the curves $g(s)=g\cdot\exp\,sT$. In the defining matrix representation it reads:
\beq
\frac{d}{ds}\left(\begin{array}{cc}
u & -\bar{v}\\
v & \bar{u}
\end{array}\right)\,\cdot\left(\exp sT\right)\mid_{s=0}\,=\,
\left(\begin{array}{cc}
u & -\bar{v}\\
v & \bar{u}
\end{array}\right)\,\cdot\,\left(T\right)
\label{liaction}
\eeq
Since $\exp sT$ is unitary, $T$ is antihermitian, and the choice of a basis  in terms of the Pauli matrices:
\beq
T_{x}=\frac{1}{2}\left(\begin{array}{cc} 0 & i \\ i & 0 \end{array}\right),
\qquad 
T_{y}=\frac{1}{2}\left(\begin{array}{cc} 0 & -1 \\ 1 & 0 \end{array}\right),
\qquad
T_{z}=\frac{i}{2}\left(\begin{array}{cc} 1 & 0 \\ 0 & -1 \end{array}\right),
\label{pau}
\eeq
gives the explicit form of the left invariant vector fields:
\begin{align} 
\label{lcf}
&L_{x}=-\frac{i}{2}\left(\bar{v}\frac{\del}{\del u}
-\bar{u}\frac{\del}{\del v}+u\frac{\del}{\del\bar{v}}-v\frac{\del}{\del\bar{u}}\right)\nonumber\\
&L_{y}=-\frac{1}{2}\left(\bar{v}\frac{\del}{\del u}
-\bar{u}\frac{\del}{\del v}-u\frac{\del}{\del\bar{v}}+v\frac{\del}{\del\bar{u}}\right)\nonumber\\
&L_{z}=\frac{i}{2}\left(u\frac{\del}{\del u}+v\frac{\del}{\del v}-
\bar{v}\frac{\del}{\del\bar{v}}-\bar{u}\frac{\del}{\del\bar{u}}\right)\nonumber\\
&L_{-}=L_{x}-iL_{y}\,=\,i\left(v\frac{\del}{\del\bar{u}}-u\frac{\del}{\del\bar{v}}\right)\nonumber\\
&L_{+}=L_{x}+iL_{y}\,=\,i\left(\bar{u}\frac{\del}{\del v}-\bar{v}\frac{\del}{\del u}\right),
\end{align}
satisfying the commutation relations:
\begin{align}
&\left[L_{z};L_{-}\right]=iL_{-}, \nn \\
&\left[L_{z};L_{+}\right]=-iL_{+}, \nn \\
&\left[L_{-};L_{+}\right]=2iL_{z}.
\label{clil}
\end{align} 
The components of the right invariant vector fields $R_{T}=\check{\mathrm{l}}(T)$ are then clearly  given in the defining matrix representation \eqref{Rv} as: 
\beq
\frac{d}{ds}\left(\exp sT\right)\cdot\left(\begin{array}{cc}
u & -\bar{v}\\
v & \bar{u}
\end{array}\right)\,\mid_{s=0}\,=\,\left(T\right)\,\cdot\,
\left(\begin{array}{cc}
u & -\bar{v}\\
v & \bar{u}
\end{array}\right)
\label{riaction} 
\eeq 
acquiring the form:
\begin{align} 
\label{rcf} 
&R_{x}=\frac{i}{2}\left(v\frac{\del}{\del
u}+u\frac{\del}{\del v}-
\bar{u}\frac{\del}{\del\bar{v}}-\bar{v}\frac{\del}{\del\bar{u}}\right)\nonumber\\
&R_{y}=-\frac{1}{2}\left(v\frac{\del}{\del u}-u\frac{\del}{\del v}-
\bar{u}\frac{\del}{\del\bar{v}}+\bar{v}\frac{\del}{\del\bar{u}}\right)\nonumber\\
&R_{z}=\frac{i}{2}\left(u\frac{\del}{\del u}-v\frac{\del}{\del v}+
\bar{v}\frac{\del}{\del\bar{v}}-\bar{u}\frac{\del}{\del\bar{u}}\right)\nonumber\\
&R_{-}=R_{x}-iR_{y}\,=\,i\left(v\frac{\del}{\del u}-\bar{u}\frac{\del}{\del\bar{v}}\right)\nonumber\\
&R_{+}=R_{x}+iR_{y}\,=\,i\left(u\frac{\del}{\del v}-\bar{v}\frac{\del}{\del\bar{u}}\right).
\end{align}
The commutation relations they satisfy are:
\begin{align}
&\left[R_{z};R_{-}\right]=-iR_{-}, \nn \\ 
&\left[R_{z};R_{-}\right]=iR_{+}, \nn \\
&\left[R_{-};R_{+}\right]=-2iR_{z}. 
\label{clir} 
\end{align} 
The quadratic Casimir of the Lie algebra $\mathfrak{su(2)}$ is written as 
\beq
C=\frac{1}{2}(L_{+}L_{-}+L_{-}L_{+})+L_{z}L_{z}=\frac{1}{2}(R_{+}R_{-}+R_{-}R_{+})+R_{z}R_{z}.
\label{Cc}
\eeq 
The set $\mathfrak{X}(S^3)$ is a free left $C^{\infty}(S^{3})$- module. Right vector fields can
be expressed in the basis of the left vector fields as  $R_{a}=J_{ab}L_{b}$. The matrix $J$ is given by:
\beq \left(\begin{array}{c} R_{-} \\ R_{z} \\ R_{+}
\end{array}\right)\,=\,
\left(\begin{array}{ccc} \bar{u}^{2} & 2\bar{u}v & -v^{2} \\ -\bar{u}\bar{v} & u\bar{u}-v\bar{v} & -uv \\
-\bar{v}^{2} & 2u\bar{v} & u^{2} \end{array}\right)\,
\left(\begin{array}{c} L_{-} \\ L_{z} \\ L_{+} \end{array} \right)
\label{Tcla}
\eeq
and its inverse matrix is:
\beq \left(\begin{array}{c} L_{-} \\ L_{z} \\
L_{+}
\end{array}\right)\,=\,
\left(\begin{array}{ccc} u^{2} & -2uv & -v^{2} \\ u\bar{v} & u\bar{u}-v\bar{v} & \bar{u}v \\
-\bar{v}^{2} & -2\bar{u}\bar{v} & \bar{u}^{2} \end{array}\right)\,
\left(\begin{array}{c} R_{-} \\ R_{z} \\ R_{+} \end{array} \right)
\eeq

A similar analysis can be performed in the study of the cotangent space $\mathfrak{X}^{*}(G)$ of a Lie group. This is a $C^{\infty}(S^{3})$-bimodule, with two basis of globally defined 1-forms, namely the left invariant $\{\tilde{\omega}_{a}\}$ dual to the set of left invariant vector fields $\{L_{a}\}$, and the right invariant $\{\tilde{\eta}_{b}\}$ dual to the set of right invariant vector fields $\{R_{b}\}$. They satisfy the invariance property:
\begin{align}
&\mathrm{l}_{g}^{*}(\tilde{\omega}_{a})=\tilde{\omega}_{a}, \nn \\
&\mathrm{r}_{g}^{*}(\tilde{\eta}_{b})=\tilde{\eta}_{b}:
\end{align}
one then immediately computes:
 \beq
R_{i}\,=\,J_{ij}L_{j}\,\,\Leftrightarrow\,\,\tilde{\eta}_{s}J_{sp}\,=\,\tilde{\omega}_{p}.
\label{LRcla}
\eeq 
The left invariant 1-forms are: 
\begin{align}
&\tilde{\omega}_{z}=-2i\left(\bar{u}du+\bar{v}dv\right) \nonumber \\
&\tilde{\omega}_{-}=-i\left(\bar{v}d\bar{u}-\bar{u}d\bar{v}\right) \nonumber \\
&\tilde{\omega}_{+}=-i\left(udv-vdu\right) 
\label{linv}
\end{align} 
with $\tilde{\omega}_{x}=(\tilde{\omega}_{-}+\tilde{\omega}_{+})$ and $\tilde{\omega}_{y}=i(\tilde{\omega}_{+}-\tilde{\omega}_{-})$. The antilinear involution  on $\Omega^{1}(S^{3})$, compatible with the antilinear involution on $C^{\infty}(S^{3})$, is given by $\tilde{\omega}_{x}^{*}=\tilde{\omega}_{x}$, $\tilde{\omega}_{y}^{*}=\tilde{\omega}_{y}$, $\tilde{\omega}_{z}^{*}=\tilde{\omega}_{z}$. The right-invariant
1-forms are: 
\begin{align}
&\tilde{\eta}_{z}=2i\left(ud\bar{u}+\bar{v}dv\right) \nonumber \\
&\tilde{\eta}_{-}=i\left(ud\bar{v}-\bar{v}du\right) \nonumber \\
&\tilde{\eta}_{+}=-i\left(\bar{u}dv-vd\bar{u}\right). 
\label{rinv}
\end{align} 
Given a complex valued smooth function  $\phi\in C^{\infty}(S^{3})$, the exterior derivative is the map $\dd:C^{\infty}(S^{3})\to\Omega^{1}(S^{3})$   defined via: 
\beq
\dd \phi(X)=X(\phi)
\eeq 
in terms of the Lie derivative $X(\phi)$ of $\phi$ along the vector field $X$.
This map acquires the form: 
\beq
\dd \phi=L_{a}(\phi)\tilde{\omega}_{a}=R_{b}(\phi)\tilde{\eta}_{b}
\label{decla}
\eeq
where now $L_{a}(\phi)$ represents the Lie derivative of $\phi$ along the vector field $L_{a}$, while $R_{b}(\phi)$ represents the Lie derivative of $\phi$ along the vector field $R_{b}$. 

From the $C^{\infty}(S^{3})$-bimodule $\Omega^{1}(S^{3})$ define the tensor product of forms as the $C^{\infty}(S^{3})$-bimodule $\{\Omega^{1}(S^{3})\}^{\otimes k}=\Omega^{1}(S^{3})\otimes_{C^{\infty}(S^{3})}\ldots\otimes_{C^{\infty}(S^{3})}\Omega^{1}(S^{3})$ ($k$ times). 
The exterior algebra coming from the differential calculus \eqref{decla} is defined as the graded associative algebra $\Omega(S^{3})=\left(\oplus_{k}\Omega^{k}(S^{3});\wedge\right)$, with $k$-forms and wedge product introduced in terms of an alternation mapping 
$\mathfrak{A}:\{\Omega^{1}(S^{3})\}^{\otimes k}\to\{\Omega^{1}(S^{3})\}^{\otimes k}$ \cite{AM}.
The wedge product is bilinear, and satisfies the identity $\alpha\wedge\beta=(-1)^{kl}\beta\wedge\alpha$ for any $k$-form $\alpha$ and $l$-form $\beta$. The complex involution is extended by requiring 
 $$
 (\alpha\wedge\beta)^{*}=(-1)^{kl}\beta^{*}\wedge\alpha^{*}.
 $$
The exterior derivative is extended to $\dd:\Omega^{k}(S^{3})\to\Omega^{k+1}(S^{3})$ as the unique $\IC$-linear mapping satisfying the conditions: 
\begin{enumerate}
\item $\dd$ is a graded $\wedge$-derivation, that is $\dd(\alpha\wedge\beta)=(\dd\alpha)\wedge\beta+(-1)^{k}\alpha\wedge\dd\beta$ for any $k$-form $\alpha$; 
\item  $\dd^{2}=\dd\circ\dd=0$; 
\item on $\phi\in\,\Omega^{0}(S^{3})\simeq C^{\infty}(S^{3})$,  it is given by $\dd \phi$ as in  \eqref{decla}.
\end{enumerate}
It is then easy to see that $\Omega^{2}(S^{3})$ is three dimensional, with a basis given by $\{\tilde{\omega}_{-}\wedge\tilde{\omega}_{+}, \tilde{\omega}_{+}\wedge\tilde{\omega}_{z}, \tilde{\omega}_{z}\wedge\tilde{\omega}_{-}\}$:
extending in a natural way via the pull back the left and right actions of the group $SU(2)$ on $\Omega^{2}(S^{3})$, it is also clear that such basis elements are left invariant. 
  From \eqref{linv} one has:
\begin{align}
&\dd\tilde{\omega}_{-}=i\tilde{\omega}_{-}\wedge\tilde{\omega}_{z}, \nn \\
&\dd\tilde{\omega}_{+}=-i\tilde{\omega}_{+}\wedge\tilde{\omega}_{z}, \nn \\
&\dd\tilde{\omega}_{z}=2i\tilde{\omega}_{-}\wedge\tilde{\omega}_{+}.
\label{d1fcl}
\end{align}
The bimodule $\Omega^{3}(S^{3})$ is one dimensional, with again a left invariant  basis 3-form given by $\{\tilde{\omega}_{-}\wedge\tilde{\omega}_{+}\wedge\tilde{\omega}_{z}\}$.  A right invariant basis of the exterior algebra $\Omega(S^{3})$ is analogously given in terms of the 1-forms $\tilde{\eta}_{a}$.

\subsection{A Laplacian operator on the group manifold $SU(2)$}\label{se:LS}

Being  $SU(2)$  a semisimple Lie group, the group manifold $S^{3}$ can be equipped with the Cartan-Killing metric originated from the Cartan decomposition of the Lie algebra $\mathfrak{su(2)}$.  
Consider now as a riemannian metric structure on $S^{3}$ the symmetric tensor 
\beq
g=\alpha(\tilde{\omega}_{x}\otimes\tilde{\omega}_{x}+\tilde{\omega}_{y}\otimes\tilde{\omega}_{y})+\tilde{\omega}_{z}\otimes\tilde{\omega}_{z},
\label{gmet}
\eeq
with $\alpha\in\,\IR^{+}$.   For $\alpha=1$ such a metric tensor coincides with the  the Cartan-Killing metric. The volume associated to the $\mathfrak{g}$-orthonormal basis and to the choice of the orientation $(x,y,z) $ is given by $\theta=\alpha\,\tilde{\omega}_{x}\wedge\tilde{\omega}_{y}\wedge\tilde{\omega}_{z}$, so that  $\theta^{*}=\theta$. Such a volume $\theta$ is a Haar volume, namely it is invariant with respect to both the left $\mathrm{l}_{g}^{*}$ and the right actions $\mathrm{r}_{g}^{*}$ of the Lie  group $SU(2)$ on itself, since an explicit calculation gives $L_{a}(\theta)=R_{a}(\theta)=0$. The Hodge duality $\star:\Omega^{k}(S^{3})\to\Omega^{3-k}(S^{3})$ which corresponds  to this volume \cite{AM} is the $C^{\infty}(S^{3})$-linear map given on the left invariant basis of the exterior algebra $\Omega(S^{3})$ by $\star(1)=\theta$, $\star(\theta)=1$, and:
\beq
\begin{array}{lcl}
\star(\tilde{\omega}_{x})=\tilde{\omega}_{y}\wedge\tilde{\omega}_{z}, & \qquad & \star(\tilde{\omega}_{y}\wedge\tilde{\omega}_{z})=\tilde{\omega}_{x}, \\
\star(\tilde{\omega}_{y})=\tilde{\omega}_{z}\wedge\tilde{\omega}_{x}, & \qquad & \star(\tilde{\omega}_{z}\wedge\tilde{\omega}_{x})=\tilde{\omega}_{y}, \\
\star(\tilde{\omega}_{z})=\alpha\,\tilde{\omega}_{x}\wedge\tilde{\omega}_{y}, & \qquad & \star(\tilde{\omega}_{x}\wedge\tilde{\omega}_{y})=\alpha^{-1}\,\tilde{\omega}_{z}. 
\end{array}
\label{clHs}
\eeq
The differential calculus on the group manifold $S^{3}$ as well as the above $\star$-Hodge duality on the exterior algebra $\Omega(S^{3})$ give a  Laplacian operator  defined as $\Box_{S^{3}}\phi=\star\dd\star\dd\phi$ on any $\phi\in\,C^{\infty}(S^{3})$. It can be written as a differential operator in terms of the left invariant vector fields:
\beq
\Box_{S^{3}}\phi=[\frac{1}{2\alpha}(L_{-}L_{+}+L_{+}L_{-})+L_{z}L_{z}]\phi
\label{clLa}
\eeq
The Laplacian operator is the Casimir of the Lie algebra $\mathfrak{su(2)}$ only if $\alpha=1$, that is only if the metric from where  it is originated is the Cartan-Killing metric.

The Hodge structure satisfies two  identities:
\begin{align}
&\star^{2}(\xi)=(-1)^{k(3-k)}\xi=\xi 
\label{clsH}
\\  
&\xi\wedge(\star\xi^{\prime})=\xi^{\prime}\wedge(\star\xi)
\label{symhs}
\end{align}
for any  $\xi,\xi^{\prime}\in\,\Omega^{k}(S^{3})$. This allows to define a symmetric bilinear map $\hs{~}{~}_{S^{3}}:\Omega^{k}(S^{3})\times\Omega^{k}(S^{3})\to C^{\infty}(S^{3})$ ($k=0,\ldots,3$) as:
\beq
\hs{\xi}{\xi^{\prime}}_{S^{3}}\theta=\xi\wedge(\star\xi^{\prime}).
\label{bifocl}
\eeq
It is clearly a symmetric tensor on $\{\mathfrak{X}^{*}(S^3)\}^{\otimes 2k}$, whose components can be expressed in terms of the components of the inverse metric $g^{-1}=g^{-1ab}L_{a}\otimes L_{b}\,\in\,\{\mathfrak{X}^{1}(S^{3})\}^{\otimes2}$ with $g^{-1ab}g_{bc}=\delta^{a}_{c}$, as
\beq
\hs{\tilde{\omega}_{i_{1}}\wedge\ldots\wedge\tilde{\omega}_{i_{k}}}{\tilde{\omega}_{j_{1}}\wedge\ldots\wedge\tilde{\omega}_{j_{k}}}_{S^{3}}=\,\sum_{\sigma}\pi_{\sigma}g^{-1i_{1}\sigma(j_{1})}\ldots g^{-1i_{k}\sigma(j_{k})}
\label{bifoclg}
\eeq
where the summation is over permutations $\sigma$ of k elements, with parity $\pi_{\sigma}$.  
Starting from the Hodge duality a second bilinear map  $\hs{~}{~}^{\sim}_{S^{3}}:\Omega^{k}(S^{3})\times\Omega^{k}(S^{3})\to C^{\infty}(S^{3})$, can be introduced as 
\beq
\hs{\xi^{\prime}}{\xi}_{S^{3}}^{\sim}\theta=\xi^{*}\wedge(\star\xi^{\prime})
\label{bifopr} 
\eeq
for any $\xi,\xi^{\prime}\in\,\Omega^{k}(S^{3})$, being hermitian $(\hs{\xi^{\prime}}{\xi}_{S^{3}}^{\sim})^*=\hs{\xi}{\xi^{\prime}}_{S^3}^{\sim}$. The Haar volume form can be used to introduce an integral on a manifold \cite{AM}, $\int_{\theta}:\Omega^3(S^3)\to\IC$; being $S^3$ compact, the volume of the group manifold can be normalised, setting $\int_{\theta}\theta=1$.  From \eqref{bifocl} and \eqref{bifopr} it is possible to define on the exterior algebra $\Omega(S^3)$ both a scalar product, 
\beq
(\xi;\xi^{\prime})_{S^{3}}=\int_{\theta}\xi\wedge(\star\xi^{\prime})=\int_{\theta}\hs{\xi}{\xi^{\prime}}_{S^{3}}\theta,
\label{symps}
\eeq
and an hermitian inner product, 
\beq
(\xi^{\prime};\xi)^{\sim}_{S^{3}}=\int_{\theta}\xi^{*}\wedge(\star\xi^{\prime})=\int_{\theta}\hs{\xi^{\prime}}{\xi}_{S^3}^{\sim}\theta.
\label{bifoprip}
\eeq
\noindent An evaluation on a non hermitian basis in $\Omega(S^3)$ presents the differences between the non vanishing terms of two bilinear forms: 
\begin{align}
&\hs{1}{1}_{S^{3}}=1; \nn \\
&\hs{\tilde{\omega}_{-}}{\tilde{\omega}_{+}}_{S^{3}}=
\hs{\tilde{\omega}_{+}}{\tilde{\omega}_{-}}_{S^{3}}=\frac{1}{2\alpha}, \qquad
\hs{\tilde{\omega}_{z}}{\tilde{\omega}_{z}}_{S^{3}}=1;  \nn \\
&\hs{\tilde{\omega}_{+}\wedge\tilde{\omega}_{z}}{\tilde{\omega}_{-}\wedge\tilde{\omega}_{z}}_{S^{3}}=
\hs{\tilde{\omega}_{-}\wedge\tilde{\omega}_{z}}{\tilde{\omega}_{+}\wedge\tilde{\omega}_{z}}_{S^{3}}=\frac{1}{2\alpha}, \qquad
\hs{\tilde{\omega}_{-}\wedge\tilde{\omega}_{+}}{\tilde{\omega}_{-}\wedge\tilde{\omega}_{+}}_{S^{3}}=\frac{1}{4\alpha^{2}};  \nn \\
&\hs{\theta}{\theta}_{S^{3}}=1;
\label{bifoclxyz}
\end{align}
while 
\begin{align}
&\hs{1}{1}^{\sim}_{S^{3}}=1; \nn \\
&\hs{\tilde{\omega}_{-}}{\tilde{\omega}_{-}}^{\sim}_{S^{3}}=
\hs{\tilde{\omega}_{+}}{\tilde{\omega}_{+}}^{\sim}_{S^{3}}=\frac{1}{2\alpha}, \qquad
\hs{\tilde{\omega}_{z}}{\tilde{\omega}_{z}}^{\sim}_{S^{3}}=1; \nn \\  
&\hs{\tilde{\omega}_{+}\wedge\tilde{\omega}_{z}}{\tilde{\omega}_{+}\wedge\tilde{\omega}_{z}}^{\sim}_{S^{3}}=
\hs{\tilde{\omega}_{-}\wedge\tilde{\omega}_{z}}{\tilde{\omega}_{-}\wedge\tilde{\omega}_{z}}^{\sim}_{S^{3}}=\frac{1}{2\alpha},\qquad
\hs{\tilde{\omega}_{-}\wedge\tilde{\omega}_{+}}{\tilde{\omega}_{-}\wedge\tilde{\omega}_{+}}^{\sim}_{S^{3}}=\frac{1}{4\alpha^{2}};  \nn \\
&\hs{\theta}{\theta}^{\sim}_{S^{3}}=1.
\label{bifoprpm}
\end{align}

\subsection{The principal bundle structure and the monopole connection}\label{se:pbc}

Consider  the one parameter subgroup of $SU(2)$ given by $\gamma_{T_{z}}(s)=\exp sT_{z}$ where $T_{z}$ is the generator in \eqref{pau}. In this specific  matrix representation it is 
\beq
\gamma_{T_{z}}(s)=\exp \left[\frac{is}{2}\left(\begin{array}{cc} 1 & 0 \\ 0 & -1 \end{array}\right)\right]=
\left(\begin{array}{cc} e^{is/2} & 0 \\ 0 & e^{-is/2} \end{array} \right),
\eeq
thus proving that $\gamma_{T_{z}}(s)\simeq \U(1)$ as a subgroup in $SU(2)$. The space of left cosets $SU(2)/\U(1)$  is the set of the orbits of the right principal action $\check{\mathrm{r}}_{\exp sT_{z}}(g)=g \exp sT_{z}$ which is free, and smooth; its infinitesimal generator coincides with the vector field $L_{z}$ \eqref{liaction}. As already mentioned the canonical projection $\pi:SU(2)\to SU(2)/\U(1)$ gives a principal bundle whose vertical field is $L_{z}$. A  formulation for a principal bundle atlas on a homogeneous space is extensively analysed in terms of local sections \cite{hus,michor}.   
This section describes in detail how  a principal bundle atlas is introduced \cite{gs87,tra} defining suitable trivialisations.  

Parametrise $S^{3}$ by 
\begin{align}
u&=\cos\theta/2\,e^{i(\varphi+\psi)/2} \nn \\
v&=\sin\theta/2\,e^{-i(\varphi-\psi)/2} \nn, 
\end{align}
with $0\leq\theta\leq\pi$ and $\phi,\psi\in\,\IR$, the Hopf map $\pi:SU(2)\to S^{2}\simeq SU(2)/\U(1)$ is defined by:
\begin{align}
&b_{z}=uu^{*}-vv^{*}=\cos\theta, \nn \\
&b_{y}=uv^{*}+vu^{*}=\sin\theta\cos\varphi, \nn \\
&b_{x}=-i(vu^{*}-uv^{*})=-\sin\theta\sin\varphi
\label{bdef}
\end{align}
with $b_{z}^{2}+b_{x}^{2}+b_{y}^{2}=1$.
It is immediate to see that $\pi(u,v)=\pi(u^{\prime},v^{\prime})$ if and only if $u^{\prime}=ue^{i\alpha}$ and $v^{\prime}=ve^{i\alpha}$ with $\alpha\in\IR$: this is also a way to recover that the projection has the standard fibre $\U(1)$. 
A choice for an open covering of the sphere $S^{2}$ is given by:
\begin{align}
&S^{2}_{(N)}=\{S^{2}:b_{z}\neq 1\}\qquad\Rightarrow\qquad\pi^{-1}(S^{2}_{(N)})=S^{3}_{(N)}=\{S^{3}:v\neq0\},\nn \\
&S^{2}_{(S)}=\{S^{2}:b_{z}\neq -1\}\qquad\Rightarrow\qquad\pi^{-1}(S^{2}_{(S)})=S^{3}_{(S)}=\{S^{3}:u\neq0\},
\label{triv} 
\end{align}
with $S^{3}_{(j)}\simeq S^{2}_{(j)}\times \U(1)$ via the diffeomorphisms:
\begin{align}
g\simeq(u,v)\in\,S^{3}_{(N)}\quad:\lambda_{N}(g)=(\pi(g);\frac{v}{\abs{v}})\in\,S^{2}_{(N)}\times \U(1), \nn \\
g\simeq(u,v)\in\,S^{3}_{(S)}\quad:\lambda_{S}(g)=(\pi(g);\frac{u}{\abs{u}})\in\,S^{2}_{(S)}\times \U(1). \nn
\end{align}
The set of transition functions associated with this trivialisation is given by 
$\lambda^{-1}_{NS}=\lambda_{SN}=\lambda_{S}\circ\lambda_{N}^{-1}:(S^{2}_{(N)}\cap S^{2}_{(S)})\times \U(1)\to \U(1)$.
Choose $b\sim(\theta,\varphi)\in\,S^{2}_{(N)}\cap S^{2}_{(S)}$. The element $(b,e^{i\alpha})\in\,(S^{2}_{(N)}\cap S^{2}_{(S)})\times \U(1)$ is mapped into
\begin{align}
&\lambda_{N}^{-1}(b,e^{i\alpha})=(u=\frac{b_{y}-ib_{x}}{\sqrt{2(1-b_{z})}}\,e^{i\alpha}; v=\sqrt{\frac{1-b_{z}}{2}}\,e^{i\alpha})\,\in S^{3}_{(N)}\qquad  \nn \\
&\qquad \Rightarrow \lambda_{S}\circ\lambda_{N}^{-1}(b; e^{i\alpha})=(b, e^{i\varphi}e^{i\alpha}) \nn .
\end{align}
This means that $\lambda_{SN}(b)\cdot e^{i\alpha}=e^{i\varphi}e^{i\alpha}$. The transition functions describe a left action of the $\U(1)$ gauge group on itself, and trivially satisfy the cocycle conditions.

For any integer $n$ there is a representation of the gauge group, 
\beq
\rho_{(n)}:\U(1)\to\IC^{*},\qquad\rho_{(n)}(e^{i\alpha})=e^{in\alpha}
\label{repU}
\eeq  
so that for any $n\in\IZ$ there is a line bundle $\mathcal{E}_{n}=SU(2)\times_{\rho_{(n)}}\IC$ associated to the principal Hopf bundle. Since the representations of the gauge group given in \eqref{repU} are defined on $\IC$, the set $\Omega^{r}(S^{3},\IC)_{\rho_{(n)}}\simeq\Omega^{r}(S^{3})$ of $\rho_{(n)}(\U(1))$-equivariant $r$-forms on the Hopf bundle can be easily described  in terms of the action of the vertical field of the bundle, giving the infinitesimal version of the definition in \eqref{Omrg} (with $r=0,\ldots,3$)
\beq
\Omega^{r}(S^{3})_{\rho_{(n)}}=\{\phi\in\,\Omega^{r}(S^{3}):\,\check{\mathrm{r}}_{k}^{*}(\phi)=\rho_{(n)}^{-1}(k)\phi\,\Leftrightarrow\, L_{z}(\phi)=-\frac{in}{2}\phi\}.
\label{Omr}
\eeq 
The sets $\Omega^{r}(S^{3})_{\rho_{(n)}}$ are  $C^{\infty}(S^{2})$-bimodule. The horizontal $\rho_{(n)}(\U(1))$-equivariant $r$-forms are  given as:
\beq
\mathfrak{L}_{n}^{(r)}=\{\phi\in\,\Omega^{r}(S^{3})_{\rho_{(n)}}:\,i_{L_{z}}(\phi)=0\}
\label{Omrh}
\eeq
for $r>0$: one   obviously has $\mcl_{n}^{(3)}=\emptyset$, while 
\beq
\mathfrak{L}_{n}^{(0)}=\Omega^{0}(S^{3})_{\rho_{(n)}}=\{\phi\in\,C^{\infty}(S^{3}):\check{\mathrm{r}}^{*}_{k}(\phi)=\phi\,\,\Leftrightarrow\,\,L_{z}(\phi)=-(in/2)\phi\}.
\label{cln0}
\eeq
With $\Gamma^{(r)}(S^{2},\mathcal{E}_{n})$ the set of
$\ce_{n}$-valued 
 $r$-forms defined on $S^{2}$,  the isomorphisms  in \eqref{isgec} can be written as  isomorphisms  of $C^{\infty}(S^{2})$-bimodule   
 \beq
 \Gamma^{(r)}(S^{2},\mathcal{E}_{n})\simeq\mcl_{n}^{(r)}.
\label{GrL}
\eeq 
They formalise  the equivalence between $r$-form valued sections  on each line bundle $\ce_{n}$ and $\rho_{(n)}(\U(1))$-equivariant horizontal $r$-forms of the principal Hopf bundle. This equivalence can be described -- as in \cite{clp} -- using the local trivialisation \eqref{triv}. A global, algebraic description of them, naturally conceived for the generalisation to the non commutative setting, is in \cite{gl01}, and it  is based on the Serre-Swan theorem\footnote{The theorem of Serre and Swan \cite{s-s}  constructs a complete equivalence between the category of (smooth) vector bundles  over a (smooth) compact manifold $\mathcal{M}$ and bundle maps, and the category of finite projective modules over the commutative algebra $C(\mathcal{M})$ of (smooth) functions over $\mathcal{M}$ and module morphisms. The space $\Gamma(\mathcal{M},\ce)$ of (smooth) sections of a vector bundle $\pi_{\ce}:\ce\to\mathcal{M}$ over a compact manifold $\mathcal{M}$ is a finite projective module over the commutative algebra $C(\mathcal{M})$ and every finite projective $C(\mathcal{M})$-module can be realised as a module of sections of a vector bundle over $\mathcal{M}$. }. 

Given $n\in\,\IZ$, consider an element $\ket{\Psint}\in\,C^{\infty}(S^{3})^{\mn+1}$ whose components are given by:
\begin{align}
&n\geq0:\qquad\ket{\Psint}_{\mu}=\sqrt{\left(\begin{array}{c} n \\ \mu \end{array}\right)}\bar{v}^{\mu}\bar{u}^{n-\mu}\,\in\,\mcl^{(0)}_{n}, \nn \\
&n\leq0:\qquad\ket{\Psint}_{\mu}=\sqrt{\left(\begin{array}{c} \mn \\ \mu \end{array}\right)}v^{\mn-\mu}u^{\mu}\,\in\,\mcl^{(0)}_{n} 
\label{kecl}
\end{align}
with $\mu=0,\ldots\mn$. Recalling the binomial expansion it is easy to compute that:
\begin{align}
&n\geq0:\qquad\hs{\Psint}{\Psint}=\sum_{\mu=0}^{n}\left(\begin{array}{c} n \\ \mu \end{array}\right) u^{n-\mu}v^{\mu}\bar{v}^{\mu}\bar{u}^{n-\mu}=(\bar{u}u+\bar{v}v)^{n}=1, \nn \\
&n\leq0:\qquad\hs{\Psint}{\Psint}=\sum_{\mu=0}^{\mn}\left(\begin{array}{c} \mn \\ \mu \end{array}\right) \bar{u}^{\mu}\bar{v}^{\mn-\mu}v^{\mn-\mu}u^{\mu}=(\bar{u}u+\bar{v}v)^{n}=1.
\label{norp}
\end{align}
The ket-bra element $\qpt=\ket{\Psint}\bra{\Psint}\in\,\mathbb{M}^{\mn+1}(C^{\infty}(S^{2}))$ is then a projector in the free finitely generated module $C^{\infty}(S^{2})^{\mn+1}$, as it satisfies the identities $(\qpt)^{\dagger}=\qpt$, $(\qpt)^{2}=\qpt$. The matrix elements of the projectors are given by $\qpt_{\mu\nu}=\ket{\Psint}_{\mu}\bra{\Psint}_{\nu}$: each projector  $\qpt$ has rank 1, because its trace is the constant unit function given by  
\beq
tr\,\qpt=\sum_{\mu=0}^{\mn}\ket{\Psint}_{\mu}\bra{\Psint}_{\mu}=1.
\label{trcl}
\eeq

Consider the set of $\rho_{(n)}(\U(1))$-equivariant map $\mcl^{(0)}_{n}$ as a left module over $C^{\infty}(S^{2})\subset C^{\infty}(S^{3})$: any equivariant map $\phi\in\,\mcl^{(0)}_{n}$ can be  written in terms of an element $\bra{f}\in\,C^{\infty}(S^{2})^{\mn+1}$ as 
$$
\phi_{f}=\hs{f}{\Psint}=\sum_{\mu=0}^{\mn}\bra{f}_{\mu}\ket{\Psint}_{\mu}. 
$$
Given the set $\Gamma^{(0)}(S^{2},\ce_{n})$ of sections of each associated line bundle $\ce_{n}$, the equivalence with the set $\mcl^{(0)}_{n}$ of $\rho_{(n)}(U(1))$-equivariant maps of the Hopf bundle is formalised via an isomorphism between $C^{\infty}(S^{2})$-left modules, represented by:
\begin{align}
\Gamma^{(0)}(S^{2},\ce_{n})\qquad&\leftrightarrow\qquad\mcl^{(0)}_{n} \nn \\
\bra{\sigma_{f}}=\bra{f}\qpt\qquad&\leftrightarrow\qquad\hs{f}{\Psint} \nn \\
\bra{\sigma_{f}}=\phi_{f}\bra{\Psint}\qquad&\leftrightarrow\qquad\phi_{f}=\hs{\sigma_{f}}{\Psint}
\label{isoc0}
\end{align}
for any $\bra{f}\in\,C^{\infty}(S^{2})^{\mn+1}$. Since from this definition it is $\bra{\sigma_{f}}\qpt=\bra{\sigma_{f}}$, this isomorphism enables to recover $\bra{\sigma_{f}}\in\,\Gamma^{(0)}(S^{2},\ce_{n})~\simeq~C^{\infty}(S^{2})^{\mn+1}\qpt$.
An explicit computation from \eqref{lcf} and \eqref{linv} gives:
\begin{align}
&L_{z}(\tilde{\omega}_{+})=i\tilde{\omega}_{+} \qquad\Rightarrow\qquad \tilde{\omega}_{+}\in\,\mcl_{-2}^{(1)};\nn \\
&L_{z}(\tilde{\omega}_{-})=-i\tilde{\omega}_{-} \qquad\Rightarrow\qquad \tilde{\omega}_{-}\in\,\mcl_{2}^{(1)},
\label{Lzom}
\end{align}
so that for any $n\in\,\IZ$  the set of $\rho_{(n)}(\U(1))$-equivariant horizontal 1-forms of the Hopf bundle is 
\beq
\mcl_{n}^{(1)}=\{\phi=\phi^{\prime}\tilde{\omega}_{-}+\phi^{\prime\prime}\tilde{\omega}_{+}:\,\phi^{\prime}\in\,\mcl^{(0)}_{n-2}\,\mathrm{and}\,\phi^{\prime\prime}\in\,\mcl_{n+2}^{(0)}\}.
\label{cln1}
\eeq
For $n=0$ one also recovers from \eqref{isotr} the equivalence $\mcl_{0}^{(1)}\simeq\Omega^{1}(S^{2})$, so to have the $C^{\infty}(S^{2})$-bimodule identification $\mcl_{n}^{(1)}\simeq\Omega^{1}(S^{2})\otimes_{C^{\infty}(S^{2})}\mcl_{n}^{(0)}$.   For $r=1$ the isomorphism in \eqref{GrL} can be written as:
\begin{align}
\Gamma^{(1)}(S^{2},\ce_{n})\simeq\Omega^{1}(S^{2})^{\mn+1}\cdot\qpt&\qquad\leftrightarrow\qquad\mcl_{n}^{(1)}\simeq\Omega^{1}(S^{2})\otimes_{C^{\infty}(S^{2})}\mcl_{n}^{(0)},\nn \\
\bra{\sigma}=\phi\bra{\Psint}&\qquad\leftrightarrow\qquad\phi=\hs{\sigma}{\Psint}.
\label{isoc1}
\end{align}
Given any $\phi\in\,\mcl_{n}^{(1)}$, set $\bra{\sigma}=\phi\bra{\Psint}\in\,\Omega^{1}(S^{2})^{\mn+1}$, so to have $\bra{\sigma}=\bra{\sigma}\,\qpt$. To write the inverse mapping, consider $\bra{\sigma}\in\,\Omega^{1}(S^{2})^{\mn+1}\qpt$ with components  $\bra{\sigma}_{\mu}\in\,\Omega^{1}(S^{2})$ in the bra-vector notation, satisfying $\bra{\sigma}_{\mu}\qpt_{\mu\nu}=\bra{\sigma}_{\nu}$. Define $\phi=\hs{\sigma}{\Psint}$: it is then straightforward to recover that $\phi\in\,\mcl_{n}^{(1)}$ and that $\bra{\sigma}_{\mu}=\phi\,\bra{\Psint}_{\mu}$.

The same path can be followed to analyse the higher order forms. One has $L_{z}(\tilde{\omega}_{-}\wedge\tilde{\omega}_{+})=0$, so the $C^{\infty}(S^{2})$-bimodule of horizontal $\rho_{(n)}(\U(1))$-equivariant 2-forms of the Hopf bundle is given by 
\beq
\mcl_{n}^{(2)}=\{\phi=\phi^{\prime\prime\prime}\tilde{\omega}_{-}\wedge\tilde{\omega}_{+}:\,\phi^{\prime\prime\prime}\in\,\mcl_{n}^{(0)}\}\simeq\Omega^{2}(S^{2})\otimes_{C^{\infty}(S^{2})}\mcl_{n}^{(0)}
\label{cln2}
\eeq
 for any $n\in\,\IZ$. It is clear that 
for $r=2$ the isomorphism in \eqref{GrL} can be written as:
\begin{align}
\Gamma^{(2)}(S^{2},\ce_{n})\simeq\Omega^{2}(S^{2})^{\mn+1}\cdot\qpt&\qquad\leftrightarrow\qquad\mcl_{n}^{(2)}\simeq\Omega^{2}(S^{2})\otimes_{C^{\infty}(S^{2})}\mcl_{n}^{(0)},\nn \\
\bra{\sigma}=\phi\bra{\Psint}&\qquad\leftrightarrow\qquad\phi=\hs{\sigma}{\Psint}.
\label{isoc1}
\end{align}

\bigskip

The most natural choice of a connection, compatible with the local trivialisation, is given via the definition, as a $\IC$-valued connection 1-form, of 
\beq
\omega=\frac{i}{2}\tilde{\omega}_{z}=(u^{*}\dd u+v^{*}\dd v).
\label{cof}
\eeq 
It globally -- i.e. trivialisation independent --  selects the horizontal part of the tangent space as the left $C^{\infty}(S^{3})$-module $H^{(\omega)}(S^{3})\subset\mathfrak{X}(S^{3})=\{L_{\pm}\}$ since
$\omega(L_{\pm})=0$. On the basis of left invariant vector fields the horizontal projection acts as $L_{\pm}^{(\omega)}=L_{\pm}$, $L_{z}^{(\omega)}=0$.

\subsection{A Laplacian operator on the base manifold $S^{2}$}\label{se:Ls2}

 The canonical isomorphism expressed in  \eqref{isotr} allows to recover the exterior algebra $\Omega(S^{2})$ on the basis of the Hopf bundle as the set of horizontal forms in $\Omega(S^{3})$ which are also invariant for the right principal action of the gauge group $U(1)$. Recalling the definition of the $C^{\infty}(S^{2})$-bimodules of $\rho_{(n)}(U(1))$-equivariant forms given in \eqref{cln1} and \eqref{cln2}, it is possible to identify
 \begin{align}
&\Omega^{0}(S^{2})=C^{\infty}(S^{2})\simeq\mcl_{0}^{(0)}; \nn \\
&\Omega^{1}(S^{2})\simeq\mcl_{0}^{(1)}=\{\phi=\phi^{\prime}\tilde{\omega}_{-}+\phi^{\prime\prime}\tilde{\omega}_{+}:\,\phi^{\prime}\in\,\mcl_{-2}^{(0)}, \,\phi^{\prime\prime}\in\,\mcl_{2}^{(0)}\}; \nn \\
&\Omega^{2}(S^{2})\simeq\mcl_{0}^{(2)}=\{f\tilde{\omega}_{-}\wedge\tilde{\omega}_{+}:\,f\in\,\mcl_{0}^{(0)}=C^{\infty}(S^{2})\},
 \label{eals2}
\end{align}  
where all such identifications are $C^{\infty}(S^{2})$-bimodule isomorphisms. 

On the basis manifold $S^{2}\simeq SU(2)/U(1)=\pi(SU(2))$, whose trivialisation is given in \eqref{triv}, consider the metric 
\beq
\check{g}=2\alpha\,(\tilde{\omega}_{-}\otimes\tilde{\omega}_{+}+\tilde{\omega}_{+}\otimes\tilde{\omega}_{-})
\label{gmet2}
\eeq
and its associated volume $\check{\theta}=\alpha\,\tilde{\omega}_{x}\wedge\tilde{\omega}_{y}=2i\alpha\,\tilde{\omega}_{-}\wedge\tilde{\omega}_{+}=i_{L_{z}}\theta$ in terms of the volume on the group manifold $S^{3}$. 
The corresponding Hodge duality is the $C^{\infty}(S^{2})$-linear map $\star:\Omega^{k}(S^{2})\to\Omega^{2-k}(S^{2})$ given by:
\beq
\begin{array}{lcl}
\star(\check{\theta})=1, & \qquad\qquad & \star(1)=\check{\theta}, \\
\star(\phi^{\prime\prime}\tilde{\omega}_{+})=i\phi^{\prime\prime}\tilde{\omega}_{+}, & \qquad\qquad &  
\star(\phi^{\prime}\tilde{\omega}_{-})=-i\phi^{\prime}\tilde{\omega}_{-},
\end{array}
\label{cHs2}
\eeq
with  $\phi^{\prime}\in\,\mcl_{-2}^{(0)}$ and $\phi^{\prime\prime}\in\,\mcl_{2}^{(0)}$.
The Laplacian operator on  $S^{2}$ can be now evaluated:
\beq
\Box_{S^{2}}f=\star\dd\star\dd f=\frac{1}{2\alpha}(L_{+}L_{-}+L_{-}L_{+})f.
\label{las2}
\eeq
It corresponds to the action of the Laplacian  $\Box_{S^{3}}$ \eqref{clLa} on the subalgebra algebra  $C^{\infty}(S^{2})\subset C^{\infty}(S^{3})$.

\begin{rema}
\label{bis2cl}
Given the Hodge duality \eqref{cHs2}, the expression \eqref{bifocl} defines a bilinear symmetric tensor $\hs{~}{~}_{S^{2}}:\Omega^{k}(S^{2})\times\Omega^{k}(S^{2})\to C^{\infty}(S^{2})$ (with $k=0,1,2$):
\beq
\hs{\xi}{\xi^{\prime}}_{S^{2}}\check{\theta}=\xi\wedge(\star\xi^{\prime}),
\label{bifocl2}
\eeq
for any  $\xi,\xi^{\prime}\in\,\Omega^{k}(S^{2})$. Its non zero terms are given by:
\begin{align}
&\hs{1}{1}_{S^{2}}=1; \nn \\
&\hs{\phi^{\prime}\tilde{\omega}_{-}}{\phi^{\prime\prime}\tilde{\omega}_{+}}_{S^{2}}=
\hs{\phi^{\prime\prime}\tilde{\omega}_{+}}{\phi^{\prime}\tilde{\omega}_{-}}_{S^{2}}
=\phi^{\prime}\phi^{\prime\prime}/2\alpha; \nn \\
&\hs{\check{\theta}}{\check{\theta}}_{S^{2}}=1:
\label{bs2}
\end{align}
such a tensor coincides with  the restriction to the exterior algebra $\Omega(S^{2})$ of the analogue tensor $\hs{~}{~}_{S^{3}}$.  The expression 
\beq
\hs{\xi}{\xi^{\prime}}^{\sim}_{S^{2}}\check{\theta}=\xi^{\prime*}\wedge(\star\xi),
\label{bifocl2h}
\eeq
with again $\xi,\xi^{\prime}\in\,\Omega^{k}(S^{2})$, defines a bilinear map on $\Omega(S^{2})$, which coincides with the 
 restriction of the bilinear map $\hs{~}{~}^{\sim}_{S^{3}}$ to $\Omega(S^{2})$:
\begin{align}
&\hs{1}{1}_{S^{2}}^{\sim}=1; \nn \\
&\hs{\phi^{\prime}\tilde{\omega}_{-}}{\psi^{\prime}\tilde{\omega}_{-}}_{S^{2}}^{\sim}=\frac{1}{2\alpha}\psi^{\prime*}\phi^{\prime}=
\hs{\phi^{\prime}\tilde{\omega}_{-}}{\psi^{\prime}\tilde{\omega}_{-}}_{S^{3}}^{\sim}, \nn \\
&\hs{\phi^{\prime\prime}\tilde{\omega}_{+}}{\psi^{\prime\prime}\tilde{\omega}_{+}}_{S^{2}}^{\sim}=\frac{1}{2\alpha}\psi^{\prime*}\phi^{\prime}=
\hs{\phi^{\prime\prime}\tilde{\omega}_{+}}{\psi^{\prime\prime}\tilde{\omega}_{+}}_{S^{3}}^{\sim}; \nn \\
&\hs{\check{\theta}}{\check{\theta}}_{S^{2}}^{\sim}=1=\hs{2i\alpha\,\tilde{\omega}_{-}\wedge\tilde{\omega}_{+}}{2i\alpha\,\tilde{\omega}_{-}\wedge\tilde{\omega}_{+}}_{S^{3}}^{\sim}
\label{co12bf}
\end{align}
for any $\phi^{\prime},\psi^{\prime}\in\,\mcl_{-2}^{(0)}$ and $\phi^{\prime\prime},\psi^{\prime\prime}\in\,\mcl_{2}^{(0)}$.
\end{rema}

\begin{rema}
\label{omps}
Introducing from the volume form $\check{\theta}$ an integral  $\int_{\check{\theta}}:\Omega^{2}(S^{2})\mapsto\IC$ with the normalisation $\int_{\check{\theta}}\theta=1$, the bilinear maps in \eqref{bifocl2} and \eqref{bifocl2h} give on the exterior algebra $\Omega(S^{2})$ a symmetric scalar product
and a hermitian inner product, setting:
\begin{align}
&(\xi;\xi^{\prime})_{S^{2}}=\int_{\check{\theta}}\xi\wedge(\star\xi^{\prime}), \label{pss2cl}
\\
&(\xi;\xi^{\prime})_{S^{2}}^{\sim}=\int_{\check{\theta}}\xi^{\prime*}\wedge(\star\xi). \label{phs2cl}
\end{align}
It is clear that they coincide  with the restrictions to $\Omega(S^{2})$ of respectively \eqref{symps} and \eqref{bifoprip}.

\end{rema}

\section{The quantum principal Hopf bundle}\label{se:qphb}

This section  describes a quantum formulation of  a  Hopf bundle. It starts with a description of the algebraic approach  to the theory of differential calculi on Hopf algebras coming from  \cite{wor89,KS97}  and then algebraically presents  the geometric structures of a principal bundle.

\subsection{Algebraic approach to the theory of differential calculi on Hopf algebras}

The first order differential forms on the smooth  group manifold $SU(2)\simeq S^{3}$ have been presented as elements in the space $\mathfrak{X}^{*}(S^{3})$, or more properly as sections of the cotangent bundle $T^{*}(S^{3})$. The set $\Omega^{1}(S^{3})\simeq\mathfrak{X}^{*}(S^{3})$  of 1-forms is a bimodule over $C^{\infty}(S^{3})$, with the exterior derivative $\dd$ satisfying the basic Leibniz rule $\dd(ff^{\prime})=(\dd f)f^{\prime}+f\dd f^{\prime}$ for any $f,f^{\prime}\in C^{\infty}(S^{3})$. Moreover, being $S^{3}$ a compact manifold, any differential form $\theta\in\Omega^{1}(S^{3})$ is necessarily of the form $\theta=f_{k}\dd f^{\prime}_{k}$ (with $k\in\IN$).  

In an algebraic setting, these properties are a definition.  Given a $\IC$-algebra with a unit $\ca$ and $\Omega$ a bimodule over $\ca$ with a linear map $\dd:\ca\to\Omega$,  $(\Omega, \dd)$ is defined a first order differential calculus over $\ca$ if  $\dd(ff^{\prime})=(\dd f)f^{\prime}+f\dd f^{\prime}$ for any $f,f^{\prime}\in\ca$ and if any element $\theta\in\Omega$ can be written as $\theta=\sum_{k}\,f_{k}\dd f^{\prime}_{k}$ with $f_{k},f^{\prime}_{k}\in\ca$.

For a $\IC$-algebra with unit $\ca$, any first order
differential calculus $(\oca{1}, \dd)$ on $\ca$ can be obtained from
the universal calculus $(\oca{1}_{un}, \delta)$. The space of universal 1-forms
is the submodule of $\ca \otimes \ca$ given by
$\oca{1}_{un} = \ker(m:\ca \otimes \ca \to \ca)$, with $m(a\otimes
b)=ab$ the multiplication map. The universal differential 
$\delta: \ca \to \oca{1}_{un}$ is $\delta a = 1 \otimes a - a \otimes 1$.
If $\cn$ is any sub-bimodule of $\oca{1}_{un}$ with
projection $\pi_{\cn}: \oca{1}_{un} \to \oca{1} = \oca{1}_{un} /\cn$, then $(\oca{1}, \dd)$, with $\dd:=\pi_{\cn}\circ\delta$, is a first order
differential calculus over $\ca$ and any such a calculus can be
obtained in this way. The projection $\pi_{\cn}:\Omega^{1}(\ca)_{un}\to\Omega^{1}(\ca)$ is $\pi_{\cn}(\sum_{i}a_i\otimes b_{i})=\sum_{i}a_i\dd b_{i}$ with associated subbimodule $\cn=\ker\pi$.

The concept of action of a group on a manifold is algebraically dualised via the notion of coaction of a Hopf algebra $\ch$ on an algebra $\ca$: if the algebra $\ca$ is covariant for the coaction of a quantum group $\ch=(\ch,\Delta, \eps, S)$, one has a notion of covariant calculi on $\ca$ 
as well, thus translating the idea of invariance of the differential calculus on a manifold for the action of a group. 
Then, let  $\ca$ be a (right, say) $\ch$-comodule algebra, with a right  coaction $\Delta_R : \ca \to \ca \otimes \ch$ which is also an algebra map. In order to state the covariance of the calculus 
$(\oca{1}, \dd)$ one needs to extend the coaction of $\ch$. A map 
$  \Delta_{R}^{(1)} : \oca{1} \to \oca{1}\otimes \ch$ is defined by the requirement
$$
\Delta_{R}^{(1)}(\dd f) = (\dd \otimes \id) \Delta_{R}(f)
$$
and bimodule structure governed by
\begin{align*}
&\Delta_{R}^{(1)}(f \dd f') = \Delta_{R}(f) \Delta_{R}^{(1)}(\dd f') , \\
&  \Delta_{R}^{(1)}((\dd f) f') = \Delta_{R}^{(1)}(\dd f) \Delta_{R}(f').
\end{align*}
The calculus is said to be right covariant if it  happens that 
$$
(\id \otimes \Delta) \Delta_{R}^{(1)} = (\Delta_{R}^{(1)} \otimes \id) \Delta_{R}^{(1)}
$$ 
and 
$$
(\id \otimes \eps) \Delta_{R}^{(1)} = 1 .
$$
A calculus is right covariant if and only if for the corresponding bimodule $\cn$ it is verified that $\Delta_{R}^{(1)}(\cn) \subset \cn \otimes \ch$, where $\Delta_{R}^{(1)}$ is defined on  $\cn$ by formul{\ae} as above with the universal derivation $\delta$ replacing the derivation $\dd$:
\beq
\Delta_{R}^{(1)}(\delta f)=(\delta\otimes \id)\Delta_{R}(f).
\label{cdeu}
\eeq

Differential calculi on a quantum group $\ch=(\ch,\Delta, \eps, S)$ were  studied in \cite{wor89}. As a quantum group  consider a Hopf $*$-algebra with an invertible antipode: the coproduct 
$\Delta: ~\ch ~\to \ch \otimes \ch$  defines  both a right and a left coaction of $\ch$ on itself:
\begin{align}
&\Delta^{(1)}_{R}(\dd h)=(\dd\otimes 1)\Delta(h), \nn \\
&\Delta^{(1)}_{L}(\dd h)=(1\otimes \dd)\Delta(h).
\label{Dl1}
\end{align}
Right and left covariant calculi on $\ch$ will be defined as before. Right covariance of the calculus implies that $\och{1}$ has a module basis $\{\eta_{a}\}$ of right invariant 1-forms, that is 1-forms for which
$$
\Delta_{R}^{(1)}(\eta_{a})= \eta_{a} \otimes 1,
$$
and left covariance of a calculus similarly implies that
$\och{1}$ has a module basis $\{\omega_{a}\}$ of left invariant 1-forms, that is 1-forms for which 
$\Delta_{L}^{(1)}(\omega_{a})=1\otimes\omega_{a}$.
In addition one has the notion of a bicovariant  calculus, namely a both left and right covariant calculus, satisfying the  compatibility condition:
$$
(\id\otimes\Delta^{(1)}_{R})\circ\Delta^{(1)}_{L}=(\Delta^{(1)}_{L}\otimes\id)\circ\Delta^{(1)}_{R}.
$$
Given the bijection 
\beq
r : \ch\otimes\ch\to\ch\otimes\ch,
\qquad r(h\otimes h^{\prime})=(h\otimes1)\Delta(h^{\prime}) , 
\label{3p2}
\eeq
one proves that $ r(\och{1}_{un})=\ch\otimes \ker\eps $. 
Then, if
$\mathcal{Q} \subset \ker\eps$ is a right ideal of $\ker\eps$, the
inverse image $\mathcal{N}_{\mathcal{Q}}=r^{-1}(\ch\otimes\mathcal{Q})$
is a sub-bimodule contained in $\och{1}_{un}$.
The differential calculus defined by such a bimodule,
$\och{1}:=\och{1}_{un}/\mathcal{N}_{\mathcal{Q}}$,
is left-covariant, and any left-covariant differential calculus can
be obtained in this way. Bicovariant calculi are in one to one correspondence with  right ideals $\mathcal{Q}\subset\ker\varepsilon$ which are in addition stable under the right adjoint coaction $\Ad$ of $\ch$ onto itself, that is $\Ad(\cq) \subset \cq \otimes \ch$.  
Explicitly, one has $\Ad=\left(\id \otimes m \right) \left(\tau\otimes
\id \right)\left(S\otimes\Delta\right)\Delta$, with $\tau$ the flip
operator, or $\Ad(h) = \co{h}{2} \otimes \left(S(\co{h}{1}) \co{h}{3} \right)$
using the Sweedler notation $\cop h =: \co{h}{1} \otimes \co{h}{2}$ with summation understood, and higher numbers for iterated coproducts.
%
%

Given the $*$-structure on $\ch$, a first order differential calculus 
$(\Omega^{1}(\ch),\dd)$ on $\ch$ is called a $*$-calculus if there exists an anti-linear involution $*:\Omega^{1}(\ch)\to\Omega^{1}(\ch)$ such that $(h_{1}(\dd h)h_{2})^*=h_{2}^{*}(d(h^*))h_{1}^{*}$ for any $h,h_{1},h_{2}\in\,\ch$. A left covariant first order differential calculus is  \cite{wor89} a $*$-calculus if and only if $(S(Q))^*\in\,\cq$ for any $Q\in\,\cq$. In such a case the $*$-structure is also compatibe with the left coaction $\Delta^{(1)}_{L}$ of $\ch$ on $\Omega^{1}(\ch)$: $\Delta^{(1)}_{L}(dh^{*})=(\Delta^{(1)}(\dd h))^*$.

\bigskip
The ideal $\mathcal{Q}$ also determines the tangent space of
the calculus. This is the complex vector space of elements  $\{X_{a}\}$  in
$\ch^{\prime}$  defined by
$$
\mathcal{X}_{\mathcal{Q}}:=
\{X\in \ch^{\prime} ~:~ X(1)=0,\,X(Q)=0, \,\, \forall \, Q\in\mathcal{Q}\},
$$
whose dimension, which coincides with the dimension of the calculus, is given by $\dim\,\mathcal{X}_{\cq}=\dim(\ker\varepsilon_{\ch}/\cq)$. If the vector space $\mathcal{X}_{\cq}$ is finite dimensional, then \cite{KS97} its elements   $X_{a}$ belong to the dual Hopf algebra $\ch^{o}\subset\ch^{\prime}$. 
Given an infinite dimensional Hopf $*$-algebra $\ch$ and the set $\ch^{\prime}$ of its linear functionals, the set $\ch^{\prime}\otimes\ch^{\prime}$ is a linear subspace of $(\ch\otimes\ch)^{\prime}$ obtained via the identification $X\otimes Y\,\in\,\ch^{\prime}\otimes\ch^{\prime}$ with the linear functional on $\ch\otimes\ch$ determined by $(X\otimes Y)(h_{1}\otimes h_{2})=X(h_{1})Y(h_{2})$. For any $X\in\,\ch^{\prime}$ consider $\Delta X$ as the element in $(\ch\otimes\ch)^{\prime}$ defined by $\Delta X(h_{1}\otimes h_{2})=X(h_{1}h_{2})$. The space $\ch^{o}\subset\ch^{\prime}$ denotes the set of linear functionals  
$X\in\,\ch^{\prime}$ for which $\Delta X\in\,\ch^{\prime}\otimes\ch^{\prime}$, i.e. there exist functionals $\{Y_{a}\},\{Z_{b}\}\in\,\ch^{\prime}$ --  with $a,b=1,\ldots,r$; $r\in\,\IN$ -- such that 
$$
X(h_{1}h_{2})=\sum_{i=1}^{r}Y_{i}(h_{1})Z_{i}(h_{2})\qquad\Leftrightarrow\qquad\Delta X=\sum_{i=1}^{r}Y_{i}\otimes Z_{i}.
$$
Dualising the structure maps from $\ch$ to $\ch^{\prime}$ via:
\begin{align}
&X_{1}X_{2}(h)=X_{1}(h_{(1)})X_{2}(h_{(2)}),\nn \\
&\varepsilon_{\ch^{\prime}}(X)=X(1),\nn \\
&(S_{\ch^{\prime}}(X))(h)=X(S(h)),\nn \\
&\idop_{\ch^{\prime}}(h)=\varepsilon(h), \nn \\
&X^{*}(h)=\overline{X(S(h)^{*})}
\label{dmH}
\end{align}
for any $X, X_{1},X_{2}\,\in\,\ch^{\prime}$ and $h,h_{1},h_{2}\,\in\,\ch$, the dual $\ch^{o}$ is proved  to be the largest Hopf $*$-subalgebra contained in $\ch^{\prime}$.  The presence of a $*$-structure on a first order left-covariant differential calculus  can be translated into a condition on the quantum tangent space: $(\Omega^{1}(\ch),\dd)$ is a left-covariant differential calculus if and only if $\mathcal{X}_{\cq}^{*}\subset\mathcal{X}_{\cq}$ with $\ch^{\prime}$ endowed by the complex structure in \eqref{dmH}.

The exterior derivative can be written as:
\beq
\dd h := \sum\nolimits_a
~(X_{a} \triangleright h) ~\omega_{a} , 
\label{ded}
\eeq
in terms of the canonical left and right $\ch^{\prime}$-module algebra structure on $\ch$ given by \cite{wor87}:
\begin{align}
&X\triangleright h:=h_{(1)}(X(h_{(2)})),\nn \\ 
& h\triangleleft X:=X(h_{(1)})h_{(2)}.
\label{deflr}
\end{align}
Left and right actions mutually commute:
$$
(X_{1}\triangleright h)\triangleleft X_{2} = X_{1}\triangleright(h\triangleleft X_{2}),
$$
and the 
$*$-structures are compatible with both actions:
\begin{align}
&X \lt h^* = ((S(X))^* \lt h)^*,\nn \\
&h^* \rt  X = (h \rt  (S(X))^*)^*,
\qquad \forall\,  X \in\,\ch^{o}, \ h \in \ch. \nn
\end{align} 
Given the two Hopf $*$-algebras $\ch=(\ch,\Delta,\varepsilon,S)$ and $\cu=(\cu,\Delta_{\cu},\varepsilon_{\cu},S_{\cu})$, they  can be dually paired. This duality is expressed by the existence of a bilinear map
$\hs{~}{~}:\cu \times \ch
\to \IC$ such that:
\begin{align}
&\hs{\Delta_{\cu}(U)}{h_{1}\otimes h_{2}}=\hs{U}{h_{1}h_{2}},\nn\\
&\hs{U_{1}U_{2}}{h}=\hs{U_{1}\otimes U_{2}}{\Delta(h)},\nn \\
&\hs{U}{1}=\varepsilon_{\cu}(U), \nn \\
&\hs{1}{h}=\varepsilon(h)
\label{extc}
\end{align}
for any $U_{a}\in\cu(\ch)$ and $h_{b}\in\ch$. The pairing is also required to be compatible with $*$-structures:
\begin{align}
&\hs{U^*}{h}=\overline{\hs{U}{(S(h))^*}}, \nn \\ 
&\hs{U}{h^*}=\overline{\hs{(S_{\cu}(U))^*}{h}}.
\end{align}
Such a dual pairing has the property that $\hs{S_{\cu}(U)}{h}=\hs{U}{S(h)}$. A dual pairing can be defined on the generators and then extended to the whole algebras following the relations \eqref{extc}: it is called non degenerate if the condition $\hs{U}{h}=0$ for any $h\in\,\ch$ implies $U=0$, and if $\hs{U}{h}=0$ for any $U\in\,\cu$ implies $h=0$. 

It comes from this analysis out that via a non degenerate dual pairing  between the two Hopf algebras $\ch$ and $\cu$ , it is possible to regard $\cu$ as a Hopf $*$-subalgebra of $\ch^{o}$, and $\ch$ as a Hopf $*$-subalgebra of $\cu^{o}$, after the identifications  $U(h)=h(U)=\hs{U}{h}$ for any $U\in\,\cu$ and $h\in\,\ch$.    
A further comparison among relations \eqref{dmH} and \eqref{extc} shows that  $\ch$ and $\ch^{o}$ are  dually paired in a natural way, with a   pairing which is non degenerate if $\ch^{o}$ separates the points in $\ch$. 

The
derivation nature of elements in $\mathcal{X}_{\mathcal{Q}}$ is
expressed by their coproduct, 
$$
\Delta(X_{a})=1\otimes X_{a}+
\sum\nolimits_b X_{b}\otimes f_{ba},
$$
with the elements  $f_{ab} \in \ch^{o}$ having
specific properties \cite{wor89}:
\begin{align*}
&\Delta(f_{ab})=f_{ac}\otimes f_{cb}, \nn \\
& \varepsilon(f_{ab})=\delta_{ab}, \nn \\
&S(f_{ab})f_{bc}=f_{ab}S(f_{bc})=\delta_{ac}.
\end{align*}
 These elements also control the commutation
relation between the basis 1-forms and elements of $\ch$:
\begin{align}
&\omega_{a} h = \sum\nolimits_b (f_{ab}\triangleright h)\omega_{b}, \nn \\ 
&h \omega_{a} = \sum\nolimits_b \omega_{b} \left( (S^{-1}(f_{ab}) )\triangleright h \right)  
\qquad \mathrm{for} \quad h \in \ch. \nn
\end{align}
For a left covariant differential calculus, the elements $X_{a}\in\mathcal{X}_{\cq}$ play the role which is classically played by the vectors tangent to  a Lie group manifold at the group identity: the first of equations \eqref{deflr} transforms them into the analogue of left invariant derivations on the Hopf algebra of functions on the group. Their dual forms $\omega_{a}$ play the role of the left invariant one forms. For a bicovariant differential calculus it is possible to define a basis of the bimodule of 1-forms which are right invariant. The right coaction of $\ch$ on  $\Omega^{1}(\ch)$ defines a matrix:
\beq
\Delta^{(1)}_{R}(\omega_{a})=\omega_{b}\otimes J_{ba}
\eeq
where $J_{ab}\in\ch$. This matrix is invertible, since $S(J_{ab})J_{bc}=\delta_{ac}$ and $J_{ab}S(J_{bc})=\delta_{ac}$; it satisfies the properties  $\Delta(J_{ab})=J_{ac}\otimes J_{cb}$, $\varepsilon(J_{ab})=\delta_{ab}$  and can be used to define a set of 1-forms:
\beq
\eta_{a}=\omega_{b}S(J_{ba})\qquad\Leftrightarrow\qquad\eta_{a}J_{ab}=\omega_{b}
\eeq
which are  right invariant:
\beq
\Delta^{(1)}_{R}(\eta_{a})=\eta_{a}\otimes 1.
\eeq
On the basis of right invariant  1-forms, the exterior derivative operator acquires the form:
\beq
\dd h=\eta_{a}(h\triangleleft Y_{a}) 
\label{dedr}
\eeq
where $Y_{a}=-S^{-1}(X_{a})$ are the analogue of the derivations associated to  right invariant vector fields. Equation \eqref{decla} is then represented, in an algebraic approach to the theory of differential calculi, by \eqref{ded} and \eqref{dedr}. The derivation nature of $Y_{a}$ as well as the commutation relation between the basis of right invariant 1-forms and elements of $\ch$ are ruled by the same elements $f_{ab}\in\cu(\ch)$ \cite{asca}:
\begin{align}
&\Delta(Y_{a})=Y_{a}\otimes 1+
\sum\nolimits_b S^{-1}(f_{ba})\otimes Y_{b}
\nonumber \\
&\eta_{a}h=(h\triangleleft S^{-2}(f_{ab}))\eta_{b}, \nn \\ 
& h\eta_{a}=\eta_{b}(h\triangleleft(S^{-1}(f_{ab})).
\nonumber
\end{align}

\subsection{Quantum principal bundles}\label{QPB}

An algebraic formalisation of the geometric structures of a principal bundle has been introduced in \cite{bm93} and refined in \cite{BM97}. A slightly different formulation of such a structure is in \cite{durI,durII}; an interesting comparison between the two approaches is in \cite{durcomm}.

Following \cite{bm93}, consider as a  total space  an algebra $\mathcal{P}$ (with multiplication $m : \mathcal{P}\otimes \mathcal{P} \to \mathcal{P}$) and as structure group  
a Hopf algebra $\ch$. Thus $\mathcal{P}$ is a right $\ch$-comodule algebra
with coaction $\Delta_{R}\,:\,\mathcal{P}\to\mathcal{P}\otimes\ch$.
The subalgebra of the right coinvariant elements, 
$
\mathcal{B}=\mathcal{P}^{\ch}\,=\,\{p\in\mathcal{P}\,:\,\Delta_{R} p
= p\otimes 1\}
$,
is the base space of the bundle.
At the `topological level' the principality of the bundle is the requirement of 
exactness of the sequence:
\beq
0\, \to \,\mathcal{P}\left(\Omega^{1}(\mathcal{B})_{un}\right)\mathcal{P}\,
 \to \,\Omega^{1}(\mathcal{P})_{un}
\,\stackrel{\chi}  \to \,\mathcal{P}\otimes \ker\varepsilon_{\ch}\,  \to \,0
\label{topes}
\eeq
with $\Omega^{1}(\mathcal{P})_{un}$ and $\Omega^{1}(\mathcal{B})_{un}$ the universal calculi and the map $\chi$  defined by
\beq
\chi:\,\mathcal{P}\otimes\mathcal{P}\, \to\,\mathcal{P}\otimes\ch,\qquad 
\chi:=\left(m \otimes \id\right)\left(\id\otimes\Delta_{R}\right)
\label{chimap},
\eeq
or $\chi(p' \otimes p) = p' \Delta_R(p)$. The exactness of this sequence is equivalent to the requirement that the analogous `canonical map'  $\mathcal{P}\otimes_{\mathcal{B}}\mathcal{P}\to\,
\mathcal{P}\otimes\ch$ (defined as the formula above) is an isomorphism. This is the definition that the inclusion $\mathcal{B} \hookrightarrow  \mathcal{P}$ be a Hopf-Galois extension \cite{sc90}. 

\begin{rema}
The surjectivity of the map $\chi$ appears as the dual translation of the classical condition  that the action of the structure group on the total space of the principal bundle is free. In the classical setting described in section \ref{se:ic}, given the principal bundle $(\mathcal{P},\Ks,[\mathcal{M}],\pi)$, 
 the condition that the right principal  action $\mathrm{r}_{k}$ is free can be written as the injectivity of the map:
$$
P\times G\,\to\,P\times_{M} P,\qquad(p,k)\,\mapsto\,(p,\mathrm{r}_{k}(p)),
$$ 
whose dualisation is  the condition of the surjectivity of the map $\chi$.
\end{rema}
\bigskip
With differential calculi on both the total algebra $\cp$ and the structure Hopf algebra $\ch$ one needs compatibility conditions that eventually lead to an exact sequence like in \eqref{topes} with the  calculi at hand replacing the universal ones. Then, let $\left(\Omega^{1}(\mathcal{P}),\dd\right)$ be 
a $\ch$-covariant differential calculus  on $\mathcal{P}$ given via the subbimodule $\mathcal{N}_{\mathcal{P}} \in \left(\Omega^{1}(\mathcal{P})_{un}\right)$, and $\left(\Omega^{1}(\ch),\dd\right)$ a bicovariant differential calculus on $\ch$ given via the $\Ad$-invariant right ideal $\mathcal{Q}_{\ch} \in \ker\varepsilon_{\ch}$. 
In order to extend the coaction $\Delta_{R}$ of $\ch$ on $\cp$ to a coaction of $\ch$ on $\Omega^{1}(\cp)$, one requires $\Delta_{R}(\cn_{\cp})\subset\cn_{\cp}\otimes \ch$.
The coaction $\Delta_{R}$ of $\ch$ on $\cn_{\cp}\subset\,\cp\otimes\cp$ is understood as a usual coaction of a Hopf algebra on a tensor product of its comodule algebras, i.e.
$$
\Delta_{R}=(\id\otimes\id\otimes\cdot)\circ(\id\otimes \tau\id)\circ(\Delta_{R}\otimes\Delta_{R}).
$$
The condition $\Delta_{R}(\cn_{\cp})\subset\cn_{\cp}\otimes\ch$ is equivalent to the condition \eqref{cdeu}.

The compatibility of the calculi are then the  requirements that 
$\chi(\mathcal{N}_{\mathcal{P}})\subseteq \mathcal{P}\otimes\mathcal{Q}_{\ch}$ and that the map 
$\sim_{ {\mathcal{N}_{{\mathcal{P}} }}} : \Omega^{1}(\mathcal{P}) \to  
\mathcal{P}\otimes(\ker\varepsilon_{\ch}/\mathcal{Q}_{\ch})$, defined by the diagram
\beq
\begin{array}{lcl}
\Omega^{1}(\mathcal{P})_{un} 
& \stackrel{{\pi_{\mathcal{N}}}} {\longrightarrow} & \Omega^{1}(\mathcal{P}) \\
\downarrow \chi  &  & \downarrow\sim_{{ \mathcal{N}_{\mathcal{P}}}} \\
\mathcal{P}\otimes \ker\varepsilon_{\ch} &
\stackrel{{\id\otimes\pi_{\cq_{\ch}}}}{\longrightarrow}&
\mathcal{P}\otimes(\ker\varepsilon_{\ch}/\mathcal{Q}_{\ch})
\end{array}
\label{qdia}
\eeq
(with $\pi_{\cn}$ and $\pi_{\cq_{\ch}}$ the natural projections) is
surjective and has kernel 
\beq 
\ker\sim_{{ 
\mathcal{N}_{\mathcal{P}}}} = \mathcal{P}\Omega^{1}(\mathcal{B})\mathcal{P}
=:\Omega^{1}_{\mathrm{hor}}(\cp). 
\label{codd}
\eeq
Here $\Omega^{1}(\mathcal{B})=\cb \dd \cb$ is the space of nonuniversal 
1-forms on $\mathcal{B}$ associated to the bimodule $\mathcal{N}_{\mathcal{B}}:= \mathcal{N}_{\mathcal{P}} \cap \Omega^{1}(\mathcal{B})_{un}$. 
These conditions ensure the exactness of the sequence:
\beq
0\,\to\,\mathcal{P}\Omega^{1}(\mathcal{B})\mathcal{P}\,\to\,
\Omega_{1}(\mathcal{P})
\,\stackrel{\sim_{\mathcal{N}_{\mathcal{P}}}} \longrightarrow\, \mathcal{P}\otimes \left(\ker\varepsilon_{\ch}/\mathcal{Q}_{\ch}\right)\,\to\,0 .
\label{des}
\eeq
The condition 
$\chi(\mathcal{N}_{\mathcal{P}})\subseteq \mathcal{P}\otimes\mathcal{Q}_{\ch}$ is needed to have a well defined map 
$\sim_{ {\mathcal{N}_{{\mathcal{P}} }}}$: with all conditions for a quantum principal bundle $(\cp, \cb, \ch; \cn_{\cp}, \cq_{\ch})$ satisfied, this inclusion implies the  equality $\chi(\mathcal{N}_{\mathcal{P}})=\mathcal{P}\otimes\mathcal{Q}_{\ch}$.
Moreover, if $(\cp, \cb, \ch)$ is a quantum principal bundle with the universal calculi, the equality $\chi(\mathcal{N}_{\mathcal{P}})=\mathcal{P}\otimes\mathcal{Q}_{\ch}$ ensures that $(\cp, \cb, \ch; \cn_{\cp}, \cq_{\ch})$ is a quantum principal bundle with the corresponding nonuniversal calculi.


Elements in the quantum tangent space  $\mathcal{X}_{\cq_{\ch}}(\ch)$ giving the calculus on the structure quantum group $\ch$ act on $\ker\varepsilon_{\ch} / \mathcal{Q}_{\ch}$
via the pairing $\hs{\cdot}{\cdot}$ between $\ch^{o}$ and $\ch$. Then, with each $\xi\in\mathcal{X}_{\cq_{\ch}}(\ch)$ one defines a map 
\beq\label{hf}
\tilde{\xi} : \Omega^{1}(\mathcal{P}) \to\mathcal{P}, \qquad \tilde{\xi}:=\left(\id\otimes\xi\right) \circ (\sim_{\mathcal{N}_{\mathcal{P}}} )
\eeq
and declare a 1-form $\omega \in \Omega^{1}(\mathcal{P})$ to be  horizontal if{}f
$\tilde{\xi}\left(\omega\right)=0$,  for all elements $\xi\in
\mathcal{X}_{\cq_{\ch}}(\ch)$.
The collection of horizontal 1-forms is easily seen to coincide with $\Omega^{1}_{\mathrm{hor}}(\cp)$ in \eqref{codd}.

\subsection{A topological quantum Hopf bundle}

As a  step toward a quantum formulation of the classical  Hopf bundle $\pi:S^{3}\to S^{2}$ this section will describe, following \cite{lrz}, a topological $\U(1)$-bundle \cite{bm93} over the standard Podle\'s sphere $\sq$ \cite{Po87},  with  total space  the manifold of the quantum group $SU_{q}(2)$.

\subsubsection{The algebras}\label{qdct}

The coordinate algebra $\ASU$ of the quantum group
$\SU$ is the $*$-algebra generated by $a$ and~$c$, with relations
\begin{align}
\label{derel}
& ac=qca\,\,\,\,\,ac^*=qc^*a\,\,\,\,\,cc^*=c^*c , \nn \\
& a^*a+c^*c=aa^*+q^{2}cc^*=1 .
\end{align}
The deformation parameter $q\in\IR$ is taken in the interval 
$0<q<1$, since for $q>1$ one gets isomorphic algebras; at $q=1$ one recovers the commutative coordinate algebra on the group manifold $\mathrm{SU(2)}$.
The Hopf algebra structure for $\ASU$ is given by the coproduct:
$$
\Delta\,\left[
\begin{array}{cc} a & -qc^* \\ c & a^*
\end{array}\right]=\left[
\begin{array}{cc} a & -qc^* \\ c & a^*
\end{array}\right]\otimes\left[
\begin{array}{cc} a & -qc^* \\ c & a^*
\end{array}\right] ,
$$
antipode:
$$
S\,\left[
\begin{array}{cc} a & -qc^* \\ c & a^*
\end{array}\right]=\left[
\begin{array}{cc} a^* & c^* \\ -qc & a
\end{array}\right],
$$
and counit:
$$\epsilon \left[
\begin{array}{cc} a & -qc^* \\ c & a^*
\end{array}\right]=\left[
\begin{array}{cc} 1 & 0 \\ 0 & 1
\end{array}\right].
$$

\bigskip
The quantum universal envelopping algebra $\su$ is the Hopf $*$-algebra
generated as an algebra by four elements $K,K^{-1},E,F$ with $K K^{-1}=1$ and  subject to
relations: 
\begin{align} 
&K^{\pm}E=q^{\pm}EK^{\pm}, \nonumber \\  
&K^{\pm}F=q^{\mp}FK^{\pm}, \nonumber \\  
&[E,F] =\frac{K^{2}-K^{-2}}{q-q^{-1}}. 
\label{relsu}
\end{align} 
The $*$-structure is
$$
K^*=K, \qquad E^*=F, \qquad F^*=E, 
$$
and the Hopf algebra structure is provided  by coproduct:
\begin{align*}
&\Delta(K^{\pm}) =K^{\pm}\otimes K^{\pm}, \\
&\Delta(E) =E\otimes K+K^{-1}\otimes E,  \\ 
&\Delta(F)
=F\otimes K+K^{-1}\otimes F;
\end{align*} 
antipode:
$$
S(K) =K^{-1}, \qquad
S(E) =-qE, \qquad
S(F) =-q^{-1}F; 
$$
and a counit:
$$
\varepsilon(K)=1, \qquad \varepsilon(E)=\varepsilon(F)=0.
$$
{}From the
relations \eqref{relsu}, the quadratic quantum Casimir element:
\beq\label{cas}
C_{q}\,:=\,\frac{qK^{2}-2+q^{-1}K^{-2}}{(q-q^{-1})^{2}}+FE-\tfrac{1}{4}
\eeq
generates the centre of $\su$. The irreducible finite dimensional $*$-representations $\sigma_{J}$ of $\su$ (see e.g. \cite{maj95})  are labelled by nonnegative
half-integers $J\in \half \IN$ (the spin); they are given by\footnote{The `$q$-number' is defined as:
\begin{equation}
[x] = [x]_q := \frac{q^x - q^{-x}}{q - q^{-1}} ,
\label{eq:q-integer}
\end{equation} 
for $q \neq 1$ and any $x \in \IR$.}
\begin{align}
\sigma_J(K)\,\ket{J,m} &= q^m \,\ket{J,m},
\nn \\
\sigma_J(E)\,\ket{J,m} &= \sqrt{[J-m][J+m+1]} \,\ket{J,m+1},
\label{eq:uqsu2-repns} \\
\sigma_J(F)\,\ket{J,m} &= \sqrt{[J-m+1][J+m]} \,\ket{J,m-1},
\nn
\end{align}
where the vectors $\ket{J,m}$, for $m = J, J-1,\dots, -J+1, -J$, form
an orthonormal basis for the $(2J+1)$-dimensional, 
irreducible $\su$-module $V_J$, and the brackets
denote the $q$-number.  Moreover, $\sigma_J$
is a $*$-representation of $\su$, with respect to the hermitian
scalar product on $V_J$ for which the vectors $\ket{J,m}$ are
orthonormal. In each representation $V_J$, the Casimir \eqref{cas}
is a multiple of the identity with constant given by:
\beq
C_{q}^{(J)}= [J+\half]^2- \tfrac{1}{4}.
\label{quca}
\eeq

\bigskip

The Hopf algebras  $\su$ and $\ASU$ are dually paired. The bilinear mapping $\hs{\cdot}{\cdot}:\su\times\ASU\mapsto\IC$ compatible with the $*$-structures, is set on the generators by: 
\begin{align}
&\langle K,a\rangle=q^{-1/2}, \quad \langle K^{-1},a\rangle=q^{1/2}, \nn\\
& \langle K,a^*\rangle=q^{1/2}, \quad \langle K^{-1},a^*\rangle=q^{-1/2}, \nn\\
&\langle E,c\rangle=1, \quad \langle F,c^*\rangle=-q^{-1}, \label{ndp}
\end{align}
with all other couples of generators pairing to~0. Since the deformation parameter $q$ runs in the real interval range $]0,1[$, this pairing is proved \cite{KS97} to be non degenerate. The canonical left and right actions of $\su$ on $\ASU$ can be recovered by:
\beq
\label{lact}
\begin{array}{lll} 
K^{\pm}\triangleright a^{s} =q^{\mp\frac{s}{2}}a^{s} & 
F\triangleright a^{s} =0 &
E\triangleright a^{s} =-q^{(3-s)/2} [s] a^{s-1} c^{*} \\
K^{\pm}\triangleright a^{* s} =q^{\pm\frac{s}{2}}a^{* s} & 
F\triangleright
a^{*s} =q^{(1-s)/2} [s] c a^{* s-1} &
E\triangleright a^{* s} =0 \\
K^{\pm}\triangleright c^{s} =q^{\mp\frac{s}{2}}c^{s} &
F\triangleright c^{s} =0 &
E\triangleright c^{s} =q^{(1-s)/2} [s]  c^{s-1} a^* \\
K^{\pm}\triangleright c^{* s} =q^{\pm\frac{s}{2}}c^{* s} &
F\triangleright c^{*s} =-q^{-(1+s)/2} [s] a c^{*s-1}  &
E\triangleright c^{* s} =0;
\end{array}
\eeq
and:
\beq
\label{ract}
\begin{array}{lll}
 a^{s}\triangleleft K^{\pm} =q^{\mp\frac{s}{2}}a^{s} & 
a^{s}\triangleleft F =q^{(s-1)/2} [s] c a^{s-1} &
a^{s}\triangleleft E =0 \\
 a^{* s}\triangleleft K^{\pm} =q^{\pm\frac{s}{2}}a^{* s} & 
a^{* s}\triangleleft F =0 &
a^{*s}\triangleleft E =-q^{(3-s)/2} [s] c^{*}a^{*s-1} \\
c^{s}\triangleleft K^{\pm} =q^{\pm\frac{s}{2}}c^{s}  & 
c^{s}\triangleleft F =0  &
c^{s}\triangleleft E =q^{(s-1)/2} [s] c^{s-1} a \\
c^{* s}\triangleleft K^{\pm} =q^{\mp\frac{s}{2}}c^{* s} &
c^{* s}\triangleleft F =-q^{(s-3)/2} [s] a^{*}c^{*s-1} &
c^{* s}\triangleleft E =0.
\end{array}
\eeq

\bigskip

Denote $\ca(\U(1)):=\IC[z,z^*] \big/ \!\!<zz^* -1>$; the map $\pi: \ASU \, \to\,\ca(\U(1))$,
\beq  \label{qprp}
 \pi\,\left[
\begin{array}{cc} a & -qc^* \\ c & a^*
\end{array}\right]=
\left[
\begin{array}{cc} z & 0 \\ 0 & z^*
\end{array}\right]
\eeq 
is a surjective Hopf $*$-algebra homomorphism, so that $\ca(\U(1))$
becomes a quantum subgroup of $SU_{q}(2)$ with a right coaction,
\beq 
\Delta_{R}:= (\id\otimes\pi) \circ \Delta \, : \, \ASU \,\mapsto\,\ASU \otimes
\ca(\U(1)) . \label{cancoa} 
\eeq 
The coinvariant elements for this
coaction, 
elements $b\in\ASU$ for which $\Delta_{R}(b)=b\otimes 1$,  
form a subalgebra of $\ASU$ which is the coordinate algebra
$\Asq$ of the standard Podle\'s sphere $\sq$. From:
\begin{align}
&\Delta_{R}(a)=a\otimes z, \nn \\
&\Delta_{R}(a^*)=a^*\otimes z^{*}, \nn \\
&\Delta_{R}(c)=c\otimes z,\nn \\
&\Delta_{R}(c^*)=c^*\otimes z^{*} 
\label{DeR} 
\end{align}
as a set of generators
for $\Asq$ one can choose: 
\beq 
\label{podgens}
B_{-} := - ac^* , \qquad 
B_{+} :=q ca^* , \qquad 
B_{0} :=  \frac{q^{2}}{1+q^{2}} - q^{2} cc^*
,\eeq 
satisfying the relations\footnote{I should like to thank T.Brzezinski, who noticed that the commutation relations among the generators $B_{j}$ of the algebra $\Asq$ written in \cite{lrz} are not correct.}:
\begin{align*}
& B_{-}B_{0} = [\frac{q^{2}-q^{4}}{1+q^{2}} B_{-}+q^{2}B_{0}B_{-}], \\
& B_{+}B_{0} =[\frac{ q^{2}-1}{q^{2}+1} B_{+}+q^{-2}B_{0}B_{+}], \\
& B_{+}B_{-}=q \left[q^{-2} B_{0} - (1+q^{2})^{-1} \right] 
\left[q^{-2} B_{0} + (1+q^{-2})^{-1} \right] , \\
& B_{-}B_{+}=q \left[B_{0} + (1+q^{2})^{-1} \right] 
\left[B_{0} - (1+q^{-2})^{-1} \right],
\end{align*}
and $*$-structure:
\[
(B_{0})^*=B_{0}, \qquad (B_{+})^*= - q B_{-} .
\]
The sphere $\sq$ is a
quantum homogeneous space of $\SU$ and the coproduct of $\ASU$
restricts to a left coaction of $\ASU$ on $\Asq$ which on
generators reads:
\begin{align*}
\Delta(B_{-})&=a^{2}\otimes B_{-}-(1+q^{-2})B_{-}\otimes B_{0}+
c^{* 2}\otimes B_{+} , \nn\\
\Delta(B_{0})&= q\, ac\otimes B_{-}+(1+q^{-2})
B_{0}\otimes B_{0}- c^*a^* \otimes B_{+} , \nn\\
\Delta(B_{+})&=q^{2}\, c^{2}\otimes B_{-}+(1+q^{-2})B_{+}\otimes
B_{0}+a^{* 2}\otimes B_{+} .
\end{align*}


\subsubsection{The associated line bundles}\label{se:avb}

The left action of the group-like element $K$ on $\ASU$ allows  \cite{maetal} to give a
vector basis decomposition $\ASU=\oplus_{n\in\IZ} \cl^{(0)}_{n}$, where 
\beq
\label{libu} 
\cl^{(0)}_{n} := \{x \in \ASU ~:~ K \lt x = q^{n/2} x \} .
\eeq 
In particular $\Asq = \cl^{(0)}_{0}$.  One also has $\cl_{n}^{(0)*} \subset \cl^{(0)}_{-n}$ and $\cl^{(0)}_{n}\cl^{(0)}_{m} \subset \cl^{(0)}_{n+m}$. 
Each $\cl^{(0)}_{n}$ is a bimodule over $\Asq$; relations \eqref{DeR} show that they can be equivalently characterised by the coaction $\Delta_{R}$ of the quantum subgroup $\ca(U(1))$ on $\ASU$:
\beq
\cl^{(0)}_{n}=\{x\in\,\ASU~:~\Delta_{R}(x)=x\otimes z^{-n}\}.
\label{clnc}
\eeq
This equation appears as the natural quantum analogue of the classical relation \eqref{cln0}, introducing  $\cl_{n}^{(0)}\subset\ASU$ as $\Asq$-bimodule of \emph{co}-equivariant elements with respect to the coaction \eqref{cancoa} of the gauge group algebra. The relation \eqref{libu} can then be read as an infinitesimal version of that  in \eqref{clnc}. The classical $\mathfrak{L}_{n}^{(0)}$ are recovered as rank 1  projective left $C^{\infty}(S^{2})$-modules: 
the analogue property in the quantum setting  was shown in \cite{SWPod}. 
Each $\cl^{(0)}_{n}$ is isomorphic to a projective left $\Asq$-module of rank 1.
These projective left $\Asq$-modules
give modules of equivariant maps or of sections of line bundles over the quantum sphere $\sq$ with winding numbers (monopole charge) $-n$.
The corresponding projections \cite{BM98,HM98} can be explicitly written. Given $n\in\,\IZ$, consider an element $\ket{\Psin}\in\,\ASU^{\mn+1}$ whose components are:
\begin{align}
&n\geq 0:\qquad\ket{\Psin}_{\mu}=\sqrt{\beta_{n,\mu}} ~ c^{* \mu}a^{* n-\mu}\,\in\,\cl^{(0)}_{n}, \nn\\
&\mathrm{where:}  \qquad \beta_{n,0}=1;\qquad
\beta_{n,\mu}=q^{2\mu}\prod\nolimits_{j=0}^{\mu-1}\left(\frac{1-q^{-2\left(n-j\right)}}{1-q^{-2\left(j+1\right)}}\right),
\quad \mu = 1, \ldots,  n 
\label{ketp}
\end{align}
\begin{align}
&n\leq0:\qquad\ket{\Psin}_{\mu}=\sqrt{\alpha_{n,\mu}}  ~ c^{\mn-\mu}a^{\mu}\,\in\,\cl^{(0)}_{n}, \nn \\
&\mathrm{where:}\qquad\alpha_{n,0}=1;\qquad\alpha_{n,\mu}=\prod\nolimits_{j=0}^{\mn-\mu-1}\left(\frac{1-q^{2\left(\mn-j\right)}}
{1-q^{2\left(j+1\right)}}\right),  \quad \mu = 1, \ldots,  \mn  
\label{ketn}
\end{align}
Using the commutation relations \eqref{derel} and the explicit form of the coefficients in \eqref{ketp} and \eqref{ketn}, it is possible to compute that:
\begin{align}
&n\geq0:\qquad 
\hs{\Psi^{(n)}}{\Psi^{(n)}} = 
\sum\nolimits_{\mu=0}^{n} 
 \beta_{n,\mu}\, a^{n-\mu} c^{\mu} c^{* \mu} a^{* n-\mu} =
 (aa^{*}+q^{2}cc^{*})^{n} = 1, \nn \\
&n\leq0:\qquad  \hs{\Psi^{(n)}}{\Psi^{(n)}} =
\sum\nolimits_{\mu=0}^{\mn} 
 \alpha_{n,\mu}\, a^{* \mu} c^{* \mn-\mu} c^{\mn-\mu} a^{\mu} = (a^{*}a+c^{*}c)^{\mn} = 1
 \label{idpsi}
 \end{align}
so that a projector $\qpp\in\,{\bb M}_{\mn+1}(\ca(S^{2}_{q}))$ can be defined as:
\beq
\qpp=\ket{\Psin}\bra{\Psin}
\label{dP}
\eeq
which is by construction an idempotent - $(\qpp)^{2}=\qpp$ - and selfadjoint operator - $(\qpp)^{\dagger}=\qpp$ - whose entries are:
\begin{align}
&n\geq0:\qquad \qpp_{\mu\nu}=\sqrt{\beta_{n,\mu}\beta_{n,\nu}}\,
c^{*\mu}a^{*n-\mu}a^{n-\nu}c^{\nu} \in \Asq, \nn \\
&n\leq0:\qquad 
\qpn_{\mu\nu}= \sqrt{\alpha_{n,\mu}\alpha_{n,\nu}}\,c^{\mn-\mu}a^{\mu}
a^{*\nu}c^{*\mn-\nu} \in \Asq .
\label{pren}
\end{align}
\bigskip 

The projections  \eqref{dP}  play a central role in the description of the quantum Hopf bundle. As a first application one can prove that the algebra inclusion $\Asq\hookrightarrow\ASU$ satisfies the topological requirements for a quantum principal bundle, when  both the algebras are equipped with the universal calculus.  
\begin{prop}
\label{pro:un}
The datum  $(\ASU,\Asq,\ca(\U(1)))$ is a quantum principal bundle. 
\begin{proof}
The proof consists of showing the exactness of the sequence 
$$
0\,\to\,\ASU\left(\Omega^{1}(\sq)_{un}\right)\ASU\, \to\, \Omega_{1}(\SU)_{un} \, \stackrel{\chi}\longrightarrow \, 
\ASU \otimes \ker\varepsilon_{\U(1)}\,\to\,0
$$
or equivalently that the map 
$\chi : \Omega^{1}(\SU)_{un} \to \ASU\otimes \ker\varepsilon_{\U(1)}$ defined as in (\ref{chimap}) -- and with the $\ca(\U(1))$-coaction on $\ASU$ given in (\ref{cancoa}) -- is surjective. 
Given an element $x\in\,\cl^{(0)}_{n}\subset\ASU$, from  \eqref{clnc} the map $\chi$ acts as:
\beq
\chi(\delta x)=\chi(1\otimes x-x\otimes 1)
=x\otimes(z^{-n}-1).
\label{chex}
\eeq
A generic element in $\ASU\otimes \ker\varepsilon_{\U(1)}$ is of the form $x\otimes\left(z^{n}-1\right)$ with $n\in\IZ$ and $x\in\ASU$. To show surjectivity of $\chi$ the strategy is  to show that $1\otimes\left(z^{n}-1\right)$ is in its image since left $\ASU$-linearity of $\chi$ will give the general result: if $\gamma\in\Omega^{1}(\SU)_{un}$ is such that $\chi(\gamma)=1\otimes\left(z^{n}-1\right)$, then $\chi(x \gamma)=x\left(1\otimes(z^{n}-1)\right)=x \otimes \left(z^{n}-1\right)$. 
Fixed now $n\in\,\IZ$, define an element $\gamma$ in $\ASU$ as $\gamma=\hs{\Psi^{(-n)}}{\delta\Psi^{(-n)}}$ following \eqref{ketp} and \eqref{ketn}. Since $\ket{\Psi^{(-n)}}\in\,\cl^{(0)}_{-n}$, one computes that:
$$
\chi(\gamma)=1\otimes(z^{n}-1),
$$
thus completing the proof.
\end{proof}
\end{prop}

Next, it is possible to identify  the spaces of equivariant maps $\cl^{(0)}_n$ -- or equivalently of \emph{coequivariant} elements $\cl^{(0)}_{n}$ -- with the left
$\Asq$-modules of sections $\ce^{(0)}_{n}=(\Asq)^{\mn+1}\qpp$.  For
this  write any element in the free module $(\Asq)^{\mn
+1}$ as $\bra{f}=(f_0, f_1,\dots, f_{\mn})$ with $f_{\mu} \in \Asq$. This
allows  to write equivariant maps as
\begin{align*}
 \phi_f :=\hs{f}{ \Psi^{(n)} }
&= \sum\nolimits_{\mu=0}^{n} f_{\mu} \sqrt{\beta_{n,\mu}}\,c^{*\mu} a^{*n-\mu}
\qquad \qquad \mathrm{for} \quad n \geq 0 ,  \\
&= \sum\nolimits_{\mu=0}^{\mn} f_{\mu} \sqrt{\alpha_{n,\mu}}\,c^{\mn-\mu} a^{\mu}
\qquad \qquad \mathrm{for} \quad n \leq 0 . 
\end{align*}
making it straightforward to establish the proposition, which generalises to the quantum setting the equivalence \eqref{isoc0}:  
\begin{prop}
\label{isoeqsec}
Given $n\in\,\IZ$, let $\ce^{(0)}_{n}:=(\Asq)^{\mn+1}\qpp$.
There is a left $\Asq$-modules isomorphism:
$$
\cl^{(0)}_{n} ~\xrightarrow{~\simeq~}~ \ce^{(0)}_{n}, \quad \phi_f  \mapsto \bra{\sigma_f} = \phi_f \bra{\Psi^{(n)} }
= \bra{f} \qpp ,
$$ 
with inverse
$$
\ce^{(0)}_{n} ~\xrightarrow{~\simeq~}~ \cl^{(0)}_{n}, \quad \bra{\sigma_f} = \bra{f} \qpp \mapsto \phi_f :=\hs{f}{ \Psi^{(n)} }.
$$
\end{prop}

\subsubsection{A Peter-Weyl decomposition of $\ASU$}\label{se:pw}

The aim of this section is to describe  the known
decomposition of the modules $\cl^{(0)}_{n}$ into representation
spaces under the  action of $\su$ \cite{KS97}. From \eqref{libu} one has a vector space decomposition
$\ASU=\oplus_{n\in\IZ} \cl^{(0)}_{n}$, with  
\beq\label{rellb}
E \lt \cl^{(0)}_{n} \subset \cl^{(0)}_{n+2}, \qquad F \lt \cl^{(0)}_{n} \subset \cl^{(0)}_{n-2} .
\eeq
On the other hand, commutativity of the left and right actions of $\su$ yields that 
$$
\cl^{(0)}_{n} \rt  h \subset \cl^{(0)}_{n}, \qquad \forall \, h\in \su .
$$
It has  already been  shown in \cite{SWPod}
that there is also a decomposition, 
\beq
\label{decoln}
\cl^{(0)}_{n}:=\bigoplus_{J=\tfrac{|n|}{2}, \tfrac{|n|}{2} +1,
\tfrac{|n|}{2} +2, \cdots}V_{J}^{\left(n\right)} ,
\eeq 
with $V_{J}^{\left(n\right)}$ the spin
$J$-re\-pre\-sen\-ta\-tion space (for the right action) of $\su$.
Altogether it gives a Peter-Weyl decomposition for $\ASU$ 
(already given in \cite{wor87}).

More explicitly, the
highest weight vector for each $V_{J}^{\left(n\right)}$ in \eqref{decoln} is
$ c^{J-n/2} a^{*J+n/2}  $:
\begin{align}\label{hwv}
& K \lt (c^{J-n/2} a^{*J+n/2}) = q^{n/2} (c^{J-n/2} a^{*J+n/2}), \nn \\
& (c^{J-n/2} a^{*J+n/2}) \rt  K = q^{J} (c^{J-n/2} a^{*J+n/2}) , \nn \\ 
&(c^{J-n/2} a^{*J+n/2}) \rt  F = 0.
\end{align}
Analogously, the lowest weight vector for each $V_{J}^{\left(n\right)}$ in \eqref{decoln} is $a^{J-n/2} c^{*J+n/2}$:
\begin{align*}
& K \lt (a^{J-n/2} c^{*J+n/2}) = q^{n/2} (a^{J-n/2} c^{*J+n/2}),  \\
& (a^{J-n/2} c^{*J+n/2}) \rt  K = q^{-J} (a^{J-n/2} c^{*J+n/2}) , \\  
&(a^{J-n/2} c^{*J+n/2}) \rt E = 0.
\end{align*}
The elements of the vector spaces $V_{J}^{\left(n\right)}$ can be obtained by acting on the highest weight vectors with the lowering operator 
$\rt  E$, since clearly $\left(c^{J-n/2} a^{*J+n/2}\right)\rt  E\in\cl^{(0)}_{n}$, or explicitly,
$$
K\lt\left[\left(c^{J-n/2} a^{*J+n/2}\right)\rt  E\right]=q^{n/2}
\left[\left(c^{J-n/2} a^{*J+n/2}\right)\rt  E\right] .
$$
To be definite, consider $n\geq0$. The first admissible $J$ is $J=n/2$; the highest weight element is $a^{*n}$ and the vector space 
$V_{n/2}^{\left(n\right)}$ is spanned by 
$\{a^{*n}\rt  E^{l}\}$ with $l=0,\dots,n+1$:
$V_{n/2}^{\left(n\right)} = \mathrm{span} \{ a^{*n},c^{*}a^{*n-1},\ldots , c^{*n} \}$. 
Keeping $n$ fixed, the other admissible values of $J$ are $J=s + n/2$ with $s\in\IN$. The vector spaces $V_{s+n/2}^{\left(n\right)}$ are spanned by $\{c^{s}a^{*s+n}\rt  E^{l}\}$ with $l=0,\ldots,2s+n+1$. 
Analogous considerations are valid when $n\leq0$. In this cases, the admissible values of $J$ are $J=s + \mn/2=s-n/2$, the highest weight vector in $V_{s-n/2}^{\left(n\right)}$ is the element $c^{s-n}a^{*s}$, and a basis is given by the action of the lowering operator $\rt  E$, that is  
$V_{s-n/2}^{\left(n\right)}= \mathrm{span} \{ \left(c^{s-n}a^{*s}\right)\rt  E^{l}, 
\; l=0,\ldots,2s-n+1\}$.


From (\ref{rellb}) one has that the left action $F \lt$ maps $\cl^{(0)}_{n}$ 
to $\cl^{(0)}_{n-2}$. If $p\geq0$, the element $a^{*p}$ is the highest weight vector in $V_{p/2}^{(p)}$ and one has that $F\lt a^{*p} \propto c a^{*p-1}$. The element $c a^{*p-1}$ is the highest weight vector in 
$V_{p/2}^{(p-2)}$ since one finds that $(c a^{*p-1}) \rt F=0$ and 
$(c a^{*p-1}) \rt K=q^{p/2} (c a^{*p-1})$. In the same vein, the elements  
$F^{t}\lt a^{*p} \propto c^{t}a^{*p-t}$ are the highest   
weight elements in $V_{p/2}^{(p-2t)} \subset \cl^{(0)}_{p-2t},\,\,t=0,\ldots,p$. Once again, a complete basis of each subspace $V_{p/2}^{(p-2t)}$ is obtained by the right action of the lowering operator $\rt  E$.

With these considerations, the algebra $\ASU$ can be partitioned into finite dimensional blocks which are the analogues of the Wigner D-functions  \cite{Mos} for the group $SU(2)$. To illustrate the meaning of this partition, start with the element $a^{*}$, the highest weight vector of the space 
$V_{1/2}^{(1)}$. Representing the left action of $F \lt$ with a horizontal arrow and the right action of  $\rt E$ with a vertical one, yields the box
$$
\begin{array}{ccc}
a^{*} & \to & c \\
\downarrow & & \downarrow \\
-qc^{*} & \to & a 
\end{array}  \;\; ,
$$
where the first column is a basis of the subspace 
$V_{1/2}^{(1)}$, while the second column is a basis of the subspace 
$V_{1/2}^{(-1)}$. Starting from $a^{*2}$ -- the highest weight vector of $V_{1}^{(2)}$ -- one gets:
$$
\begin{array}{ccccc}
a^{*2} & \to & q^{-1/2}\left[2\right]c a^{*} & \to & \left[2\right]c^{2} \\
\downarrow & & \downarrow & & \downarrow \\
-q^{1/2}\left[2\right]c^{*}a^{*} & \to & \left[2\right]\left(a a^{*} - c c^{*}\right) & \to & \left[2\right]^2 q^{1/2}ca \\
\downarrow & & \downarrow & & \downarrow \\
q^{2}\left[2\right]c^{*2} & \to & -q^{1/2}\left[2\right]^2 a c^{*} & \to & \left[2\right]^2 a^{2} 
\end{array} \;\; .
$$
The three columns of this box are bases for the subspaces $V_{1}^{(2)}$, $V_{1}^{(0)}$,$V_{1}^{(-2)}$, respectively. The recursive structure is clear. For a positive integer $p$, one has a box $W_{p}$ made up of $\left(p+1\right)\times\left(p+1\right)$ elements. Without explicitly computing the coefficients, one gets:
$$
\begin{array}{ccccccccccc}
a^{*p} & \to & ca^{*p-1} & \to & \ldots & \to & c^{t}a^{*p-t} & \to & \ldots & \to & c^{p} \\
\downarrow & & \downarrow & & \ldots & & \downarrow & & \ldots & & \downarrow \\
c^{*}a^{*p-1} & \to & \ldots & \to & \ldots & \to & \ldots & \to & \ldots & \to & ac^{p-1} \\
\downarrow & & \downarrow & & \ldots & & \downarrow & & \ldots & & \downarrow \\
\ldots & \to & \ldots & \to & \ldots & \to & \ldots & \to & \ldots & \to & \ldots \\
\downarrow & & \downarrow & & \ldots & & \downarrow & & \ldots & & \downarrow \\
c^{*s}a^{*p-s} & \to & \ldots & \to & \ldots & \to & \ldots & \to & \ldots & \to & a^{s}c^{p-s} \\
\downarrow & & \downarrow & & \ldots & & \downarrow & & \ldots & & \downarrow \\
\ldots & \to & \ldots & \to & \ldots & \to & \ldots & \to & \ldots & \to & \ldots \\
\downarrow & & \downarrow & & \ldots & & \downarrow & & \ldots & & \downarrow \\
c^{*p} & \to & ac^{*p-1} & \to & \ldots & \to & a^{t}c^{*p-t} & \to & \ldots & \to & a^{p} \\
\end{array} \;\; .
$$
The space $W_{p}$ is the direct sum of representation spaces for the right action of $\su$, 
$$
W_{p}=\oplus_{t=0}^{p}V_{p/2}^{(p-2t)}, 
$$
and on each $W_{p}$ the quantum Casimir $C_{q}$ acts is the same manner from both the right and the left, with eigenvalue (\ref{quca}), that is 
$C_{q}\lt w_{p} = w_{p}\rt  C_{q} = \left([\frac{p+1}{2}]^{2}-\tfrac{1}{4}\right)w_{p}$, for all $w_{p}\in W_{p}$.
The Peter-Weyl decomposition for the algebra $\ASU$ is given as
$$
\ASU = \oplus_{p\in\IN}W_{p} = 
\oplus_{p\in\IN}\left(\oplus_{t=0}^{p}V_{p/2}^{(p-2t)}\right).
$$
A compatible basis with this decomposition is  given by elements
\beq
w_{p:t,r}:=F^{t}\lt a^{*p}\rt  E^{r}\,\in\,W_{p}
\label{ws}
\eeq
for $t,r=0,1\ldots,p$. In order to get elements in the Podle\'s sphere subalgebra $\Asq\simeq\cl^{(0)}_{0}$ out of a highest weight vector $a^{*p}$ we need $p=2l$ to be even and left action of $F^{l}$: $F^{l}\lt a^{*2l} \propto c^{l}a^{*l} \in \Asq$. Then, the right action of $E$ yields a spherical harmonic decomposition,
\beq
\Asq = \oplus_{l \in\IN}V_{l}^{(0)},
\label{decsp}
\eeq
with a basis of $V_{l}^{(0)}$ given by the vectors $F^{l}\lt a^{*2l}\rt  E^{r}$, for $r=0,1,\dots,2l$.




\subsection{A quantum Hopf bundle with non-universal differential calculi}
Once described how the inclusion $\Asq\hookrightarrow\ASU$ has the structure of a topological quantum principal bundle, the aim of this section is to describe non-universal differential calculi on the algebras $\ASU, \Asq, \ca(\U(1))$, and to show that these are compatible \cite{bm93,BM97}. 

\subsubsection{The left-covariant 3D calculus on $SU_{q}(2)$}\label{se:lcc}

The first differential calculus defined on the quantum group $\SU$
is the left-covariant one developed in \cite{wor87}. It is
three dimensional with corresponding ideal $\mathcal{Q}_{\SU}\subset\ker\varepsilon_{\SU}$
generated by the 6 elements $
\{a^*+q^{2}a-(1+q^{2});c^{2};c^*c;c^{*2};(a-1)c;(a-1)c^*\} $.
Its quantum tangent space  turns out to be, in terms of the non degenerate pairing \eqref{ndp},  the vector space over the complex $\mathcal{X}_{\SU}\subset\cu_{q}(\mathfrak{su(2)})$, whose basis is
\begin{align}
&X_{-}
=q^{-1/2}FK , \nn \\
& X_{+} =q^{1/2}EK ,\nn \\  
&X_{z} =\frac{1-K^{4}}{1-q^{-2}}; 
\label{Xq}
\end{align}
their coproducts result: 
\begin{align}
&\cop X_z = 1\otimes X_z + X_z \otimes K^4, \nn \\
&\cop X_\pm = 1\otimes X_\pm + X_\pm \otimes K^2 .
\end{align}
The differential 
$\dd
: \ASU \to \Omega^1(\SU)$ is 
\beq\label{exts3} 
\dd x=
(X_{+}\triangleright x) \,\omega_{+} + (X_{-}\triangleright x)
\,\omega_{-} + (X_{z}\triangleright x) \,\omega_{z}, 
\eeq 
for all $x\in\ASU$. This equation gives a basis for the dual space of 1-forms $\Omega^{1}(\ASU)$,
\begin{align}
&\omega_{z} =a^*\dd a+c^*\dd c , \nn \\
& \omega_{-} =c^*\dd
a^*-qa^*\dd c^*, \nn \\
&\omega_{+} =a\dd c-qc \dd a , 
\label{q3dom}
\end{align}
of left-covariant forms, that is
$\Delta_{L}^{(1)}(\omega_{s})=1\otimes\omega_{s}$, with $\Delta_{L}^{(1)}$ the (left)
coaction of $\ASU$ onto itself extended to forms \eqref{Dl1}. The above relations \eqref{q3dom} can be inverted to
\begin{align}
&\dd a=-qc^*\omega_{+}+a\omega_{z} , \nn \\
&\dd a^*=-q^{2}a^*\omega_{z}+c\omega_{-},\nn\\
&\dd c=a^*\omega_{+}+c\omega_{z} , \nn \\
&\dd c^*=-q^{2}c^*\omega_{z}-q^{-1}a\omega_{-} .\label{dfro}
\end{align}
A direct computation shows that $(S(\cq_{\SU}))^*\subset\cq_{\SU}$. This differential calculus is then a $*$-calculus, with  $\omega_{-}^*=-\omega_{+}$ and
$\omega_{z}^*=-\omega_{z}$. The bimodule structure is:
\begin{align}\label{bi1}
&\omega_{z}\phi=q^{2n}\phi\omega_{z}, \nn \\
&\omega_{\pm}\phi=q^{n}\phi\omega_{\pm}
\end{align}
for any $\phi\in\,\cl^{(0)}_{n}$.
Higher dimensional forms can be defined
in a natural way by requiring compatibility for commutation
relations and that $\dd^2=0$. Consider the tensor product
$\{\Omega(\SU)\}^{\otimes2}=\Omega^{1}(\SU)\otimes_{\ASU}\Omega^{1}(\SU)$. 
A consistent alternation mapping on $\{\Omega(\SU)\}^{\otimes2}$, generalising the alternation mapping in the classical formalism, can be introduced only if the quantum differential calculus is bicovariant. The strategy to define  a wedge product  comes then from 
Lemma 15 in chapter 14 in  \cite{KS97}, where it is proved that $\mathcal{S}_{\cq_{\SU}}(x)=\sum_{a,b}\hs{X_{a}X_{b}}{x}\omega_{a}\otimes \omega_{b}$ for any $x\in\,\cq_{\SU}$ generates a two-sided ideal in $\{\Omega(\SU)\}^{\otimes2}$.   The bimodule of exterior differential 2-forms results to be the quotient 
\beq
\Omega^{2}(\SU)\simeq\{\Omega^{1}(\SU)\}^{\otimes2}/\ASU\{\mathcal{S}_{\cq}\}\ASU.
\label{2fc}
\eeq
 The wedge product $\wedge:\Omega^{1}(\SU)\times\Omega^{1}(\SU)\to\Omega^{2}(\SU)$ embodies the commutation relations among 1-forms: from the six generators in $\cq_{\SU}$ the elements  generating $\mathcal{S}_{\cq}$ can be written as   
\begin{align}
\label{commc3}
&\omega_{+}\wedge\omega_{+}=\omega_{-}\wedge\omega_{-}=\omega_{z}
\wedge\omega_{z}=0, \nn\\
&\omega_{-}\wedge\omega_{+}+q^{-2}\omega_{+}\wedge\omega_{-}=0,
\nn\\
& \omega_{z}\wedge\omega_{-}+q^{4}\omega_{-}\wedge\omega_{z}=0, \nn \\
& \omega_{z}\wedge\omega_{+}+q^{-4}\omega_{+}\wedge\omega_{z}=0.  
\end{align}
Such commutation rules also show that the bimodule $\Omega^{2}(\SU)$ is 3 dimensional, the three basis 2-forms being exact, since one has
\begin{align} 
\label{dformc3} 
&\dd
\omega_{z} =-\omega_{-}\wedge\omega_{+} , \nn \\ 
& \dd \omega_{+}
=q^{2}(1+q^{2})\omega_{z}\wedge\omega_{+} , \nn \\ 
&\dd \omega_{-}
=-(1+q^{-2})\omega_{z}\wedge\omega_{-};
\end{align} 
the commutation relations moreover clarify that this left covariant calculus has a unique top form $\omega_{-}\wedge\omega_{+}\wedge\omega_{z}$. The $*$-structure is extended to $\Omega^{m+n}$ by $(\alpha\wedge\beta)^*=(-1)^{mn}\beta^{*}\wedge\alpha^{*}$ with $\alpha\in\,\Omega^{m}$ and $\beta\in\Omega^{n}$. This definition is compatible with \eqref{commc3}.

The left covariance of the differential calculus allows to extend  to higher order forms in a natural way the left coaction $\Delta^{(1)}_{L}$ of $\ASU$ on $\Omega^{1}(\SU)$. An element $\eta\in\,\{\Omega^{1}(\SU)\}^{\otimes k}$  can always be written as $\eta=x_{a_{1}~\ldots~ a_{k}}~\omega_{a_{1}}~\otimes~\ldots~\otimes~\omega_{a_{k}}$ in terms of the left invariant forms $\omega_{j}$ in \eqref{q3dom}. Define 
$$
\Delta^{(k)}_{L}(\eta)=x_{a_{1}\ldots a_{k}(1)}\otimes x_{a_{1}\ldots a_{k}(2)}\omega_{a_{1}}\otimes\ldots\otimes\omega_{a_{k}},
$$    
from  the Sweedler notation for the coproduct $\Delta(x_{a_{1}\ldots a_{k}})$. One proves that this definition is consistent on the exterior algebra $\Omega^{k}(\SU)$, as $\Delta^{(2)}_{L}(\mathcal{S}_{\cq})\subset 1\otimes\mathcal{S}_{\cq}$, 
and that 
$\Delta^{(k)}_{L}(\dd\eta)=(1\otimes \dd)\Delta^{(k-1)}_{L}(\eta)$ 
for any $\eta\in\,\Omega^{k}(\SU)$ with $k=1,2,3$. The relations \eqref{dformc3} show then that $\Omega^{2}(\SU)$ has a basis of exact left invariant forms, given by $\dd\omega_{j}$; it is also clear that $\omega_{-}\wedge\omega_{+}\wedge\omega_{z}$ is a left-invariant $3$-form.

\subsubsection{The calculus on the structure group}\label{se:csg}

The strategy adopted in \cite{bm93} consists in defining the calculus on $\U(1)$
via the Hopf projection $\pi$ in \eqref{qprp}. Out of
the $\mathcal{Q}_{\SU}$ which determines the left covariant calculus
on $\SU$, one defines a right ideal
$\mathcal{Q}_{\U(1)}=\pi(\mathcal{Q}_{\SU})\subset\ker\,\varepsilon_{\U(1)}$ for the calculus on
$\U(1)$. 

This specific $\mathcal{Q}_{\U(1)}$ results  generated by the
element $\xi=(z^{-1}-1)+q^{2}(z-1)$, and the differential calculus is then characterised by the quotient $\ker\varepsilon_{\U(1)}/\cq_{\U(1)}$. Any term in $\ker\varepsilon_{\U(1)}$ can be written as $\varphi=u(z-1)=\sum_{j\in\,\IZ}\,u_{j}z^{j}(z-1)$, with $u=\sum_{j\in\,\IZ}\,u_{j}z^{j}\,\in\,\ca(\U(1))$ and $u_{j}\in\,\IC$, so that  the elements $\varphi(j)=z^{j}(z-1)$ define a vector space basis over $\IC$ of $\ker\varepsilon_{\U(1)}$. The basis elements $\varphi(j)$ can be written in terms of the element $\xi$, via the two identities:
\begin{align}
&j\geq0,\qquad\varphi(j)=z^{j}(z-1)=\xi\left(\sum_{m=1}^{j}\,q^{-2m}z^{j-m+1}\right)+q^{-2j}(z-1), \nn\\
&j\leq0, \qquad\varphi(j)=z^{-\mj}(z-1)=-\xi\left(\sum_{m=0}^{\mj-1}\,q^{2m}z^{1+m-\mj}\right)+q^{2\mj}(z-1),
\label{idk}
\end{align}
which can be proved by induction on $j$. Define a map $\lambda:\ker\varepsilon_{\U(1)}\to\ker\varepsilon_{\U(1)}$  setting on the basis elements 
$\lambda(\varphi(j))=q^{-2j}(z-1)$, and linearly extending it to:
\beq
\lambda:\,u(z-1)=\sum_{j\in\,\IZ}\,u_{j}z^{j}(z-1)\qquad\mapsto\qquad\sum_{j\in\,\IZ}\,u_{j}q^{-2j}(z-1).
\label{lam}
\eeq 
It is clear that $\lambda$ describes the choice of a representative element out of the equivalence class $[u(z-1)]\in\,\ker\varepsilon_{\U(1)}/\cq_{\U(1)}$, since it is possible to see that $\ker\lambda=\cq_{\U(1)}$. To prove this assertion, one first directly computes that $\lambda(\xi)=0$, then since $\lambda$ is linear one recovers that $\lambda(u\xi)=\lambda(u(q^{2}(z-1)+(z^{-1}-1)))=q^{2}\lambda(u(z-1))+\lambda(u(z^{-1}-1))$, so to have:
\begin{align}
\lambda(u\xi)&=q^{2}\lambda(u(z-1))+\lambda\left(\sum_{j\in\,\IZ}\,u_{j}z^{j}(z^{-1}-1)\right)\nn\\
&= q^{2}\lambda(u(z-1))+\lambda\left(-\sum_{j\in\,\IZ}\,u_{j}z^{j-1}(z^{-1})\right)\nn \\
&=q^{2}\left(\sum_{j\in\,\IZ}\,u_{j}q^{-2j}(z-1)\right)-\sum_{j\in\,\IZ}\,u_{j}q^{-2(j-1)}(z-1)=0,
\label{uxi}
\end{align}
thus proving that $\cq_{\U(1)}\subset\ker\lambda$.
To prove the inverse inclusion, consider an element $\check{u}=u(z-1)\in\,\ker\varepsilon_{\U(1)}$, and write it as:
\begin{align}
u(z-1)&=\sum_{j\in\,\IZ}\,u_{j}z^{j}(z-1)\nn \\
&=\sum_{j\in\,\IN}\,u_{j}z^{j}(z-1)\,+ \,\sum_{j\in\,\IN}\,u_{-j}z^{-j}(z-1)\nn\\
&=\sum_{j\in\,\IN}\,u_{j}(\alpha(j)\xi+q^{-2j}(z-1))\,+ \,\sum_{j\in\,\IN}\,u_{-j}(\beta(-j)\xi+q^{2j}(z-1))\nn\\
\end{align}
where $\alpha(j)=\sum_{m=1}^{j}\,q^{-2m}z^{j-m+1}$ and $\beta(-j)=\sum_{m=0}^{\mj-1}\,q^{2m}z^{1+m-\mj}$ are the terms proportional to $\xi$ in \eqref{idk} for positive and negative values of $j\in\,\IZ$. The previous  sum can be rewritten as:
\begin{align} &u(z-1)=\xi\left(\sum_{j\in\,\IN}\,u_{j}\alpha(j)\,+ \,\sum_{j\in\,\IN}\,u_{-j}\beta(-j)\right)+\,\sum_{j\in\,\IZ}\,u_{j}q^{-2j}(z-1).\nn
\end{align}
From the definition \eqref{lam}, it is $\lambda(\check{u})=0 \leftrightarrow\sum_{j\in\,\IZ}\,u_{j}q^{-2j}=0$, so the last lines proves that $\ker\lambda\subset\cq_{\U(1)}$. 
\begin{lemm}
\label{pu}
Given  the ideal $\cq_{\U(1)}\subset\ker\varepsilon_{\U(1)}$ generated by the element $\xi=(z^{-1}-1)+q^{2}(z-1)$,  it is $\ker\varepsilon_{\U(1)}/\cq_{\U(1)}\simeq\IC$.
\begin{proof}
Define a map $\tilde{\lambda}:\ker\varepsilon_{\U(1)}\to\IC$  setting, on the basis elements $\varphi(j)\in\,\ker\varepsilon_{\U(1)}$, $\tilde{\lambda}(\varphi(j))=q^{-2j}$ and extending it to $\ker\varepsilon_{\U(1)}$ by linearity. The properties of the map $\lambda$ defined in \eqref{lam} clarify that $\ker\tilde{\lambda}=\cq_{\U(1)}$,  so to give a well defined  map $\tilde{\lambda}:\ker\varepsilon_{\U(1)}/\cq_{\U(1)}\to\IC$. It is immediate to see that  $\tilde{\lambda}$ is an isomorphism of vector spaces, thus describing  the  equivalence: with $w\in\,\IC$,  the  map $\tilde{\lambda}^{-1}(w)=w\in\,[w(z-1)]\,\subset\ker\varepsilon_{\U(1)}$ represents the inverse of the map $\tilde{\lambda}$. 
\end{proof}
\end{lemm}

This result shows that the differential calculus generated by the specific $\cq_{\U(1)}$
is 1D, while a direct computation shows that it is bicovariant. 
As a basis element for its quantum tangent space one can consider 
\beq\label{vvf} 
X=X_{z} =\frac{1-K^{4}}{1-q^{-2}},
\eeq 
with dual left-invariant 1-form given by $\omega_{z}$. This calculus turns out to have a $*$-structure, with $\omega_{z}^*=-\omega_{z}$. Explicitly, one has 
$\omega_{z} = z^{*} \dd z$ with 
\begin{align*}
&\dd z=z\omega_{z} , \\
&\dd z^{*} =-q^{2}z^{*} \omega_{z}; 
\end{align*}
and  noncommutative $\ca(\U(1))$-bimodule  relations
\begin{align}
& z \dd z=q^{2}(\dd z) z; \nn \\
&\omega_{z} z = q^{-2} z\omega_{z}, \nn \\
&\omega_{z} z^{*} = q^{2}z^{*}\omega_{z}. \nn 
\end{align}

\subsubsection{The standard 2D calculus on $\sq$}\label{se:cals2}

The restriction of the above 3D calculus to the sphere $\sq$ yields
the unique left covariant 2-dimensional calculus on the latter \cite{ma05}. An
evolution of this approach has led \cite{SW04} to a description of
the unique 2D calculus of $\sq$ in term of a Dirac operator. The `cotangent bundle' $\Omega^1(\sq)$ is shown to be isomorphic
to the direct sum $\cl^{(0)}_{-2}\oplus\cl^{(0)}_2$, that is the line bundles
with winding number $\pm 2$. Since the element $K$ acts as the
identity on $\Asq$, the differential \eqref{exts3} becomes, when restricted to the latter, 
\begin{align*}
\dd f &= (X_{-}\triangleright f) \,\omega_{-} + (X_{+}\triangleright f) \,\omega_{+}  \\
& =  (F \triangleright f) \,
(q^{-1/2}\omega_{-}) + (E \triangleright f) \,(q^{1/2}\omega_{+}) ,  \qquad \mathrm{for} \, f\in\Asq .
\end{align*}
These leads to break the exterior derivative into a holomorphic and an
anti-holomorphic part, $\dd = \delb + \del$, with:
\begin{align*}
&\delb f=\left(X_{-}\triangleright f\right)\omega_{-} = (F \triangleright f) \,
(q^{-1/2}\omega_{-}) , \nn \\
&\del f=\left(X_{+}\triangleright f\right)\omega_{+} = (E \triangleright f)
\,(q^{1/2}\omega_{+}) , \qquad \mathrm{for} \quad f\in\Asq .
\end{align*}
An explicit computation on the generators \eqref{podgens} of $\sq$
yields: 
$$
\begin{array}{lll}
\delb \bm= q^{-1}\, a^{2} \, \omega_{-}, 
& \delb\bz=q \, ca  \, \omega_{-}, 
& \delb\bp= q\, c^{2} \, \omega_{-} , \\
~\\
\del\bp= q^{2} \,a^{*2} \, \omega_{+}, 
& \del\bz= -q^{2}\, c^{*}a^{*} \, \omega_{+}, 
& \del\bm= q^{2}\, c^{*2} \, \omega_{+} .
\end{array}
$$
The above shows that:
$\Omega^1({\sq})=\Omega^{1}_{-}(\sq) \oplus\Omega^{1}_{+}(\sq)$ where
$\Omega^{1}_{-}(\sq)\simeq \cl^{(0)}_{-2} \simeq \delb(\Asq)$ is the
$\Asq$-bimodule generated by: 
$$
\{\delb\bm,\delb\bz,\delb\bp\}=\{a^{2},ca,c^{2}\}\,\omega_{-} =
q^{2}\omega_{-}\{a^{2},ca,c^{2}\} 
$$
and $\Omega^{1}_{+}(\sq)\simeq
\cl^{(0)}_{+2}\simeq \del(\Asq)$ is the one generated by:
$$
\{\del\bp,\del\bz,\del\bm\}=\{a^{*2},c^{*}a^{*},c^{*2}\}\,
\omega_{+}= q^{-2} \omega_{+} \{a^{*2},c^{*}a^{*},c^{*2}\} .
$$
That these two modules of forms are not free is also expressed by
 the existence of relations among the differential:
$$
\del\bz= q^{-1} \bm\del\bp - q^{3} \bp\del\bm , \qquad 
\delb\bz= q \bp\delb\bm - q^{-3} \bm\delb\bp .
$$

Writing  any 1-form as $\alpha = \phi^{\prime}
\omega_{-} + \phi^{\prime\prime} \omega_{+} \in \cl^{(0)}_{-2} \omega_{-} \oplus \cl^{(0)}_{+2}
\omega_{+}$, the product of 1-forms is
\beq
(\phi^{\prime} \omega_{-} + \phi^{\prime\prime} \omega_{+} ) \wedge (\psi^{\prime} \omega_{-} + \psi^{\prime\prime} \omega_{+} ) =
(q^{-2} \phi^{\prime\prime} \psi^{\prime} - \phi^{\prime}\psi^{\prime\prime} )  \omega_{+}\wedge\omega_{-},
\label{p12c}
\eeq
while the exterior derivative acts as:
\begin{align}
\dd(\phi^{\prime}\omega_{-}+\phi^{\prime\prime}\omega_{+})&=
(\dd\phi^{\prime})\wedge\omega_{-}+\phi^{\prime}\dd\omega_{-}+
(\dd\phi^{\prime\prime})\wedge\omega_{+}+\phi^{\prime\prime}\dd\omega_{+} \nn
\\
&=(X_{+}\lt\phi^{\prime})\omega_{+}\wedge\omega_{-}+\{(X_{z}\lt\phi^{\prime})\omega_{z}\wedge\omega_{-}+\phi^{\prime}\dd\omega_{-}\} \nn \\
&\qquad+
(X_{-}\lt\phi^{\prime\prime})\omega_{-}\wedge\omega_{+}+\{(X_{z}\lt\phi^{\prime\prime})\omega_{z}\wedge\omega_{+}+\phi^{\prime\prime}\dd\omega_{+}\} \nn \\
&=\{(X_{-}\lt\phi^{\prime\prime})-q^{2}(X_{+}\lt\phi^{\prime})\}\omega_{-}\wedge\omega_{+},
\label{d1fq}
\end{align}
since the terms in curly brackets vanish: $\{(X_{z}\lt\phi^{\prime})\omega_{z}\wedge\omega_{-}+\phi^{\prime}\dd\omega_{-}\}=\{(X_{z}\lt\phi^{\prime\prime})\omega_{z}\wedge\omega_{+}+\phi^{\prime\prime}\dd\omega_{+}\}=0$ from \eqref{dformc3} and \eqref{libu}.
It is then clear that the calculus on the quantum sphere is 2D, and that $\Omega^{2}(\sq)=\Asq\omega_{-}\wedge\omega_{+}=\omega_{-}\wedge\omega_{+}\Asq$, as 
both
$\omega_{\pm}$ commute with elements of $\Asq$ and so does
$\omega_{-}\wedge\omega_{+}$. 

\begin{rema}
{}From \eqref{dformc3} it is natural to ask that $\dd \omega_-=\dd
\omega_+=0$ when restricted to $\sq$. Then, the exterior derivative
of any 1-form  $\alpha = \phi^{\prime} \omega_{-} + \phi^{\prime\prime} \omega_{+} \in \cl^{(0)}_{-2}
\omega_{-} \oplus \cl^{(0)}_{+2} \omega_{+}$ is  given by:
\begin{align}\label{d1f}
\dd \alpha & = \dd (\phi^{\prime} \omega_{-} + \phi^{\prime\prime} \omega_{+}) \nn \\
& = \del \phi^{\prime} \wedge \omega_{-} + \delb \phi^{\prime\prime} \wedge \omega_{+} \nn \\
&= (X_+ \lt \phi^{\prime}  - q^{-2}  X_- \lt \phi^{\prime\prime}) \, \omega_{+}\wedge\omega_{-} \nn \\
& =
q^{-1/2} ( E \lt \phi^{\prime}  - q^{-1} F \lt \phi^{\prime\prime} ) \, \omega_{+}\wedge\omega_{-} ,
\end{align}
since $K\lt$ acts as $q^{\mp}$ on $\cl^{(0)}_{\mp 2}$. Notice that in the
above equality, both $E \lt \phi^{\prime}$ and $F \lt \phi^{\prime\prime}$ belong to $\Asq$, as it
should be.
\end{rema}

 The above results can be summarised in the following
proposition, which is the natural generalisation of the description in \eqref{eals2} of the classical exterior algebra on the sphere manifold $S^{2}$.
\begin{prop}\label{2dsph}
The 2D differential calculus on the sphere $\sq$ is given by:
$$
\Omega({\sq}) = \Asq \oplus \left(\cl^{(0)}_{-2} \omega_{-}
\oplus \cl^{(0)}_{+2} \omega_{+} \right) \oplus \Asq \omega_{+}\wedge\omega_{-} ,
$$
with multiplication rule
$$
\big(f_0; \phi^{\prime},\phi^{\prime\prime};f_2\big)\big(g_0;\psi^{\prime},\psi^{\prime\prime};g_2\big)=\big(f_0g_0; f_0 \psi^{\prime}+\phi^{\prime} g_0,f_0 \psi^{\prime\prime}+\phi^{\prime\prime}g_0;
f_0g_2+f_2g_0+q^{-2} \phi^{\prime\prime}\psi^{\prime}  - \phi^{\prime}\psi^{\prime\prime} \big),
$$
and  exterior derivative $\dd = \delb + \del$:
$$
\begin{array}{ll}
f \mapsto (q^{-1/2} F\lt f, q^{1/2} E \lt f ), & \mathrm{for} \, f\in\Asq , \\
~ & \\
(\phi^{\prime},\phi^{\prime\prime}) \mapsto q^{-1/2} ( E \lt \phi^{\prime}  - q^{-1} F \lt \phi^{\prime\prime} ) , & \mathrm{for} \,
(\phi^{\prime},\phi^{\prime\prime})\in\cl^{(0)}_{-2} \oplus \cl^{(0)}_{+2} .
\end{array}
$$
\end{prop}


\subsubsection{The compatibility conditions between the calculi}\label{se:ccc}

Given the 3D left-covariant differential calculus on $\SU$ described in section \ref{se:lcc}, as well the 1D bicovariant differential calculus on the gauge group algebra $\U(1)$ in section \ref{se:csg}, the `principal bundle compatibility' of these calculi is established by showing that the sequence $\eqref{des}$ is exact. For the case at hand, this sequence becomes
\begin{multline*}
0\, \to\, \ASU\left(\Omega^{1}(\sq)\right)\ASU\, \to \\ \to\,\Omega^{1}(\SU)\,\stackrel{\sim_{\mathcal{N}_{\SU}}}  \longrightarrow \,\ASU\otimes \ker\varepsilon_{\U(1)}/\mathcal{Q}_{\U(1)}\,\to\,0 ,
\end{multline*}
where $\mathcal{Q}_{\U(1)}$ is the ideal given in section~\ref{se:csg} that defines the calculus on $\ca(\U(1))$ and the map $\sim_{\mathcal{N}_{\SU}}$ is defined as in the diagram \eqref{qdia} which now acquires the form:
$$
\begin{array}{lcl}
\Omega^{1}(\SU)_{un} & \stackrel{ \pi_{\cn_{\SU}} }{\longrightarrow} 
& \Omega^{1}(\SU) \\
\downarrow \chi  &  & \downarrow \sim_{\cn_{\SU}} \\
\ASU \otimes \ker\varepsilon_{\U(1)} & \stackrel{ \id\otimes\pi_{\cq_{\U(1)}} }{\longrightarrow} 
&\ASU\otimes (\ker\varepsilon_{\U(1)}/\mathcal{Q}_{\U(1)}) \,.
\end{array}
$$ 
Having a quantum homogeneous bundle, that is a quantum bundle whose total space is a Hopf algebra and whose fiber is a Hopf subalgebra of it, with the differential calculus on the fiber obtained from the corresponding projection, for the above sequence to be exact it is enough \cite{BM97} to check two conditions. The first one is  
$$
(\id\otimes\pi)\circ\Ad(\mathcal{Q}_{\SU}) \,\subset\, \mathcal{Q}_{\SU}\otimes \ca(\U(1))
$$
with $\pi: \ASU \to \ca(\U(1))$ the projection in \eqref{qprp}. This is 
easily established by a direct calculation and using the explicit form of the elements in $\mathcal{Q}_{\SU}$. The second condition amounts to the statement that the kernel of the projection $\pi$ can be written as a right $\ASU$-module of the kernel of $\pi$ itself restricted to the base algebra $\Asq$. Then, one needs to show that 
$
\ker\pi\,\subset\, (\ker\pi|_{\sq}) \ASU 
$,
the opposite implication being obvious. With $\pi$ defined in \eqref{qprp}, one has that 
$$
\ker\pi=\{ cf, \,c^{*}g, \quad \mathrm{with} \quad f,g \in\ASU\}. 
$$
Then $cf=c (a^{*}a+c^{*}c) f = ca^{*} (af) + c^{*}c (cf)$, with both $ca^{*}$ and $c^{*}c$ in $\ker\pi|_{\sq}$. The same holds for elements of the form $c^{*}g$, and the inclusion follows. 

The analysis of the map $\sim_{\cn_{\SU}}:\Omega^{1}(\SU)\to\ASU\otimes\ker\varepsilon_{\U(1)}/\cq_{\U(1)}$ shows that  $\omega_{\pm}\in\,\Omega^{1}(\ASU)$ are indeed the generators of the horizontal forms of the principal bundle, being in the $\ker\sim_{\cn_{\SU}}$. From \eqref{chex} one recovers:
\begin{align*}
&\chi(\delta a)=a\otimes(z-1), \nn \\ 
&\chi(\delta a^{*}) =a^{*}\otimes(z^{*}-1),
\nn \\
&\chi(\delta c)=c\otimes(z-1), \nn \\ 
&\chi(\delta c^{*})=c^{*}\otimes(z^{*}-1) .
\end{align*} 
Given the two generators $\omega_{\pm}$ and the specific $\cq_{\SU}$ which determines the 3D calculus, corresponding universal 1-forms can be taken to be:
\begin{align*}
&\omega_{+}=a\dd c-qc\dd a\qquad\Rightarrow\qquad (a\delta c-qc\delta a) \in[\pi_{\cn_{\SU}}]^{-1}(\omega_{+}), \\
&\omega_{-}=c^{*}\dd a^{*}-qa^{*}\dd c^{*}\qquad\Rightarrow\qquad (c^* \delta a^*- qa^*\delta c^*) \in[\pi_{\cn_{\SU}}]^{-1}(\omega_{-}).
\end{align*}
The action of the canonical map then gives:
\begin{align*}
&\chi(a\delta c-qc\delta a)=(ac-qca)\otimes(z-1)=0, \\ 
&\chi(c^* \delta a^*- qa^*\delta c^*) = (c^*a^* - qa^*c^*)\otimes(z^{*}-1)=0,
\end{align*}
which means that 
\beq
\sim_{\cn_{\SU}}(\omega_{\pm})=0
\label{omho}
\eeq
For the third generator $\omega_{z}$, one shows in a similar fashion that
\beq
\sim_{\cn_{\SU}}(\omega_{z})=1\otimes(\pi_{\cq_{\U(1)}}(z-1)).
\label{omzv}
\eeq 
{}From these it is possible to
conclude that the elements $\omega_\pm$ generate the $\ASU$-bimodule of horizontal forms, while   from \eqref{vvf} one has that the vector 
$X=X_z=(1-q^{-2})^{-1}(1-K^{4})$ is the dual generator to the calculus on the structure Hopf algebra $\ca(\U(1))$. For the corresponding `vector field' $\tilde{X}$ on $\ASU$ as in \eqref{hf}, one has that 
$\tilde{X}(\omega_\pm) = \langle{X},{\sim_{\cn_{\SU}}(\omega_{\pm})}\rangle=0$, while  
$\tilde{X}(\omega_z) = \langle{X},{\sim_{\cn_{\SU}}(\omega_{z})}\rangle=1$. 
These results identify $\tilde{X}$ as a vertical vector field.

\section{A $\star$-Hodge duality on $\Omega(\SU)$ and a Laplacian on $\SU$}\label{se:hssl}

In classical differential geometry a metric structure $\mathrm{g}$ on a N-dimensional manifold $\mathcal{M}$ enables to define a Hodge duality $\star:\Omega^{k}(\mathcal{M})\to\Omega^{N-k}(\mathcal{M})$ on the exterior algebra $\Omega(\mathcal{M})$. The strategy is to consider the volume form $\theta\in\,\Omega^{N}(\mathcal{M})$ associated to a $\mathrm{g}$-orthonormal basis; this corresponds to   the choice of an orientation. Via the Hodge duality it becomes possible to introduce in $\Omega(\mathcal{M})$ both a symmetric bilinear product and a sesquilinear inner product.

The algebraic formulation of geometry of quantum groups, that has been described, presents  no metric tensor. The strategy to introduce a Hodge duality on the exterior algebra $\Omega(\ch)$ coming from a N-dimensional differential calculus on a Hopf algebra $\ch$ is then reversed with respect to the strategy used in the classical setting. The path consists in  defining  a suitable bilinear product on $\Omega(\ch)$ and considering  a volume N-form, from which to induce a $\star$-Hodge structure, using an equation like the one in \eqref{bifopr} as a  definition.

\bigskip

The following description of the quantum formulation of a Hodge duality originates from \cite{kmt}.   
Assume that $\ch$ is a $*$-Hopf algebra equipped with a left covariant calculus $(\Omega^{1}(\ch),\dd)$, with $N$ the dimension of the calculus such that $\dim\,\Omega_{inv}^{N}(\ch)=1$, $\dim\Omega_{inv}^{k}(\ch)=\dim_{inv}^{N-k}(\ch)$.
Suppose also that $\ch$ admits a Haar state $h:\ch\to\IC$, that is a unital linear functional on $\ch$ for which $(id\otimes h)\Delta x=(h\otimes id)\Delta x=h(x)1$ for any $x\in\,\ch$, where $1$ is used to emphasise the unit of the algebra. Suppose further that $h$ is positive, that is $h(x^{*}x)\geq0$ for all $x\in\,\ch$; it is known that the Haar state is unique and automatically faithful: if $h(x^{*}x)=0$, then necessarily $x=0$. One can endow $\ch$ with an inner product derived from $h$, setting: 
\beq 
(x^{\prime};x)_{\ch}=h(x^{*}x^{\prime})
\label{inph}
\eeq
for any $x,x^{\prime}\in\,\ch$. The whole exterior algebra can be endowed with an inner product, defined on a left invariant basis and then extended via the requirement of left invariance,
\beq
(x^{\prime}\omega^{\prime};x\,\omega)_{\ch}=h(x^{*}x^{\prime})(\omega^{\prime},\omega)_{\ch}
\label{inpo}
\eeq
for any $x,x^{\prime}\in\,\ch$ and left invariant forms $\omega,\omega^{\prime}$ in $\Omega(\ch)$.  An inner product is said graded if the spaces $\Omega^{k}(\ch)$ are pairwise orthogonal.

Out of $\Omega^{N}(\ch)$ choose a left invariant hermitian basis element $\theta=\theta^{*}$, which will be called the volume form of the calculus. A linear functional $\int_{\theta}:\Omega(\ch)\to\IC$ -- called the integral on $\Omega(\ch)$ associated to the volume form $\theta\in\,\Omega^{N}(\ch)$ -- is defined by setting $\int_{\theta}\eta=0$ if $\eta$ is a $k$-form with $k<N$, and $\int_{\theta}\eta=h(x)$ if $\eta=x\,\theta$ with $x\in\,\ch$. The differential calculus will be said non-degenerate if, whenever $\eta\in\,\Omega^{k}(\ch)$ and $\eta^{\prime}\wedge\eta=0$ for any $\eta^{\prime}\in\,\Omega^{N-k}(\ch)$, then necessarily $\eta=0$. This property reflects itself in the property of left-faithfulness of the functional $\int_{\theta}$: starting from a non degenerate calculus, it is possible to prove that, if $\eta$ is an element in $\Omega^{k}(\ch)$ for which $\int_{\theta}\eta^{\prime}\wedge\eta=0$ for all $\eta^{\prime}\in\,\Omega^{N-k}(\ch)$, then it is $\eta=0$.   

\begin{prop}
\label{Lge}
Given the exterior algebra $\Omega(\ch)$ coming from a left covariant, non degenerate calculus $(\Omega^{1}(\ch),\dd)$, there exists a unique left $\ch$-linear bijective operator $L:\Omega^{k}(\ch)\to\Omega^{N-k}(\ch)$ for $k=0,\ldots,N$, such that 
\beq
\int_{\theta}\eta^{*}\wedge L(\eta^{\prime})=(\eta^{\prime};\eta)_{\ch}
\label{Lop} 
\eeq
on any $\eta,\eta^{\prime}\in\,\Omega^{k}(\ch)$.
\end{prop}
The proof of this result is in \cite{kmt}, where the  operator $L$ is  called  a Hodge operator. With a left-invariant inner product which is positive definite, i.e. $(\omega,\omega)_{\ch}>0$ for any exterior form $\omega$, the operator $L$
 does not yet define a  $\star$-Hodge structure on $\Omega(\ch)$, since its square does not satisy the natural requirement \eqref{clsH}.   
It is then used to define a new graded left invariant inner product setting on a basis of left invariant forms $\omega\in\, \Omega(\ch)$:  
\begin{align}
&(\omega;\omega^{\prime})_{\ch}^{\natural}=(\omega;\omega^{\prime})_{\ch}, \qquad &\mathrm{on}\,\Omega^{k}(\ch), \,k<N/2; \nn \\
&(\omega;\omega^{\prime})_{\ch}^{\natural}=(L^{-1}(\omega);L^{-1}(\omega^{\prime}))_{\ch}, \qquad &\mathrm{on}\,\Omega^{k}(\ch), \,k>N/2.
\label{nip}
\end{align}
If $N$ is odd, these relations completely define a new left invariant graded inner product on the exterior algebra $\Omega(\ch)$; notice also that assuming the relation \eqref{inph} means that $(1;1)_{\ch}=1$, from which one has  $L(1)=\theta$, so to obtain  in  \eqref{nip} that $(\theta;\theta)^{\natural}_{\ch}=(1;1)_{\ch}=1$.  

In analogy with  \eqref{Lop} define a new Hodge  operator $L^{\natural}:\Omega^{k}(\ch)\to\Omega^{N-k}(\ch)$ via the inner product given in \eqref{nip} as 
\beq
\int_{\theta}\eta^{*}\wedge L^{\natural}(\eta^{\prime})=(\eta^{\prime};\eta)_{\ch}^{\natural}.
\label{Lopd}
\eeq    
Due to the left-faithfulness of the integral, it is clear that $L^{\natural}$ is a well defined bijection, which satisfies the identity $L=L^{\natural}$ when restricted to $\Omega^{k}(\ch)$ with $k<N/2$. Such an operator $L^{\natural}$ is also proved to satisfy $(L^{\natural})^{2}=(-1)^{k(N-k)}$: this is the reason why one can define  a $\star$-Hodge structure on $\Omega(\ch)$ as:
\beq
\star:\Omega^{k}(\ch)\to\Omega^{N-k}(\ch)\,\qquad \star(\eta)=L^{\natural}(\eta).
\label{HsH}
\eeq
The relation \eqref{Lopd} appears as the quantum version of the classical relation \eqref{bifoprip}, which is now used as a definition for the Hodge duality.

If the dimension of the calculus is given by an even $N=2m$, a more specific procedure is needed, The same procedure as before gives a $\star$-Hodge operator on $\Omega^{k}(\ch)$ for $k\neq m$ via the inner product \eqref{nip}. Using the volume form $\theta\in\,\Omega^{N}(\ch)$  set now a sesquilinear form
\beq
\hs{\eta^{\prime}}{\eta}=\int_{\theta}\eta^{*}\wedge\eta^{\prime},
\label{ses}
\eeq
which is non-degenerate by the faithfulness of the integral $\int_{\theta}$. 
The $\ch$-bimodule $\Omega^{m}(\ch)$ has a basis of $\left(\begin{array}{c} 2m \\ m\end{array} \right)$ left invariants elements $\omega_{a}$. The restriction of \eqref{ses} to elements $\omega_{a}$ defines a sesquilinear form on the vector space $\Omega^{m}_{\mathrm{inv}}$: this form is hermitian  if $(-1)^{m^{2}}=1$, and anti-hermitian if $(-1)^{m^{2}}=-1$, so it can be 'diagonalised'. There exists a basis $\check{\omega}_{j}\in\,\Omega^{m}(\ch)$ such that one has $\hs{\check{\omega}_{a}}{\check{\omega}_{b}}=\pm\delta_{ab}$ if it is hermitian, and $\hs{\check{\omega}_{a}}{\check{\omega}_{b}}=\pm i\delta_{ab}$ if it is anti-hermitian. It is then possible to use such a basis to define a left $\ch$-linear operator $\mathfrak{L}:\Omega^{m}(\ch)\to\Omega^{m}(\ch)$ setting on the basis
\beq
\mathfrak{L}(\check{\omega}_{a})=(-1)^{m^{2}}\hs{\check{\omega}_{a}}{\check{\omega}_{a}}\,\check{\omega}_{a}.
\label{mL}
\eeq
(no sum on $a$). This map is a bijection, and satisfies $\mathfrak{L}^{2}=(-1)^{m^{2}}$, so a $\star$-Hodge structure on $\Omega^{m}_{\mathrm{inv}}(\ch)$ can be defined as: 
\beq
\star(\check{\omega}_{a})=\mathfrak{L}(\check{\omega}_{a}),
\label{HL}
\eeq
and extended on any $\eta\in\,\Omega^{m}(\ch)$ by the requirement of left linearity, thus giving a complete constructive procedure for a $\star$-Hodge structure on $\Omega(\ch)$.  The Hodge operator $\mathfrak{L}:\Omega^{m}(\ch)\to\Omega^{m}(\ch)$ is then used to introduce a left invariant inner product on $\Omega^{m}(\ch)$, defined by:
\beq
(\omega_{a};\omega_{b})^{\natural}_{\ch}=\int_{\theta}\omega_{b}^{*}\wedge\mathfrak{L}(\omega_{a}),
\label{inpm}
\eeq
on a basis of left invariant $\{\omega_{a}\}$ 2-forms, and then extended via the requirement of left invariance as in \eqref{inpo}.    It is easy to see  that the definition eventually gives the inner product 
\beq
(\check{\omega}_{a};\check{\omega}_{b})^{\natural}_{\ch}=\delta_{ab}.
\label{inpme}
\eeq

\subsection{A $\star$-Hodge structure on $\Omega(\SU)$} 

This section describes how the outlined  procedure yields   a left invariant inner product on the exterior algebra $\Omega(SU_{q}(2))$ generated by the left covariant 3D calculus from  section \ref{se:lcc}, and the way it gives rise to a $\star$-Hodge structure. Such a $\star$-Hodge structure will be then used to define a Laplacian operator on $\ASU$, which is completely diagonalised. 

The Hopf algebra $\ASU$ has a Haar state $h:\ASU\to\IC$,  which is  positive, unique and authomatically faithful.  From  the Peter-Weyl decomposition of $\ASU$ in terms of the vector space basis elements $w_{p:r,t}\in\,W_{p}$ \eqref{ws},  the Haar state is determined by setting:
$$
h(1)=1\,\qquad h(w_{p:r,t})=0\,\,\,\forall\,p\geq0.
$$
The algebraic relations \eqref{derel} among the generators of $\ASU$ makes it then possible to prove that the only non trivial action of $h$ on $\ASU$ can also be written as:
$$
h((cc^{*})^{k})=(\sum_{j=0}^{k}\,q^{2j})^{-1}=\frac{1}{1+q^{2}+\ldots+q^{2k}},
$$ 
with $k\in\,\IN$. One can define  on $\ASU$ an inner product derived from $h$, setting:
\beq
(x^{\prime},x)_{ \SU}=h(x^{*}x^{\prime})
\eeq
with $x,x^{\prime}\in\,\ASU$. The differential  3D calculus being left covariant, the set of $k$-forms $\Omega^{k}(\SU)$  has a basis of left invariant forms. The exterior algebra  $\Omega(\SU)$ is endowed with an inner product, defined  on a left invariant basis and extended via the requirement of left invariance:
$$
(x^{\prime}\omega^{\prime},x\,\omega)_{\SU}=h(x^{\star}x^{\prime})(\omega^{\prime},\omega)_{\SU}
$$
for all $x,x^{\prime}$ in $\ASU$ and $\omega, \omega^{\prime}\in\,\Omega(\SU)$ left invariant forms. 
Assume the top form $\theta=\alpha^{\prime}\omega_{-}\wedge\omega_{+}\wedge\omega_{z}$ as volume form, with $\alpha^{\prime}\in\,\IR$ so that $\theta^{\star}=\theta$. The integral on the exterior algebra $\Omega(\SU)$ associated to the volume form $\theta$ is defined by $\int_{\theta}\eta=0$ if $\eta$ is a $k$-form with $k<2$, and $\int_{\theta}\eta=h(x)$ if $\eta=x\,\theta$. This integral is left-faithful. 

Set a left invariant graded inner product by assuming that the only non-zero products among left invariant forms are:
\begin{align}
&(1,1)_{\SU}=1, \nn \\
&(\theta,\theta)_{\SU}=1;
\label{sp03}
\end{align}
while in  $\Omega^{1}(\SU)$ are:
\begin{align}
&(\omega_{-},\omega_{-})_{\SU}=\beta, \nn \\
&(\omega_{+},\omega_{+})_{\SU}=\nu, \nn \\
&(\omega_{z},\omega_{z})_{\SU}=\gamma
\label{sp1}
\end{align}
with $\beta, \nu, \gamma\in\,\IR$, and:
\begin{align}
&(\omega_{-}\wedge\omega_{+},\omega_{-}\wedge\omega_{+})_{\SU}=1, \nn \\
&(\omega_{+}\wedge\omega_{z},\omega_{+}\wedge\omega_{z})_{\SU}=1, \nn \\
&(\omega_{z}\wedge\omega_{-},\omega_{z}\wedge\omega_{-})_{\SU}=1
\label{sp2}
\end{align}
in $\Omega^{2}(\SU)$. This choice comes as the most natural in order to mimic the properties of the classical inner product \eqref{bifoprpm}, coming from the classical Hodge structure \eqref{clHs} originated from the metric \eqref{gmet}. 
 The Hodge operator defined in  \eqref{Lop} is: 
\begin{align}
&L(1)=\alpha^{\prime}\,\omega_{-}\wedge\omega_{+}\wedge\omega_{z}, \nn \\
&L(\omega_{-})=-\alpha^{\prime}\beta q^{-6}\,\omega_{z}\wedge\omega_{-}, \nn \\
&L(\omega_{+})=-\alpha^{\prime}\nu\,\omega_{+}\wedge\omega_{z}, \nn \\
&L(\omega_{z})=-\alpha^{\prime}\gamma\,\omega_{-}\wedge\omega_{+}, \nn \\
&L(\omega_{-}\wedge\omega_{+})=-\alpha^{\prime}\,\omega_{z}, \nn \\
&L(\omega_{+}\wedge\omega_{z})=-\alpha^{\prime}\,\omega_{+}, \nn \\
&L(\omega_{z}\wedge\omega_{-})=-\alpha^{\prime}\,\omega_{-}, \nn \\
&L(\omega_{-}\wedge\omega_{+}\wedge\omega_{z})=\alpha^{\prime-1}.
\label{fL}
\end{align}
The Hodge operator $L$  is  used to define a new graded left invariant inner product on $\Omega(\SU)$, as:
\begin{align} 
&(\omega^{\prime},\omega)_{\SU}^{\natural}=(\omega^{\prime},\omega)_{\SU} &\qquad\mathrm{on}\,\,\,\Omega^{k}(\SU),\,k=0,1;\nn \\
&(\omega^{\prime},\omega)_{\SU}^{\natural}=(L^{-1}(\omega^{\prime}),L^{-1}(\omega))_{\SU}&\qquad\mathrm{on}\,\,\,\Omega^{k}(\SU),\,k=2,3,
\label{dLd}
\end{align}
on the basis of left invariant forms. On $\Omega^{k}(\SU))$ -- with $k=2,3$ -- one has:
\begin{align}
&(\omega_{-}\wedge\omega_{+},\omega_{-}\wedge\omega_{+})^{\natural}_{\SU}=\alpha^{\prime-2}\gamma^{-1}, \nn \\
&(\omega_{+}\wedge\omega_{z},\omega_{+}\wedge\omega_{z})^{\natural}_{\SU}=\alpha^{\prime-2}\nu^{-1}, \nn \\
&(\omega_{z}\wedge\omega_{-},\omega_{z}\wedge\omega_{-})^{\natural}_{\SU}=q^{12}\alpha^{\prime-2}\beta^{-1}, \nn \\
&(\theta,\theta)^{\natural}_{\SU}=1.
\label{sp23n}
\end{align}
Associated to this new inner product there is in analogy a new unique left $\ASU$-linear operator $L^{\natural}:\Omega^{k}(\SU)\to\Omega^{3-k}(\SU)$ defined by $\int_\theta\,\eta^{*}\wedge L^{\natural}(\eta^{\prime})=(\eta^{\prime},\eta)^{\natural}$, which is a bijection. This operator     
is such that $(L^{\natural})^{2}=(-1)^{k(3-k)}=1$, so following \eqref{HsH} one has a $\star$-Hodge structure on the exterior algebra $\Omega(\SU)$:
\beq
\star:\Omega^{k}(\SU)\to\Omega^{3-k}(\SU)\qquad\star(\eta)=L^{\natural}(\eta),
\eeq
given by: 
\begin{align}
&\star(1)=\theta=\alpha^{\prime}\,\omega_{-}\wedge\omega_{+}\wedge\omega_{z}, \nn \\
&\star(\omega_{-})=-\alpha^{\prime}\beta q^{-6}\,\omega_{z}\wedge\omega_{-}, \nn \\
&\star(\omega_{+})=-\alpha^{\prime}\nu\,\omega_{+}\wedge\omega_{z}, \nn \\
&\star(\omega_{z})=-\alpha^{\prime}\gamma\,\omega_{-}\wedge\omega_{+}, \nn \\
&\star(\omega_{-}\wedge\omega_{+})=-\alpha^{\prime-1}\gamma^{-1}\,\omega_{z}, \nn \\
&\star(\omega_{+}\wedge\omega_{z})=-\alpha^{\prime-1}\nu^{-1}\,\omega_{+}, \nn \\
&\star(\omega_{z}\wedge\omega_{-})=-\alpha^{\prime-1}\beta^{-1}q^{6}\,\omega_{-}, \nn \\
&\star(\omega_{-}\wedge\omega_{+}\wedge\omega_{z})=\alpha^{\prime-1}.
\label{sH}
\end{align}

\begin{rema}
\label{relh}
The definition of the graded left invariant inner product $(\cdot,\cdot)^{\natural}_{\SU}$ in \eqref{dLd} shows that, in order to   have a $\star$-Hodge structure on the exterior algebra $\Omega(\SU)$ generated by the 3D calculus, it is sufficient  to choice an hermitian volume form and  a graded left invariant inner product \emph{only} on $\Omega^{k}(\SU)$ for $k=0,1$. This is a general aspect: given  a Hopf $*$-algebra $\ch$, equipped with a finite odd $N$ dimensional left covariant differential calculus,  the formalism developed in \cite{kmt} shows that what one  needs is an hermitian volume form and a graded left invariant inner product on $\Omega^{k}(\ch)$ for $k<N/2$.     
\end{rema}

\subsubsection{A Laplacian operator on $\ASU$}

Given a differential calculus and a $\star$-Hodge structure on the Hopf algebra $\ASU$ it is possible to define a scalar Laplacian operator $\Box_{\SU}:\ASU\to\ASU$ as $\Box_{\SU}\phi=\star\dd\star\dd\phi$ for any $\phi\in\,\ASU$. This Laplacian can be written down by a computation  on the basis of the left invariant forms of the calculus:
\begin{align}
\dd\phi&=(X_{+}\lt\phi)\omega_{+}+(X_{-}\lt\phi)\omega_{-}+(X_{z}\lt\phi)\omega_{z}; \nn \\
\star\dd\phi&=-\alpha^{\prime}[\nu(X_{+}\lt\phi)\omega_{+}\wedge\omega_{z}+\beta q^{-6}(X_{-}\lt\phi)\omega_{z}\wedge\omega_{-}+\gamma(X_{z}\lt\phi)\omega_{-}\wedge\omega_{+}]. \nn
\end{align}
The last line comes from \eqref{sH} and the left linearity of the $\star$-Hodge on the exterior algebra $\Omega(\SU)$. By \eqref{dformc3} the derivative $\dd$ acts on the previous 2-form as:
\begin{align} 
\dd\star\dd\phi&=-\alpha^{\prime}[\nu(X_{-}X_{+}\lt\phi)(\omega_{-}\wedge\omega_{+}\wedge\omega_{z})+\beta q^{-6}(X_{+}X_{-}\lt\phi)(\omega_{+}\wedge\omega_{z}\wedge\omega_{-})+\gamma(X_{z}X_{z}\lt\phi)(\omega_{z}\wedge\omega_{-}\wedge\omega_{+})] \nn \\
&=-\alpha^{\prime}\{[\nu X_{-}X_{+}+\beta X_{+}X_{-}+\gamma X_{z}X_{z}]\lt\phi\}(\omega_{-}\wedge\omega_{+}\wedge\omega_{z}), \nn
\end{align}
where the commutation rules \eqref{commc3} have been used.  The last of \eqref{sH} finally gives the Laplacian operator the expression:
\beq
\star\dd\star\dd\phi=-[\nu X_{-}X_{+}+\beta X_{+}X_{-}+\gamma X_{z}X_{z}]\lt\phi
\label{la3}
\eeq
in terms of the left action of the quantum vector fields of the calculus. The  expression \eqref{la3} shows that $\Box_{\SU}:\cl_{n}\to\cl_{n}$. This operator can be diagonalised. One has to recall the decomposition \eqref{decoln} of the modules $\cl_{n}$ for the right action of $\su$:  this right action leaves invariant the eigenspaces of the  Laplacian  since left and right actions of $\su$ on $\ASU$ do commute. On each irreducible subspace $V_{J}^{(n)}$ \eqref{decoln} for the right action of $\su$ one has a basis $\phi_{n,J,l}=(c^{J-n/2}a^{\star J+n/2})\triangleleft E^{l}=w_{2J:J-\frac{n}{2},l}$ (with $l=0,\ldots,2J$)  of eigenvectors \eqref{ws} for the Laplacian. The spectrum of the Laplacian does not depend on the integer $l$:  an explicit computation shows that
\begin{align}
X_{z}X_{z}\triangleright \phi_{n,J,l}&=q^{2(n+1)}[n]^{2}\phi_{n,J,l}, \nonumber \\
X_{+}X_{-}\triangleright \phi_{n,J,l}&=q^{n-1}([J-\frac{n}{2}][J+1+\frac{n}{2}]+[n]) \phi_{n,J,l}, \nonumber \\
X_{-}X_{+}\triangleright \phi_{n,J,l}&=q^{n+1}([J-\frac{n}{2}][J+1+\frac{n}{2}]) \phi_{n,J,l}.
\label{egv}
\end{align}
The spectrum of the Laplacian \eqref{la3} is then given as $\Box_{\SU}\phi_{n,J,l}=\lambda_{n,J,l}\phi_{n,J,l}$ with:
\beq
\lambda_{n,J,l}=-q^{n}\{\nu q[J-\frac{n}{2}][J+1+\frac{n}{2}]+\beta q^{-1}([J-\frac{n}{2}][J+1+\frac{n}{2}]+[n])+\gamma q^{n+2}[n]^{2}\}.
\label{la3e}
\eeq

\section{A $\star$-Hodge structure on $\Omega(\sq)$ and a Laplacian operator on~$\Asq$}\label{se:hlss}

The way the $\star$-Hodge structure \eqref{sH} has been introduced on $\Omega(\SU)$  comes  from the analysis in \cite{kmt}. The aim of this section is to  extend that procedure in order to introduce a $\star$-Hodge structure on $\Omega(\sq)$. The strategy is to directly follow the same path, and to apply  to the differential calculus $\Omega(\sq)$ the same procedure, explicitly checking its consistency in the new setting.

\subsection{A $\star$-Hodge structure on $\Asq$}

The differential calculus on the quantum sphere $\sq$ has been described in section \ref{se:cals2} and fully presented in proposition \ref{2dsph}. It is a 2D left covariant calculus: as a volume form  consider $\check{\theta}=i\alpha^{\prime\prime}\omega_{-}\wedge\omega_{+}$.
\begin{lemm}
\label{2dnd}
The 2D calculus $\Omega({\sq})$ from  proposition \ref{2dsph} is non degenerate. 
\begin{proof}
The proof of this lemma is direct. To be definite, consider a $0$-form $\eta=f$ with $f\in\,\Asq\simeq\cl_{0}^{(0)}$, so to have a product 
$$
\eta^{\prime}\wedge\eta=f^{\prime}(\omega_{-}\wedge\omega_{+})f=f^{\prime}f\,\omega_{-}\wedge\omega_{+}
$$ 
from the commutation rules in \eqref{bi1}, where $\eta^{\prime}=f^{\prime}\,\omega_{-}\wedge\omega_{+}$ with $f^{\prime}\in\,\cl_{0}^{(0)}$. One has $\eta^{\prime}\wedge\eta=0\,\Leftrightarrow\,f^{\prime}f=0$: such a relation is satisfied for any $f^{\prime}\in\,\cl_{0}^{(0)}$ iff $f=0$.

Consider now the 1-form $\eta=x\,\omega_{-}$ with $x\in\,\cl_{-2}^{(0)}$, so to have a product 
$$
\eta^{\prime}\wedge\eta=(x^{\prime}\omega_{-}+y^{\prime}\omega_{+})\wedge x\,\omega_{-}=-y^{\prime}x\,\omega_{-}\wedge\omega_{+}
$$ 
where $(x^{\prime},y^{\prime})\in\,(\cl_{-2}^{(0)},\cl_{2}^{(0)})$. The relation $\eta^{\prime}\wedge\eta=0\,\Leftrightarrow\,y^{\prime}x=0$ is satisfied for any $y^{\prime}\in\,\cl_{2}^{(0)}$ iff $x=0$. The remaining cases can be analogously analysed,  thus proving the claim.
\end{proof}
\end{lemm}

The restriction of the Haar state $h$ to $\Asq$ yields a faithful, invariant -- that is $h(f\rt X)=h(f)\varepsilon(X)$  for $f\in\,\Asq$ and $X\in\,\su$ -- state on $\Asq$, allowing the definition of an integral 
$\int_{\check{\theta}}:\Omega(\sq)\to \IC$ given by:  
\begin{align}
&\int_{\check{\theta}}f=0, \qquad &\mathrm{on}\,f\in\,\Asq, \nn \\
&\int_{\check{\theta}}\eta=0,\qquad &\mathrm{on}\,\eta\in\,\Omega^{1}(\sq), \nn \\
&\int_{\check{\theta}}f\,\omega_{-}\wedge\omega_{+}=-i\alpha^{\prime\prime-1}\,h(f). 
\label{ints2}
\end{align}
\begin{lemm}
\label{inlf}
The integral $\int_{\check{\theta}}$ defined in \eqref{ints2} is left-faithful. 
\begin{proof}
The proof of this result is also direct. Consider, to be definite, the 1-form $\eta=x\omega_{-}$ with $x\in\,\cl_{-2}^{(0)}$, and a generic $\eta^{\prime}=x^{\prime}\omega_{-}+y^{\prime}\omega_{+}\,\in\Omega^{1}(\sq)$. The relation $\int_{\check{\theta}}\eta^{\prime}\wedge\eta=0$ for any $\eta^{\prime}\in\,\Omega^{1}(\sq)$ is equivalent to the condition $h(y^{\prime}x)=0\, \,\,\forall\,y^{\prime}\in\,\cl_{2}^{(0)}$. Since this last equality must be valid for any $y^{\prime}\in\cl_{2}^{(0)}$, choosing  $y^{\prime}=x^{*}$, it results $h(x^{*}x)=0$: the faithfulness of the Haar state $h$ then gives $x=0$. The claim of the lemma is proved by an analogous analysis on the remaining cases.  
\end{proof}
\end{lemm}

The restriction to $\Omega(\sq)$ of the left invariant graded
product  \eqref{dLd} on $\Omega(\SU)$, which is the one compatible with the $\star$-Hodge structure,  gives a left $\Asq$-invariant graded inner product:
\begin{align}
&(1,1)_{\sq}=1; \nn \\
&(x^{\prime}\omega_{-}+y^{\prime}\omega_{+},x\,\omega_{-}+y\,\omega_{+})_{\sq}=h(x^{*}x^{\prime})\beta+h(y^{*}y^{\prime}) \nu; \nn \\
&(\omega_{-}\wedge\omega_{+},\omega_{-}\wedge\omega_{+})_{\sq}=\alpha^{\prime-2}\gamma^{-1},
\label{inp2}  
\end{align}
with , $x,x^{\prime}\in\,\cl_{-2}^{(0)}$ and  $y,y^{\prime}\in\,\cl_{2}^{(0)}$. Recalling proposition \ref{Lge} -- namely  equation \eqref{Lop} -- and the results proved in lemmas \ref{2dnd} and \ref{inlf},  a left $\Asq$-linear  Hodge operator $L:\Omega^{k}(\sq)\to\Omega^{2-k}(\sq)$ can be defined  for $k=0,2$. From the first line in the inner product relation \eqref{inp2} one has $L(1)=\check{\theta}$, while the third gives $L(\check{\theta})=\alpha^{\prime\prime2}\alpha^{\prime-2}\gamma^{-1}$. It is evident that for such an Hodge operator it is $L^{2}\neq1$, which is a natural requirement for a $\star$-Hodge structure on $\Omega^{k}(\sq)$ for $k=0,2$. On the exterior algebra $\Omega(\SU)$ this problem was solved by changing the inner product via the definition \eqref{dLd}, and proving that the new Hodge operator does satisy all the required properties 
to have a consistent $\star$-Hodge. Following an analogous path, define 
\begin{align}
&(1,1)^{\natural}_{\sq}=1, \nn \\
&(x^{\prime}\omega_{-}+y^{\prime}\omega_{+},x\,\omega_{-}+y\,\omega_{+})^{\natural}_{\sq}=
(x^{\prime}\omega_{-}+y^{\prime}\omega_{+},x\,\omega_{-}+y\,\omega_{+})_{\sq}, \nn \\
&(\check{\theta},\check{\theta})_{\sq}^{\natural}=(L^{-1}(\check{\theta}),L^{-1}(\check{\theta}))_{\sq}=1,
\label{inp2ns}  
\end{align}
where the inner products on 1-forms amounts to a different labelling of the inner product in \eqref{inp2}. The Hodge operator on $\Omega^{k}(\sq)$ for $k=0,2$ relative to such a new inner product is given by $L^{\natural}(1)=\check{\theta}$ and $L^{\natural}(\check{\theta})=1$. But now the inner product has changed: the requirement that the inner product $(~,~)^{\natural}_{\SU}$ on the exterior algebra $\Omega(\SU)$ fixed -- via a restriction, as given in \eqref{inp2} --  the inner product $(~,~)_{\sq}$ on the exterior algebra  
$\Omega(\sq)$ implies 
that the condition 
\beq
(\check{\theta},\check{\theta})_{\sq}^{\natural}=(\check{\theta},\check{\theta})_{\SU}^{\natural}
\label{1con}
\eeq
\emph{has to be imposed}, giving 
\beq
\alpha^{\prime\prime2}\alpha^{\prime-2}\gamma^{-1}=1
\label{con1}
\eeq
 as a constraint among the parameters.  The constraint  \eqref{1con} can be interpreted as the quantum analogue of fixing the classical metric on the basis $S^{2}$ of the Hopf bundle as the contraction of the Cartan-Killing metric on $S^{3}\sim SU(2)$, since that choice in the classical formalism, as stressed in remark \ref{omps}, gives the equality of the inner product on $\Omega(S^{2})$ defined in \eqref{phs2cl} with the restriction of the inner product on $\Omega(S^{3})$ given in \eqref{bifoprip}.

The differential calculus on $\sq$ is even dimensional with $N=2$, so on $\Omega^{1}(\sq)$ define a sesquilinear form:
\beq
\hs{\eta^{\prime}}{\eta}=\int_{\check{\theta}}\eta^{*}\wedge\eta^{\prime}=i\alpha^{\prime\prime-1}
\{h(y^{*}y^{\prime})-q^{2}h(x^{*}x^{\prime})\}
\label{sesq2}
\eeq
where $\eta=x\,\omega_{-}+y\,\omega_{+}$ and $\eta^{\prime}=x^{\prime}\omega_{-}+y^{\prime}\omega_{+}$, with $x,x^{\prime}\in\,\cl_{-2}^{(0)}$ and $y,y^{\prime}\in\,\cl_{2}^{(0)}$. The quantum sphere $\sq$ is a quantum homogeneous space and not a Hopf algebra, so there is no left-invariant basis in $\Omega^{1}(\sq)$: neverthless such a sesquilinear form  can be "diagonalised", as 
\begin{align}
&\hs{x\,\omega_{-}}{x\,\omega_{-}}=-iq^{2}\alpha^{\prime\prime-1}\,h(x^{*}x); \nn \\
&\hs{y\,\omega_{+}}{y\,\omega_{+}}=i\alpha^{\prime\prime-1}\,h(y^{*}y),
\label{sesd2}
\end{align}
where the faithfulness of the Haar state ensures that the coefficients on the right hand side of these expressions never vanish. The general result from \cite{kmt} -- recalled in \eqref{mL} -- is no longer valid on a quantum homogeneous space: the diagonalisation in \eqref{sesd2} suggests indeed a way to define a Hodge operator.  Since $\alpha^{\prime\prime}$ can be both positive or negative, define
\begin{align}
&x\,\theta_{-}=q^{-1}\left(\frac{|\alpha^{\prime\prime}|}{h(x^{*}x)}\right)^{1/2}x\,\omega_{-},\nn \\
&y\,\theta_{+}=\left(\frac{|\alpha^{\prime\prime}|}{h(y^{*}y)}\right)^{1/2}y\omega_{+}
\label{xtyt}
\end{align}
so to have from \eqref{sesd2}:
\begin{align}
\hs{x\,\theta_{-}}{x\,\theta_{-}}=-i\frac{|\alpha^{\prime\prime}|}{\alpha^{\prime\prime}}, \nn \\
\hs{y\,\theta_{+}}{y\,\theta_{+}}=i\frac{|\alpha^{\prime\prime}|}{\alpha^{\prime\prime}}.
\label{sest}
\end{align}
In the same way as in \eqref{mL}, define a left $\Asq$-linear operator
 $\mathfrak{L}:\Omega^{1}(\sq)\to\Omega^{1}(\sq)$ setting:
\begin{align}
&\mathfrak{L}(x\,\theta_{-})=i\frac{|\alpha^{\prime\prime}|}{\alpha^{\prime\prime}}x\,\theta_{-}, \nn \\  
&\mathfrak{L}(y\,\theta_{+})=-i\frac{|\alpha^{\prime\prime}|}{\alpha^{\prime\prime}}y\,\theta_{+}.  
\label{matL}
\end{align}
Such an operator clearly satisfies the condition $\mathfrak{L}^{2}=-1$ for any value of $\alpha^{\prime\prime}$. It is not yet a consistent Hodge operator: it has to be compatible with the left invariant inner product on $\Omega^{1}(\sq)$ obtained in \eqref{inp2ns} as a restriction of the analogue on $\Omega^{1}(\SU)$. 
From the relation \eqref{inpm}, this compatibility  must be imposed:
\beq
(\eta^{\prime},\eta)_{\sq}^{\natural}=\int_{\check{\theta}}\,\eta^{*}\wedge\mathfrak{L}(\eta^{\prime}).
\label{Hos2}
\eeq
This condition is fulfilled if and only if the parameters in this formulation satisfy:
\beq
|\alpha^{\prime\prime}|\beta=q^{2},
\label{con2}
\eeq
\beq
|\alpha^{\prime\prime}|\nu=1.
\label{con3}
\eeq 
The $\star$-Hodge structure on $\Omega(\sq)$ is defined as a left $\Asq$-linear operator whose action is given by:
\begin{align}
&\star(1)=i\alpha^{\prime\prime}\,\omega_{-}\wedge\omega_{+}, \nn \\
&\star(x\,\omega_{-})=i\frac{|\alpha^{\prime\prime}|}{\alpha^{\prime\prime}} (x\,\omega_{-}), \nn \\
&\star(y\,\omega_{+})=-i\frac{|\alpha^{\prime\prime}|}{\alpha^{\prime\prime}} (y\,\omega_{+}), \nn \\
&\star(i\omega_{-}\wedge\omega_{+})=\alpha^{\prime\prime-1},
\label{sHs}
\end{align} 
with the parameters $\alpha^{\prime},\alpha^{\prime\prime},\beta,\nu,\gamma$ satisfying the constraints \eqref{con1}, \eqref{con2}, \eqref{con3}.
\begin{rema}
The $\star$-Hodge structure \eqref{sHs} differs from the one in \cite{ma05}, because in that paper  the $\star$-Hodge structure was required to satisfy the relation $\star^{2}=1$, while the path followed here is to remain consistent with the requirement that $\star^{2}=(-1)^{k(N-k)}$ on $k$-forms from a $N$-dimensional calculus.
\end{rema}
The definition \eqref{sHs} of the Hodge duality is still not complete. The constraints among the parameters involve the absolute value of $\alpha^{\prime\prime}$, so one still needs to choose their relative signs.  In the classical setting the only parameter was $\alpha\in\,\IR$, and it has been chosen positive so to give a riemannian metric $g$ in the analysis of section \ref{se:LS}.  As it is clear from \eqref{bifocl} and from the definition \eqref{bifopr}, the positivity of the metric implies the positivity of the  symmetric form $\hs{~}{~}_{S^{3}}$ \eqref{bifocl} and of the sesquilinear inner product $\hs{~}{~}_{S^{3}}^{\sim}$ \eqref{bifopr}: the signature of the metric tensor implies the signature of both the bilinear forms

In the quantum setting, having no metric tensor,  the choice of the relative signs of the parameters is equivalent to choose the signature of the left-invariant inner product \eqref{sp1} on $\Omega^{1}(\SU)$: this  will encode  a specific  metric signature.  

The natural choice for a riemannian signature is, from \eqref{sp1} and  \eqref{sp23n}, given by $\beta,\nu,\gamma\in\,\IR_{+}$.  This choice turns out to be 
 compatible with \eqref{con1}, \eqref{con2} and \eqref{con3} for every  $\alpha^{\prime}$ and $\alpha^{\prime\prime}$.  From \eqref{con2} and \eqref{con3} one also has  that:
\beq
\label{conr}
\beta=q^{2}\nu.
\eeq
This relation has a number of interesting and important consequences, described in the next propositions. 
\begin{prop}
\label{prost2}
The $\star$-Hodge structure given as a left $\Asq$-linear map $\star:\Omega^{k}(\sq)\to\Omega^{2-k}(\sq)$ for $k=0,1,2$ and defined by \eqref{sHs},  
has the property\footnote{In the classical formalism, the $\star$-Hodge structure on an exterior algebra coming from a N dimensional differential calculus  $\star:\Omega^{k}(\ch)\mapsto\Omega^{N-k}(\ch)$ satisfies the identity \eqref{symhs}: 
$$
\eta\wedge(\star\eta^{\prime})=\eta^{\prime}\wedge(\star\eta)
$$
to which the identity \eqref{stwe} reduces in the classical limit.}
\beq
\star(\eta)\wedge\eta^{\prime}=(-1)^{k(2-k)}\eta\wedge\star(\eta^{\prime})
\label{stwe}
\eeq
for any $\eta,\eta^{\prime}\in\,\Omega^{k}(\sq)$. 
\begin{proof}
The relation is trivially satisfied for $k=0,2$. Consider now the two elements $\eta=x\,\omega_{-}+y\,\omega_{+}$ and $\eta^{\prime}=x^{\prime}\omega_{-}+y^{\prime}\omega_{+}$ in $\Omega^{1}(\sq)$, which means  $x,x^{\prime}\in\,\cl_{-2}^{(0)}$ and $y,y^{\prime}\in\,\cl_{2}^{(0)}$  by proposition \ref{2dsph}. The multiplication rule from the same proposition gives:
\begin{align}
(\star\eta)\wedge\eta^{\prime}&=i\alpha^{\prime\prime}(\beta\,xy^{\prime}+\nu\,yx^{\prime})\omega_{-}\wedge\omega_{+}, \nn \\
\eta\wedge(\star\eta^{\prime})&=-i\alpha^{\prime\prime}(q^{-2}\beta\,yx^{\prime}+q^{2}\nu\,xy^{\prime})\omega_{-}\wedge\omega_{+}.
\label{we1}
\end{align}
The two expression are equal -- up to the sign, which is the claim of the proposition --  from \eqref{conr}. 
\end{proof}
\end{prop}

\begin{prop}
\label{rlst}
The left $\Asq$-linear $\star$-Hodge map defined  by \eqref{sHs} is right $\Asq$-linear: given $\eta\in\,\Omega(\sq)$, it is $\star(\eta f)=\star(\eta)f$ for any $f\in\,\Asq$.
\begin{proof}
The 2D differential calculus on the quantum sphere $\sq$ has the specific property, coming from the bimodule structure \eqref{bi1} of $\Omega^{1}(\SU)$ -- where one has  $\omega_{\pm}\phi=q^{n}\phi\,\omega_{\pm}$ for any $\phi\in\,\cl_{n}^{(0)}$ --  that $\omega_{\pm}f=f\omega_{\pm}$ with $f\in\,\cl_{0}^{(0)}\simeq\Asq$. The claim of the proposition is trivial for $\eta\in\Omega^{0}(\sq)\simeq\Asq$. For a 1-form $\eta=x\omega_{-}+y\omega_{+}$ in $\Omega^{1}(\sq)$, one has:
$$
\star(\eta f)=\star((x\omega_{-}+y\omega_{+})f)=\star(xf\omega_{-}+yf\omega_{+})=i\alpha^{\prime\prime}\nu(xf\omega_{-}-yf\omega_{+})=i\alpha^{\prime\prime}\nu(x\omega_{-}-y\omega_{+})f=\star(\eta)f.
$$
An analogue chain of equalities is valid for $\eta=f^{\prime}\omega_{-}\wedge\omega_{+}\in\,\Omega^{2}(\sq)$, with $f^{\prime}\in\,\Asq$. 
\end{proof}
\end{prop}

In the same way it is possible to prove the following identities, which will be explicitly used in the analysis of the gauged Laplacian operator, and which slightly generalise the last proposition.

\begin{lemm}
\label{cala}
Given the left $\Asq$-linear $\star$-Hodge map defined  by \eqref{sHs},  with $\phi\in\,\cl_{n}^{(0)}$, $\phi^{\prime}\in\,\cl_{-n}^{(0)}$ and $\eta\in\,\Omega^{1}(\sq)$ one has:
\begin{align*}
&\star(\phi^{\prime}\eta\phi)=\phi^{\prime}(\star\eta)\phi,\\
&\star(\phi^{\prime}(\omega_{-}\wedge\omega_{+})\phi)=q^{2n}\phi^{\prime}\{\star(\omega_{-}\wedge\omega_{+})\}\phi.
\end{align*}
\begin{proof}
With $\phi^{\prime}\eta\phi\in\,\Omega^{1}(\sq)$, and again $\eta=x\omega_{-}+y\omega_{+}$, it is explicitly:
$$
\star(\phi^{\prime}\eta\phi)=\star(\phi^{\prime}q^{n}(y\phi\omega_{+}+x\phi\omega_{-}))=-iq^{n}\alpha^{\prime\prime}\nu\,\phi^{\prime}(y\phi\omega_{+}-x\phi\omega_{-})=-i\alpha^{\prime\prime}\nu\,\phi^{\prime}(y\omega_{+}-x\omega_{-})\phi=\phi^{\prime}(\star\eta)\phi.
$$ 
$$
\star(\phi^{\prime}(\omega_{-}\wedge\omega_{+})\phi)=q^{2n}\star(\phi^{\prime}\phi(\omega_{-}\wedge\omega_{+}))=q^{2n}\phi^{\prime}\phi\star(\omega_{-}\wedge\omega_{+})=q^{2n}\phi^{\prime}\{\star(\omega_{-}\wedge\omega_{+})\}\phi,
$$
where the last equality is evident, since $\star(\omega_{-}\wedge\omega_{+})\in\,\IC$.
\end{proof}
\end{lemm}

\subsection{A Laplacian operator on $\Asq$}

Using the 2D differential calculus on the Podle\'s sphere $\sq$ and the $\star$-Hodge structure on $\Omega(\sq)$ it is natural to define a Laplacian operator $\Box_{\sq}:\Asq\to\Asq$ as $\Box_{\sq} f=\star\dd\star\dd f$ on any $f\in\,\Asq$. An explicit computation using the properties of the exterior algebra $\Omega(\sq)$ represented in proposition \ref{2dsph} gives:
\begin{align}
\dd f&=(X_{+}\lt f)\omega_{+}+(X_{-}\lt f)\omega_{-}, \nn \\
\star\dd f&=-i\alpha^{\prime\prime}[\nu(X_{+}\lt f)\omega_{+}-q^{-2}\beta(X_{-}\lt f)\omega_{-}], \nn \\  
\dd\star\dd f&=-i\alpha^{\prime\prime}[\nu X_{-}X_{+}+\beta X_{+}X_{-}]\lt f\,(\omega_{-}\wedge\omega_{+}), \nn \\
\star\dd\star\dd f&=-[\nu X_{-}X_{+}+\beta X_{+}X_{-}]\lt f.
\label{la2d}
\end{align}
The relation \eqref{libu} shows that such a Laplacian operator can be seen as an operator $\Box_{\sq}:\cl_{0}^{(0)}\to\cl_{0}^{(0)}$. In particular, from \eqref{la3}, the Laplacian $\Box_{\sq}$ is the restriction of the Laplacian $\Box_{\SU}$ to the subalgebra $\Asq\subset\ASU$. A basis of the eigenvector spaces $\cl_{0}^{(0)}=\oplus_{J\in\,\IN}V_{J}^{(0)}$ coming from \eqref{decoln} is given by elements $\phi_{0,J,l}=c^{J}a^{*J}\rt E^{l}=w_{2J:J,l}$, so that formulas \eqref{egv} drive to a spectrum of this Laplacian on $\sq$ as:
\begin{align}
\Box_{\sq} \phi_{0,J,l}&=-(q\nu+q^{-1}\beta)\{[J][J+1]\}\phi_{0,J,l} \nn \\
&=-2q\nu\{[J][J+1]\}\phi_{0,J,l}.
\label{spe2d}
\end{align}

\begin{rema}
\label{comp}
Equations \eqref{la3} and  \eqref{la2d} show that the classical relations between the Laplacians $\Box_{SU(2)}$ and $\Box_{S^{2}}$, coming from the Hodge duality associated to the metric tensor $g$ \eqref{gmet} related to the Cartan-Killing metric,  is then reproduced in the quantum formalism, in the specific realisation of the quantum Hopf bundle that has been described. The constraints among the 5 real parameters used in the analysis of the Hodge duality can be written as:
\begin{align}
&\gamma=\alpha^{\prime\prime2}\alpha^{\prime-2}, \nn \\
&\nu=|\alpha^{\prime\prime}|^{-1}, \nn \\
&\beta=q^{2}\nu.
\label{consum}
\end{align}
The parameters $\alpha^{\prime},\alpha^{\prime\prime}$ are the coefficients of the volume forms. The analysis of the classical limit of this formulation is in section \ref{limcla}. The choice:
\begin{align}
&\lim_{q\rightarrow 1}\alpha^{\prime}=-4\alpha, \nn \\
&\lim_{q\rightarrow 1}\alpha^{\prime\prime}=-2\alpha
\label{lcalpha}
\end{align}
gives \eqref{la3} and \eqref{la2d} in the classical limit. Being $\alpha$ a positive real number, it seems natural to assume $\alpha^{\prime}$ and $\alpha^{\prime\prime}$ negative real numbers. This also gives $\nu=-\alpha^{\prime\prime-1}$ from the second relation in  \eqref{consum}, so to have a Hodge duality \eqref{sHs} which is now:
\begin{align}
&\star(1)=\check{\theta}=i\alpha^{\prime\prime}\,\omega_{-}\wedge\omega_{+}, \nn \\
&\star(x\,\omega_{-})=-ix\,\omega_{-}, \nn \\
&\star(y\,\omega_{+})=iy\,\omega_{+}, \nn \\
&\star(i\omega_{-}\wedge\omega_{+})=\alpha^{\prime\prime-1},
\label{sHsb}
\end{align} 
giving, if \eqref{lcalpha} is satisfied, the Hodge duality \eqref{cHs2} in the classical limit.
\end{rema}
\bigskip
\bigskip

\section{Connections on the Hopf bundle}\label{se:conn}

The structure of a  quantum principal bundle $(\cp,\cb,\ch;\cn_{\cp},\cq_{\ch})$  with compatible differential 
calculi, given the  total space algebra $\cp$ on which   the gauge group Hopf algebra $\ch$ coacts, has been described in section \ref{QPB}. The compatibility conditions ensure the exactness of the sequence \eqref{des}:
\beq
0\,\to\,\mathcal{P}\Omega^{1}(\mathcal{B})\mathcal{P}\,\to\,
\Omega_{1}(\mathcal{P})
\,\stackrel{\sim_{\mathcal{N}_{\mathcal{P}}}} \longrightarrow\, \mathcal{P}\otimes \left(\ker\varepsilon_{\ch}/\mathcal{Q}_{\ch}\right)\,\to\,0 .
\label{des2}
\eeq
with the map $\sim_{\cn_{\cp}}$ defined via the commutative diagram \eqref{qdia}. Among the compatibility conditions, the requirement that $\Delta_{R}\cn_{\cp}\subset\cn_{\cp}\otimes\ch$ -- giving a right covariance of the differential structure on $\cp$ -- allows to extend the coaction $\Delta_{R}$ of $\ch$ on $\cp$ to a coaction of $\ch$ on 1-forms, $\Delta_{R}^{(1)}:\Omega^{1}(\cp)\to\Omega^{1}(\cp)\otimes\ch$, defining $\Delta_{R}^{(1)}\circ\dd=(\dd\otimes 1)\circ\Delta_{R}$. 

Note that $\Ad(\ker\varepsilon_{\ch})\subset(\ker\varepsilon_{\ch})\otimes\ch$. If the right ideal $\cq_{\ch}$ is $\Ad$-invariant (which is equivalent to say that the differential calculus on $\ch$ is bicovariant),  it is possible to define a right-adjoint coaction $\Ad^{(R)}:\ker\varepsilon_{\ch}/\cq_{\ch}\,\to\,\ker\varepsilon_{\ch}/\cq_{\ch}\otimes\ch$ by the commutative diagram
$$
\begin{array}{lcl}
\ker\varepsilon_{\ch} 
& \stackrel{{\pi_{\cq_{\ch}}}} {\longrightarrow} & \ker\varepsilon_{\ch}/\cq_{\ch} \\
\downarrow \Ad  &  & \downarrow \Ad^{(R)}   \\
\ker\varepsilon_{\ch}\otimes\ch &
\stackrel{{\pi_{\cq_{\ch}}}\otimes\id}{\longrightarrow}&
(\ker\varepsilon_{\ch}/\mathcal{Q}_{\ch})\otimes\ch
\end{array}
$$
Together with the right coaction $\Delta_{R}$ of $\ch$ on $\cp$, such a right-adjoint coaction $\Ad^{(R)}$ allows to define a right coaction $\Delta_{R}^{(\Ad)}$ of $\ch$ on $\cp\otimes\ker\varepsilon_{\ch}/\cq_{\ch}$ as a coaction of a Hopf algebra on the tensor product of its comodules. This coaction is explicitly given by the relation:
\beq
\Delta_{R}^{(\Ad)}(p\otimes\pi_{\cq_{\ch}}(h))=p_{(0)}\otimes\pi_{\cq_{\ch}}(h_{(2)})\otimes p_{(1)}(Sh_{(1)})h_{(3)},
\label{Adre}
\eeq
adopting the Sweedler notation for the coaction  as $\Delta_{R}(p)=p_{(0)}\otimes p_{(1)}$.

It is now possible to define a connection  
on the quantum principal bundle as a right invariant splitting of the sequence \eqref{des2}. Given a left $\cp$-linear map $\sigma:\cp\otimes(\ker\varepsilon_{\ch}/\cq_{\ch})\to\Omega^{1}(\cp)$ such that
\begin{align}
&\Delta_{R}^{(1)}\circ\sigma=(\sigma\otimes id)\Delta^{(\Ad)}_{R} ,\nn \\
&\sim_{\cn_{\cp}}\circ\sigma=id, 
\label{si}
\end{align}
then the map $\Pi:\Omega^{1}(\cp)\to\Omega^{1}(\cp)$ defined by $\Pi=\sigma\circ\sim_{\cn_{\cp}}$ is a right invariant left $\cp$-linear projection, whose kernel coincides with the horizontal forms $\mathcal{P}\Omega^{1}(\mathcal{B})\mathcal{P}$: 
\begin{align}
\Pi^{2}&=\Pi ,\nn \\
\Pi(\mathcal{P}\Omega^{1}(\mathcal{B})\mathcal{P})&=0,\nn\\
\Delta_{R}^{(1)}\circ\Pi&=(\Pi\otimes id)\circ\Delta_{R}^{(1)}.
\label{Pi}
\end{align}
The image of the projection $\Pi$ is the set of vertical 1-forms of the principal bundle. A connection on a principal bundle can also be written in terms of a connection 1-form, which is a map $\omega:\ch\to\Omega^{1}(\cp)$. Given a right invariant splitting $\sigma$ of the exact sequence \eqref{des2},  define the connection 1-form as $\omega(h)=\sigma(1\otimes\pi_{\cq_{\ch}}(h-\varepsilon_{\ch}(h)))$ on $h\in\,\ch$. Such a  connection 1-form has the following properties:
\begin{align}
\omega(\cq_{\ch})&=0, \nn \\
\sim_{\cn_{\cp}}(\omega(h))&=1\otimes\pi_{\cq_{\ch}}(h-\varepsilon_{\ch}(h))\qquad\forall\,h\in\,\ch,\nn \\
\Delta_{R}^{(1)}\circ\omega&=(\omega\otimes\id)\circ \Ad,\nn \\
\Pi(\dd p)&=(id\otimes\omega)\Delta_{R}(p)\qquad\forall\, p\in\,\cp.
\label{ome} 
\end{align}
 Conversely if $\omega$ is a linear map $\ker\varepsilon_{\ch}\to\Omega^{1}(\cp)$ that satisfies the first three conditions in \eqref{ome}, then there exists a unique connection on the principal bundle, such that $\omega$ is its connection 1-form. In this case, the splitting of the sequence \eqref{des2} is given by:
\beq
\sigma(p\otimes[h])=p\omega([h])
\label{siom}
\eeq
with $[h]$ in $\ker\varepsilon_{\ch}/\cq_{\ch}$, while the projection $\Pi$ is given by:
\beq
\Pi=m\circ(id\otimes\omega)\circ\sim_{\cn_{\cp}}
\label{piom}
\eeq
The general proof of these results is in \cite{bm93}. This  section explicitly describes the connections
on the quantum Hopf bundle with the compatible differential calculi presented in sections \ref{se:lcc} and \ref{se:csg}.

\subsection{Vertical subspaces on the quantum Hopf bundle}
The right coaction $\Delta_{R}^{(1)}:\Omega^{1}(\SU)\to\Omega^{1}(\SU)\otimes\ca(\U(1))$ of the gauge group algebra $\ca(\U(1))$ on the set of 1-forms on the total space algebra of the bundle, whose consistency is allowed by the   compatibility conditions between the 3D left covariant calculus on $\ASU$ and the 1D bicovariant calculus on $\ca(\U(1))$, gives:  
\begin{align}
&\Delta_{R}^{(1)}\omega_{z}=\omega_{z}\otimes 1, \nn \\
&\Delta_{R}^{(1)}\omega_{\pm}=\omega_{\pm}\otimes z^{\pm2}.
\label{rifo}
\end{align}
From the analysis on the 1D calculus on $\ca(\U(1))$ performed in section \ref{se:ccc} and the result of lemma \ref{pu}, a connection on the quantum Hopf bundle is given via a splitting map $\sigma:\ASU\otimes(\ker\varepsilon_{\U(1)}/\cq_{\U(1)})\to\Omega^{1}(\SU)$, which can be defined recalling the isomorphism $\tilde{\lambda}:\ker\varepsilon_{\U(1)}/\cq_{\U(1)}\to\IC$. Given $w\in\,\IC$  set:
\beq
\sigma(1\otimes w)=\sigma(w\otimes 1)=w(\omega_{z}+U\omega_{+}+V\omega_{-});
\label{si3i}
\eeq
and extend by the requirement of left $\ASU$-linearity, so to have: 
\begin{align}
\sigma(1\otimes[\varphi(j)])&=q^{-2j}(\omega_{z}+U\omega_{+}+V\omega_{-}), \nn \\
\sigma(\phi\otimes[\varphi(j)])&=q^{-2j}\phi(\omega_{z}+U\omega_{+}+V\omega_{-}),
\label{si3}
\end{align}
where $\phi\in\,\ASU$ and the requirement of right covariance \eqref{si} selects -- from \eqref{rifo} --  $U\in\,\cl^{(0)}_{2}$ and $V\in\,\cl^{(0)}_{-2}$. The projection $\Pi$ associated to this connection is easily seen to be:
\begin{align}
\Pi(\omega_{\pm})&=\sigma(\sim_{\cn_{\SU}}(\omega_{\pm}))=0 , \nn \\
\Pi(\omega_{z})&=\sigma(\sim_{\cn_{\SU}}(\omega_{z}))=\sigma(1\otimes[\varphi(0)])=\omega_{z}+U\omega_{+}+V\omega_{-} .
\label{Pi3}
\end{align}
In this expression the 1-forms $\omega_{\pm}$ are recovered as horizontal \eqref{omho}, a notion  depending only on the compatibility conditions between  the differential calculi, while  a choice of a connection is equivalent to the choice of the vertical part of $\Omega^{1}(\SU)$. The set of connections for the quantum Hopf bundle corresponds to the set of the possible choices of 1-forms on the basis of the bundle as  $\mathrm{a}=U\omega_{+}+V\omega_{-}\in\, \Omega^{1}(\sq)$, so that the second line in \eqref{Pi3} can be written as
\beq
\Pi(\omega_{z})=\omega_{z}+\mathrm{a}.
\label{Pi3b}
\eeq
The connection one form \eqref{ome} $\omega:U(1)\to\Omega^{1}(\SU)$ is given by:
\begin{align}
\omega(z^{j})&=\sigma(1\otimes[z^{j}-1]) \nn \\
&=\left(\frac{1-q^{-2j}}{1-q^{-2}}\right)(\omega_{z}+U\omega_{+}+V\omega_{-})
=\left(\frac{1-q^{-2j}}{1-q^{-2}}\right)(\omega_{z}+\mathrm{a}).\label{ome3}
\end{align}
Given the projection $\Pi$ and the connection 1-form $\omega$, it is possible to compute the lhs and the rhs of the last line in \eqref{ome}.  On the basis of left invariant differential forms and  using the explicit form of the quantum vector fields in \eqref{Xq}, with $\phi\in\,\cl^{(0)}_{j}$ one has:
\begin{align}
\Pi(\dd\phi)=\Pi((X_{j}\lt\phi)\omega_{j})&=(X_{j}\lt\phi)\Pi(\omega_{j}) \nn \\
&=\left(\frac{1-q^{2j}}{1-q^{-2}}\right)\phi(\omega_{z}+U\omega_{+}+V\omega_{-});
\label{pme}
\end{align}
and also: 
\begin{align}
(id\otimes\omega)\Delta_{R}(\phi)&=(id\otimes\omega)(\phi\otimes z^{-j}) \nn \\
&=\left(\frac{1-q^{2j}}{1-q^{-2}}\right)\phi(\omega_{z}+U\omega_{+}+V\omega_{-})=\Pi(\dd\phi).
\label{pmer}
\end{align}
The monopole connection corresponds to the choice $U=V=0\Leftrightarrow\mathrm{a}=0$, so to have $\Pi_{0}(\omega_{z})=\omega_{z}$ and the monopole connection 1-form $\omega_{0}(z^{j})=[(1-q^{-2j})/(1-q^{-2})]\omega_{z}$ \cite{brp,csta}. 
With a connection, one has the notion of covariant derivative $D:\ASU\to\Omega^{1}(\ASU)$ of equivariant maps. Given $\phi\in\,\cl^{(0)}_{n}$, define
\beq
D\phi=(1-\Pi)\dd\phi.
\label{Dphi}
\eeq
The covariant derivative $D\phi$ is clearly an horizontal 1-form: the adjective "covariant" refers to the behaviour under the coaction of the gauge group algebra, as one directly \eqref{clnc} shows that:
\beq
\Delta_{R}\phi=\phi\otimes z^{-j}\qquad\Leftrightarrow\qquad\Delta_{R}^{(1)}(D\phi)=D\phi\otimes z^{-j},
\label{Dcp}
\eeq
from the right invariance \eqref{Pi} of the projection $\Pi$. In terms of the connection 1-form the covariant derivative can be written, using  \eqref{pmer}, as :
\begin{align}
D\phi=(1-\Pi)\dd\phi&=\dd \phi-\Pi(\dd\phi)\nn \\ 
&=\dd \phi-\phi\wedge\omega(z^{-j})
\label{Dom}
\end{align}
on a $\phi\in\cl^{(0)}_{j}$. It is then immediate to recover that, for any $f\in\,\cl^{(0)}_{0}\simeq\Asq$, one has $Df=\dd f$. 
\begin{rema}
\label{shor}
Given any $\phi\in\,\cl_{n}^{(0)}$, from \eqref{Dom} and \eqref{Pi3b}, the covariant derivative can be written as:
$$
D\phi=\{(X_{+}\lt\phi)-(X_{z}\lt\phi)U\}\omega_{+}+\{(X_{-}\lt\phi)-(X_{z}\lt\phi)V\}\omega_{-}.
$$
It is an easy computation using the $\ASU$-bimodule properties \eqref{bi1}  of $\Omega^{1}(\SU)$ to prove that $D\phi\simeq\Omega^{1}(\sq)\cdot\ASU$ for any connection represented by $\mathrm{a}\in\,\Omega^{1}(\sq)$. This means that any connection on this quantum Hopf bundle is a strong connection, following the analysis in \cite{haj96}. 
\end{rema}

\subsection{Covariant derivative on the associated line bundles}
A covariant derivative, or a  connection,  on the left $\Asq$-module $\zce_{n}$ is a $\IC$-linear map
\beq
\nabla:\Omega^{k}(\sq)\otimes_{\Asq}\zce_{n}\to\Omega^{k+1}(\sq)\otimes_{\Asq}\zce_{n},  
\label{cdevl}
\eeq
defined for any $k\geq0$ and satisfying a left Leibniz rule:
$$
\nabla(\alpha\bra{\sigma})=\dd\alpha\wedge\bsigma+(-1)^{m}\alpha\wedge(\nabla\bsigma)
$$
for any $\alpha\in\,\Omega^{m}(\sq)$ and $\bsigma\in\,\Omega^{k}(\sq)\otimes_{\Asq}\zce_{n}$. A connection is completely determined by its restriction $\nabla:\zce_{n}\to\Omega^{1}(\sq)\otimes_{\Asq}\zce_{n}$ and then extended by the Leibniz rule. Connections always exist on projective modules: the canonical (Levi-Civita, or Grassmann) connection on a left projective $\Asq$-module $\zce_{n}$ is given as
\beq
\nabla_{0}\bsigma=(\dd\bsigma)\qpp;
\label{grass}
\eeq
the space $C(\zce_{n})$ of all connections on $\zce_{n}$ is an affine space modelled on $\mathrm{Hom}_{\Asq}(\zce_{n},\zce_{n}\otimes_{\Asq}\Omega^{1}(\sq))$, so that any connection can be written as:   
\beq
\nabla\bsigma=(\dd\bsigma)\qpp+(-1)^{k}\bsigma\wedge \An
\label{gapo}
\eeq
with $\bsigma\in\,\Omega^{k}(\sq)\otimes_{\Asq}\zce_{n}$ and $\An\in{\bb M}_{\mn+1}\otimes_{\Asq}\Omega^{1}(\sq)$ -- which is called the gauge potential of the connection $\nabla$ -- subject to the condition $\An=\An\qpp=\qpp \An$. The composition 
$$
\nabla^{2}=\nabla\circ\nabla\quad:\quad\Omega^{k}(\sq)\otimes_{\Asq}\zce_{n}\to\Omega^{k+2}(\sq)\otimes_{\Asq}\zce_{n}
$$ 
is $\Omega(\sq)$-linear. This map can be explicitly calculated: given $\bsigma\in\,\Omega^{k}(\sq)\otimes_{\Asq}\zce_{n}$, from \eqref{gapo}  one has 
\begin{align}
\nabla^{2}\bsigma&=\dd(\nabla\bsigma)\qpp+(-1)^{k+1}(\nabla\bsigma)\wedge \An \nn \\
&=\dd\{(\dd\sigma)\qpp+(-1)^{k}\bsigma\wedge \An\}\qpp+(-1)^{k+1}\{(\dd\bsigma)\qpp+(-1)^{k}\bsigma\wedge \An\}\wedge \An \nn \\
&=\dd\{(\dd\bsigma)\qpp\}\qpp+(-1)^{k}(\dd\bsigma\wedge \An)\qpp+(\bsigma\wedge\dd \An)\qpp\nn \\
&\qquad\qquad\qquad+(-1)^{k+1}\{(\dd\bsigma)\qpp\wedge \An\}-\bsigma\wedge \An\wedge \An \nn \\
&=\bsigma\{-(\dd\qpp\wedge\dd\qpp)\qpp+(\dd \An)\qpp-\An\wedge \An\}.
\label{nab2}
\end{align}
The restriction of the map $\nabla^{2}$ to $\zce_{n}$, seen as an element in $\Omega^{2}(\sq)\otimes_{\Asq}\zce_{n}$,  is the curvature $F_{\nabla}$ of the given connection.  

The left $\Asq$-module isomorphism between $\cl_{n}^{(0)}$ and $\ce_{n}^{(0)}$ described in proposition \ref{isoeqsec} allows for the definition of an hermitian structure on each projective left module $\ce_{n}^{(0)}$, $\{~;~\}:\ce_{n}^{(0)}\times\ce_{n}^{(0)}\to\Asq$ given as:
\beq
\{\bsigma_{\phi};\bsigma_{\phi^{\prime}}\}=\phi\phi^{\prime*},
\label{herms}
\eeq
with $\phi,\phi^{\prime}\in\,\cl_{n}^{(0)}$. Such an hermitian structure satisfies the relations:
\begin{align*}
&\{f\bsigma_{\phi};f^{\prime}\bsigma_{\phi^{\prime}}\}=f\phi(f^{\prime}\phi^{\prime})^*, \\
&\{\bsigma_{\phi},\bsigma_{\phi}\}\geq0, \quad
\{\bsigma_{\phi},\bsigma_{\phi}\}=0\,\Leftrightarrow\,\bsigma=0.
\end{align*}
  
\bigskip

The left $\Asq$-module isomorphism between $\cl^{(0)}_{n}$ and $\zce_{n}$ also enables to relate the concept of connection on the quantum Hopf bundle to that of covariant derivative on the associated line bundles.  As first step, define the  $\Asq$-bimodule:
\beq
\cl_{n}^{(1)}=\{\phi\in\,\Omega^{1}_{\mathrm{hor}}(\SU)\simeq\ASU\Omega^{1}(\sq)\ASU\,:\,\Delta_{R}^{(1)}\phi=\phi\otimes z^{-n}\}
\label{ln1}
\eeq
and introduce the notations:
$$
\ce_{n}^{(k)}=\Omega^{k}(\sq)\otimes_{\Asq}\zce_{n}.
$$
The maps:
\begin{align}
 &\cl_{n}^{(1)}\stackrel{\simeq}{\longrightarrow}\ce^{(1)}_{n}  \,:\qquad\phi\mapsto\bsigma_{\phi}=\phi\bra{\Psin}, \nn \\
 &\ce^{(1)}_{n}\stackrel{\simeq}{\longrightarrow}\cl_{n}^{(1)}\,:\qquad \bsigma\mapsto\phi=\hs{\sigma}{\Psin}
\label{isoeqsec1}
\end{align}
give  left $\Asq$-module isomorphisms (in this notation the explicit dependence on $\bra{f}\in\,\Asq^{\mn+1}$ as in proposition \ref{isoeqsec} has been dropped). Via this isomorphism, any connection on the quantum Hopf bundle -- represented by a projection $\Pi$ \eqref{Pi3} or by a connection 1-form \eqref{ome3} --  induces a gauge potential $\An$ on any associated line bundle $\zce_{n}$. 
\begin{prop}
\label{pro:naD}
Given the left $\Asq$-isomorphism $\cl_{n}^{(0)}\simeq\zce_{n}$ described in proposition 
 \ref{isoeqsec}, as well as the analogue left $\Asq$-module isomorphism $\cl_{n}^{(1)}\simeq\ce_{n}^{(1)}$ described in \eqref{isoeqsec1}, there is an equivalence between the set of connections on the quantum Hopf bundle via a projection $\Pi$ in $\Omega(\SU)$ as in \eqref{Pi3}, and the set of covariant derivative $\nabla\in\,C(\zce_{n})$ on any associated line bundle. 
With  $\phi\in\,\cl_{n}^{(0)}$ so that $\bra{\sigma_{\phi}}=\phi\bra{\Psin}\in\,\zce_{n}$, the equivalence is given by $D\phi=(\nabla\bra{\sigma}_{\phi})\ket{\Psin}$.
\begin{proof}
 Choose $\phi\in\,\cl^{(0)}_{n}$, so to have $\sigma_{\phi}=\phi\bra{\Psin}$ and from the definition in \eqref{gapo} express a covariant derivative on $\zce_{n}$ via a gauge potential as:
\begin{align}
\nabla\bsigma_{\phi}&=\dd\left(\phi\bra{\Psin}\right)\ket{\Psin}\bra{\Psin}+\phi\bra{\Psin}\An \label{nasi1} \\
&=\{\dd\phi-\phi[\hs{\Psin}{\dd\Psin}-\bra{\Psin}\An\ket{\Psin}]\}\bra{\Psin} \label{nasi2}
\end{align}
since $\An=\An\qpp$. On the other hand, being $\phi\in\,\cl^{(0)}_{n}$ one has:
\begin{align*}
D\phi=(1-\Pi)\dd\phi&=\dd\phi-(X_{z}\lt\phi)\Pi(\omega_{z})\nn \\
&=\dd\phi-\left(\frac{1-q^{2n}}{1-q^{-2}}\right)\phi\,\Pi(\omega_{z}), 
\end{align*}
with $D\phi\in\,\cl_{n}^{(1)}$ from \eqref{Dcp}. 
By the isomorphism \eqref{isoeqsec1}, equating $D\phi=(\nabla\bsigma_{\phi})\ket{\Psin}$ defines the gauge potential $\An$ as:
\begin{align}
\hs{\Psin}{\dd\Psin}-\bra{\Psin}\An\ket{\Psin}&=\frac{1-q^{2n}}{1-q^{-2}}(\omega_{z}+U\omega_{+}+V\omega_{-})\nn \\
&=\frac{1-q^{2n}}{1-q^{-2}}(\omega_{z}+\mathrm{a})=\omega(z^{-n}):
\label{omzA}
\end{align}
an explicit calculation shows that $\hs{\Psin}{\dd\Psin}=[(1-q^{2n})/(1-q^{-2})]\omega_{z}$, so the previous expression becomes:
\beq 
\bra{\Psin}\An\ket{\Psin}=-\frac{1-q^{2n}}{1-q^{-2}}(U\omega_{+}+V\omega_{-}), 
\label{eqA}
\eeq
which is solved by
\begin{align}
\An&=-\frac{1-q^{2n}}{1-q^{-2}}\ket{\Psin}(U\omega_{+}+V\omega_{-})\bra{\Psin}\nn \\
&=-\frac{1-q^{2n}}{1-q^{-2}}\ket{\Psin}\A\bra{\Psin}.
\label{gapo3}
\end{align}
This solution is unique. Being the set of connection an affine space, any different gauge potential,  solution of equation \eqref{eqA}, should be $\check{\mathrm{A}}^{(n)}=\An+\Apn$ where $\An$ is given in \eqref{gapo3} and $\Apn$ must satisfy $\bra{\Psin}\Apn\ket{\Psin}=0$, with $\Apn=\qpp \Apn \qpp=\qpp \Apn=\Apn\qpp=\Apn$. One directly has:
\begin{align*}
&\bra{\Psin}\Apn\ket{\Psin}=0\\
&\qquad\qquad\,\Rightarrow\,\, 0=\ket{\Psin}\bra{\Psin}\Apn \ket{\Psin}\bra{\Psin}=\qpp \Apn\qpp=\Apn.
\end{align*}
The complete equivalence claimed in the proposition comes by \eqref{eqA}, which gives for any gauge potential $\An$ a 1-form $\mathrm{a}\in\,\Omega^{1}(\sq)$, suitable to define a connection as in \eqref{Pi3b}.
\end{proof}
\end{prop}

The form of the gauge potential \eqref{gapo3} shows that the monopole connection $\Pi_{0}(\omega_{z})=\omega_{z}$ corresponds to the Grassmann, or canonical covariant derivative $\nabla_{0}\bsigma=(\dd\bsigma)\qpp$ on the line bundles $\zce_{n}$, having $\An=0$ for any $n\in\,\IZ$.  
A connection on the quantum Hopf bundle is defined compatible with the hermitian structure \eqref{herms} on each module of sections of the associated line bundle if 
$$
\dd\{\bsigma_{\phi};\bsigma_{\phi^{\prime}}\}=\{\nabla\bsigma_{\phi};\bsigma_{\phi^{\prime}}\}+\{\bsigma_{\phi};\nabla\bsigma_{\phi^{\prime}}\}.
$$
It is easy to compute that this condition amounts to have a connection \eqref{Pi3b} satisfying the condition 
$$
\mathrm{a}^{*}=-\mathrm{a}.
$$  

\bigskip

The compatibility between the differential calculi allows to  extend the concept of right coaction of the gauge group algebra on the whole exterior algebra $\Omega(\SU)$, introducing a right coaction $\Delta_{R}^{(k)}:\Omega^{k}(\SU)\to\Omega^{k}(\SU)\otimes\ca(\U(1))$ by induction as  
\beq
\Delta_{R}^{(k)}\circ\dd=(\dd\otimes\id)\circ\Delta_{R}^{(k-1)}. 
\label{dre}
\eeq
It becomes now  natural to define the  $\Asq$-bimodule:
\beq
\cl_{n}^{(2)}=\{\phi\in\,\ASU\Omega^{2}(\sq)\ASU\,:\,\Delta_{R}^{(2)}\phi=\phi\otimes z^{-n}\};
\label{ln2}
\eeq
so that the  maps:
\begin{align}
 &\cl_{n}^{(2)}\stackrel{\simeq}{\longrightarrow}\ce^{(2)}_{n}\,:\qquad\phi\mapsto\bsigma_{\phi}=\phi\bra{\Psin}, \nn \\
 &\ce^{(2)}_{n}\stackrel{\simeq}{\longrightarrow}\cl_{n}^{(2)}\,:\qquad \bsigma\mapsto\phi=\hs{\sigma}{\Psin}
\label{isoeqsec2}
\end{align}
are  left $\Asq$-module isomorphisms, generalising the isomorphisms given in proposition \ref{isoeqsec} and in \eqref{isoeqsec1}. In the formulation of \cite{bm93}, the elements in $\cl_{n}^{(k)}$ are strongly tensorial forms.

Recall that the covariant derivative $\nabla$ is defined in \eqref{cdevl}
as an operator $\nabla:\ce^{(k)}_{n}\to\ce^{(k+1)}_{n}$ for $k=0,1,2$, since the differential calculus on $\Asq$ is 2 dimensional; the covariant derivative $D$ has been defined by \eqref{Dphi} only on the $\Asq$-bimodule $\cl^{(0)}_{n}$, while the proposition \ref{pro:naD} shows the equivalence between $D:\cl^{(0)}_{n}\to\cl_{n}^{(1)}$ and $\nabla:\zce_{n}\to\ce^{(1)}_{n}$. The isomorphism \eqref{isoeqsec2} allows then to extend the covariant derivative to $D:\cl_{n}^{(1)}\to\cl_{n}^{(2)}$, defining:
\beq
D\phi=(\nabla\bsigma_{\phi})\ket{\Psin}
\label{D12}
\eeq
for any $\phi\in\,\cl_{n}^{(1)}$ with $\bsigma_{\phi}=\phi\bra{\Psin}\in\,\ce_{n}^{(1)}=\Omega^{1}(\sq)\otimes_{\Asq}\zce_{n}$. Such an operator can be represented in terms of the connection \eqref{ome3} 1-form $\omega$.  From the Leibniz rule one has: 
$$
\dd(\phi\bra{\Psin})=(\dd\phi)\bra{\Psin}+(-1)^{k}\phi\,\dd\bra{\Psin},
$$
with $\phi\in\,\cl_{n}^{(k)}$. This identity gives the next proposition.

\begin{prop}
\label{Dhom}
Given $\phi\in\cl_{n}^{(1)}$ , so that 
$\bsigma_{\phi}=\phi\bra{\Psin}\in\,\ce^{(1)}_{n}$, the action of the operator $D:\cl_{n}^{(1)}\to\cl_{n}^{(2)}$ defined by \eqref{D12} can be written as:
\beq
D\phi=\dd\phi+\phi\wedge\omega(z^{-n})
\eeq
\begin{proof}
The proposition is proved by a direct computation. Start from $\phi\in\,\cl_{n}^{(1)}$, so that from\eqref{gapo} one has $\nabla\bsigma_{\phi}=(\dd\bsigma_{\phi})\qpp-\bsigma_{\phi}\wedge\An$, so that :
\begin{align}
D\phi&=(\nabla\bsigma_{\phi})\ket{\Psin}\nn \\
&=(\dd\bsigma_{\phi})\ket{\Psin}-\bsigma_{\phi}\wedge\An\ket{\Psin}
\nn \\
& =\dd(\phi\bra{\Psin})\ket{\Psin}-\phi\wedge\bra{\Psin}\An\ket{\Psin} \nn \\
&=\dd\phi+\phi\wedge\hs{\Psin}{\dd\Psin}-\phi\wedge\bra{\Psin}\An\ket{\Psin}=\dd\phi+\phi\wedge\omega(z^{-n}),
\end{align}
where the last equality comes from \eqref{omzA}, expressing the gauge potential $\An$ in terms of the connection 1-form $\omega$. 
\end{proof}
\end{prop}

\bigskip

To give the curvature $F_{\nabla}$ of the given connection \eqref{nab2} a more explicit form, one can make use of two further relations. The first one, involving the projectors $\qpp$ only, comes from \cite{lrz}, while the second is proved again by direct calculation. 

\begin{lemm}
\label{lepro}
Let $\qpp$ denote the projection given in \eqref{dP}. With
the  2D calculus on $\sq$ of section~\ref{se:cals2} one finds: 
\begin{align*}
& \dd \qpp \wedge\, \dd \qpp\,  \qpp\, = -q^{-n-1} [n]
~\qpp\, \omega_{+}\wedge\omega_{-} , \\
& \qpp\,\dd \qpp \wedge\, \dd \qpp\, = -q^{-n-1} [n]
 ~\qpp\, \omega_{+}\wedge\omega_{-} 
.
\end{align*}
\end{lemm}

\begin{lemm}
\label{ledA}
Given for any $n\in\,\IZ$ the projectors $\qpp$ as in \eqref{dP} and the expression of the gauge potential $\An$ as in \eqref{gapo3}, one has:
\beq
\qpp\dd \An\qpp=-\left(\frac{1-q^{2n}}{1-q^{-2}}\right)\ket{\Psin}\dd(U\omega_{+}+V\omega_{-})\bra{\Psin}.
\label{prid}
\eeq
\begin{proof}
Setting 
$$
\A^{(n)}=\bra{\Psin}\An\ket{\Psin}=-\{(1-q^{2n})/(1-q^{-2})\}(U\omega_{+}+V\omega_{-})=-\frac{1-q^{2n}}{1-q^{-2}}\mathrm{a},
$$ 
the expression \eqref{prid} can be written as the sum of  three terms, from the Leibniz rule satisfied by the exterior derivation $\dd$:
\begin{align}
&\qpp\dd\An\qpp\nn \\
&\qquad=\left[\ket{\Psin}\hs{\Psin}{\dd\Psin}\A^{(n)}\bra{\Psin}\right]+\left[\ket{\Psin}(\dd\A^{(n)})\bra{\Psin}\right]-
\left[\ket{\Psin}\A^{(n)}\hs{\dd\Psin}{\Psin}\bra{\Psin}\right] \nn \\
&\qquad=\left(\frac{1-q^{2n}}{1-q^{-2}}\right)\ket{\Psin}\left[\omega_{z}\wedge\A^{(n)}+\A^{(n)}\wedge\omega_{z}\right]\bra{\Psin}+\left[\ket{\Psin}\dd\A^{(n)}\bra{\Psin}\right], \nn
\end{align}
where the second equality comes from the identities $\hs{\Psin}{\dd\Psin}=-\hs{\dd\Psin}{\Psin}=\{(1-q^{2n})/(1-q^{-2})\}\omega_{z}$, while the $\ASU$-bimodule relations \eqref{bi1} of 1-forms in $\Omega^{1}(\SU)$, as well as  commutation relations among them \eqref{commc3}, give:
$$
\omega_{z}\wedge(U\omega_{+}+V\omega_{-})=q^{4}U\omega_{z}\wedge\omega_{+}+q^{-4}V\omega_{z}\wedge\omega_{-}=-(U\omega_{+}+V\omega_{-})\wedge\omega_{z}, 
$$
so that $\omega_{z}\wedge\A^{(n)}+\A^{(n)}\wedge\omega_{z}=0$ and the identity claimed in \eqref{prid} is verified. 
\end{proof}
\end{lemm}
\begin{rema}
\label{reoza}
The identity $\omega_{z}\wedge\A^{(n)}+\A^{(n)}\wedge\omega_{z}=0$ also shows that the 1-form $\omega_{z}$ anti-commutes with every 1-form  in $\Omega^{1}(\sq)$.  
\end{rema}

\begin{prop}
\label{Fna}
Given the covariant derivative $\nabla:\ce^{(k)}_{n}\to\ce^{(k+1)}_{n}$
from  \eqref{gapo} with a gauge potential \eqref{gapo3} $\An=-(1-q^{2n})(1-q^{-2})^{-1}\ket{\Psin}\A\bra{\Psin}$, the operator $\nabla^{2}:\ce^{(0)}_{n}\to\ce^{(2)}_{n}$
can be written as:
\beq
\nabla^{2}\bsigma=\bsigma\wedge F_{\nabla}=-\bsigma\wedge\{\ket{\Psin}q^{n+1}[n](\omega_{-}\wedge\omega_{+}-\dd\A+q^{n+1}[n]\A\wedge\A)\bra{\Psin}\}.
\label{Fa}
\eeq
\begin{proof}
From the general expression \eqref{nab2},  the action of the operator $\nabla^{2}$ on a $\bsigma\in\,\zce_{n}$ is linear, and given by the sum of three terms. The first one, recalling the result of the lemma \ref{lepro} and the commutation rules \eqref{bi1} and \eqref{commc3}, is:
\begin{align}
-(\dd\qpp\wedge\dd\qpp)\qpp&=q^{-n-1}[n]\qpp\omega_{+}\wedge\omega_{-}\nn \\
&=-q^{1-n}[n]\ket{\Psin}\bra{\Psin}\omega_{-}\wedge\omega_{+}\nn \\
&=-q^{n+1}[n]\ket{\Psin}\omega_{-}\wedge\omega_{+}\bra{\Psin}.
\label{Fn1}
\end{align}
Since one has $\bsigma\qpp=\bsigma$, being  elements in the projective modules $\zce_{n}$, the other two terms in \eqref{gapo} are:
\begin{align}
\qpp\dd\An\qpp&=-\left(\frac{1-q^{2n}}{1-q^{-2}}\right)\ket{\Psin}\dd\A\bra{\Psin}\nn \\
&=q^{n+1}[n]\ket{\Psin}\dd\A\bra{\Psin}, \nn \\
-\An\wedge\An&=-\left(\frac{1-q^{2n}}{1-q^{-2}}\right)^{2}\ket{\Psin}\A\wedge\A\bra{\Psin}\nn \\
&=-q^{2(n+1)}[n]^{2}\ket{\Psin}\A\wedge\A\bra{\Psin}\nn.
\end{align}
The sum of these three lines gives the curvature $F_{\nabla}\in\,{\bb M}_{\mn+1}\otimes_{\Asq}\Omega^{2}(\sq)$ the expression:
\beq
F_{\nabla}=-\ket{\Psin}q^{n+1}[n](\omega_{-}\wedge\omega_{+}-\dd\A+q^{n+1}[n]\A\wedge\A)\bra{\Psin}.
\label{Fas}
\eeq
\end{proof}
\end{prop}

The  isomorphism \eqref{isoeqsec2}  allows to formulate the curvature  as a linear map $D^{2}:\cl^{(0)}_{n}\to\cl_{n}^{(2)}$, defined by:
\beq
D^{2}\phi=(\nabla^{2}\bsigma_{\phi})\ket{\Psin}
\label{curD}
\eeq   
for a given  $\phi=\hs{\sigma}{\Psin}$.  This operator can also be written in terms of the connection 1-form $\omega$.

\begin{prop}
\label{D2ome}
The operator $D^{2}:\cl^{(0)}_{n}\to\cl_{n}^{(2)}$ defined in \eqref{curD} can be written as
\beq
D^{2}\phi=-\phi\wedge\{\dd\omega(z^{-n})+\omega(z^{-n})\wedge\omega(z^{-n})\}=\phi\wedge \left(\bra{\Psin}F_{\nabla}\ket{\Psin}\right)
\label{D2om}
\eeq
on any $\phi\in\,\cl^{(0)}_{n}$.
\begin{proof}
The proof is a direct application of the result in propositions \ref{Dom} and \ref{Dhom}. It is $D\phi=\dd\phi-\phi\wedge\omega(z^{-n})$ with $\phi\in\,\cl^{(0)}_{n}$, so that:
\begin{align}
D^{2}\phi&=\dd(D\phi)+(D\phi)\wedge\omega(z^{-n})\nn \\
&=-\dd(\phi\wedge\omega(z^{-n}))+(\dd\phi-\phi\wedge\omega(z^{-n}))\wedge\omega(z^{-n})\nn \\
&=-\phi\wedge(\dd\omega(z^{-n})+\omega(z^{-n})\wedge\omega(z^{-n})).
\nn
\end{align}
The relation \eqref{omzA} can be rewritten as $\omega(z^{-n})=-q^{1+n}[n](\omega_{z}+\A)$, so to have:
\begin{align}
&\dd\omega(z^{-n})=-q^{1+n}[n](\dd\omega_{z}+\dd\A)=q^{1+n}[n](\omega_{-}\wedge\omega_{+}-\dd\A), \nn \\
&\omega(z^{-n})\wedge\omega(z^{-n})=\{q^{1+n}[n]\}^{2}(\omega_{z}+\A)\wedge(\omega_{z}+\A)=q^{2(1+n)}[n]^{2}\A\wedge\A, \nn
\end{align}
where the last equality in the second line comes from the remark \ref{reoza}. It  becomes then clear to recover from \eqref{Fa}
$$
D^{2}\phi=-\phi\wedge q^{1+n}[n]\{\omega_{-}\wedge\omega_{+}-\dd\A+q^{1+n}[n]\A\wedge\A\}=\phi\wedge \left(\bra{\Psin}F_{\nabla}\ket{\Psin}\right),
$$
meaning that the action of the operator $D^{2}$ can be represented by the 2-form $\left(\bra{\Psin}F_{\nabla}\ket{\Psin}\right)\in\cl_{0}^{(2)}$.
\end{proof}
\end{prop}

\begin{rema}
\label{rempi}
Recall from \eqref{Dphi} that, given $\phi\in\,\cl^{(0)}_{n}$, the covariant derivative $D:\cl^{(0)}_{n}\to\cl_{n}^{(1)}$ has been defined in terms of the projector $\Pi$ associated to the connection  as:
$$
D\phi=(1-\Pi)\dd\phi.
$$
Given the left $\Asq$-module isomorphisms $\cl_{n}^{(k)}\simeq\Omega^{k}(\sq)\otimes_{\Asq}\zce_{n}=\ce_{n}^{(k)}$, the proposition \ref{pro:naD} shows that any connection written as  a projector $\Pi$ as in \eqref{Pi3} induces a gauge potential $\An$, so to have a covariant derivative $\nabla:\Omega^{k}(\sq)\otimes_{\Asq}\ce_{n}\to\Omega^{k+1}(\sq)\otimes_{\Asq}\ce_{n}$. The operator $D$ is then extended in \eqref{D12} as  $D:\cl_{n}^{(1)}\to\cl_{n}^{(2)}$  in terms of the operator $\nabla$, without using the projector $\Pi$. This definition is perfectly consistent, but it seems natural to understand whether it is possible to define $D:\cl_{n}^{(1)}\to\cl_{n}^{(2)}$ via the projector $\Pi$, and even whether it is possible to extend the domain of such a covariant derivative operator $D$ from the set of  horizontal forms $\cl_{n}^{(k)}$ to the whole exterior algebra $\Omega(\SU)$, in analogy to the classical case \eqref{2p1}.

Given $\phi\in\cl_{n}^{(1)}$, the most natural definition of a covariant derivative seems to be:
\beq
\label{Dc}
\check{D}\phi=(1-\Pi)\dd\phi,
\eeq
with the horizontal projector $(1-\Pi)$ extended to $\Omega^{2}(\SU)$ by assuming a compatibility with the wedge product 
$$
\Omega^{2}(\SU)=\{\Omega^{1}(\SU)\otimes_{\ASU}\Omega^{1}(\SU)\}/\mathcal{S}_{\cq}
=\Omega^{1}(\SU)\wedge\Omega^{1}(\SU)
$$
so to have:
\beq
(1-\Pi)\Omega^{2}(\SU)=\{(1-\Pi)\Omega^{1}(\SU)\}\wedge\{(1-\Pi)
\Omega^{1}(\SU)\}.  
\label{num}
\eeq
It is easy to see that such a compatibility \emph{does not} exist. To be definite, consider an example. Choose $\omega_{+}\in\,\cl_{-2}^{(1)}$, so that $\dd\omega_{+}=q^{2}(1+q^{2})\omega_{z}\wedge\omega_{+}=-(1+q^{-2})\omega_{+}\wedge\omega_{z}$ by the commutation properties of the $\wedge$ product \eqref{commc3}. Compute now:
\begin{align}
 q^{2}(1+q^{2})(1-\Pi)\{\omega_{z}\wedge\omega_{+}\}&=q^{2}(1+q^{2})\{(1-\Pi)\omega_{z}\}\wedge\{(1-\Pi)\omega_{+}\}=q^{2}(1+q^{2})V\omega_{-}\wedge\omega_{+}, \nn \\
-(1+q^{-2})(1-\Pi)\{\omega_{+}\wedge\omega_{z}\}&=-(1+q^{-2})\{(1-\Pi)\omega_{+}\}\wedge\{(1-\Pi)\omega_{z}\}=(1+q^{-2})V\omega_{-}\wedge\omega_{+}, \nn 
\end{align}
The two expressions are different: the problem is that, for the given 3D calculus on $\ASU$, one has 
\beq
(1-\Pi)\mathcal{S}_{\cq}\nsubseteq\mathcal{S}_{\cq}.
\label{1mpsq}
\eeq
Consider the 6 relations \eqref{commc3} generating $\mathcal{S}_{\cq}$. An explicit calculation  shows that, from the three of them not involving $\omega_{z}$, one has:
\begin{align}
&\{(1-\Pi)\omega_{+}\}\wedge\{(1-\Pi)\omega_{+}\}=0, \nn \\
&\{(1-\Pi)\omega_{-}\}\wedge\{(1-\Pi)\omega_{-}\}=0, \nn \\
&\{(1-\Pi)\omega_{-}\}\wedge\{(1-\Pi)\omega_{+}\}+q^{-2}\{(1-\Pi)\omega_{+}\}\wedge\{(1-\Pi)\omega_{-}\}=0, \nn
\end{align}
while from the remaining terms:
\begin{align}
&\{(1-\Pi)\omega_{z}\}\wedge\{(1-\Pi)\omega_{-}\}+q^{4}\{(1-\Pi)\omega_{-}\}\wedge\{(1-\Pi)\omega_{z}\}=(1-q^{4})U\omega_{+}\wedge\omega_{-}, \nn \\
&\{(1-\Pi)\omega_{z}\}\wedge\{(1-\Pi)\omega_{+}\}+q^{-4}\{(1-\Pi)\omega_{+}\}\wedge\{(1-\Pi)\omega_{z}\}=(1-q^{-4})V\omega_{-}\wedge\omega_{+}, \nn \\
&\{(1-\Pi)\omega_{z}\}\wedge\{(1-\Pi)\omega_{z}\}=\A\wedge\A. \nn 
\end{align}
These computations show that only in the case of the monopole connection -- that is $\A=0$ -- it is
\beq
(1-\Pi_{0})\mathcal{S}_{\cq}\subseteq\mathcal{S}_{\cq}:
\label{1mpsq0}
\eeq
 only in the case of the monopole connection it is consistent to set 
$$
(1-\Pi_{0})\Omega^{2}(\SU)=\{(1-\Pi_{0})\Omega^{1}(\SU)\}\wedge\{(1-\Pi_{0})
\Omega^{1}(\SU)\} 
$$
and to define 
\beq
D_{0}:\Omega^{k}(\SU)\mapsto\Omega^{k+1}(\SU), \qquad 
D_{0}\phi=(1-\Pi_{0})\dd\phi
\eeq
The operator $D_{0}$ is a 'covariant' operator: given $\phi\in\Omega^{k}(\SU)$ such that $\Delta_{R}^{(k)}\phi=\phi\otimes z^{-n}$, it is $\Delta_{R}^{(k+1)}(D_{0}\phi)=D_{0}\phi\otimes z^{-n}$, and moreover $D_{0}\phi\in\,\cl_{n}^{(k)}$: $D_{0}\phi$ is horizontal. Note that $\cl_{n}^{(3)}=\emptyset$, as the calculus on $\sq$ is 2D.  It becomes an easy computation to prove that the restriction $D_{0}:\cl_{n}^{(k)}\to\cl_{n}^{(k+1)}$ acquires the form:
\beq
D_{0}\phi=(1-\Pi_{0})\dd\phi=\dd\phi-(-1)^{k}\phi\wedge\omega_{0}(z^{-n}).
\label{horD}
\eeq
This relation is the quantum analogue of  the classical \eqref{2p2}. The classical covariant derivative of an equivariant differential form $\phi$ can be expressed in terms of the connection 1-form $\omega$ only if such $\phi$ is horizontal. In this quantum formulation, the classical 
condition that $\phi$ is horizontal and equivariant has been translated into the condition $\phi\in\cl_{n}^{(k)}$.

\end{rema}

\section{A gauged Laplacian on the quantum Hopf bundle}\label{se:gL}

With a covariant derivative $\nabla$ acting on the left $\Asq$-projective modules $\ce_{n}^{(k)}=\Omega^{k}(\sq)\otimes_{\Asq}\ce_{n}$ and the  $\star$-Hodge structure on the exterior algebra $\Omega(\sq)$ introduced in section \ref{se:hlss} it is possible to define a gauged Laplacian operator $\Box_{\nabla}:\zce_{n}\to\zce_{n}$ as:
\beq
\Box_{\nabla}\bsigma=\star\nabla\star\nabla\bsigma
\label{deBn}
\eeq
on any $\bsigma\in\,\zce_{n}$. From the left $\Asq$-linearity of the $\star$-Hodge map, and the relation \eqref{gapo}, one has:
\begin{align}
\nabla\star\nabla\bsigma&=\dd\{\star(\nabla\bsigma)\}\qpp-(\star\nabla\bsigma)\wedge\An \nn \\
&=\dd\{\star[(\dd\bsigma)\qpp]+\bsigma\wedge(\star\An)\}\qpp-\{(\star[(\dd\bsigma)\qpp]\wedge\An+\bsigma\wedge(\star\An)\wedge\An\}\nn \\
&=\dd\{\star[(\dd\bsigma)\qpp]\}\qpp +\dd\{\bsigma\wedge(\star\An)\}\qpp-\star\{(\dd\bsigma)\qpp\}\wedge\An-\bsigma\wedge(\star\An)\wedge\An
\label{1nsn}
\end{align}
The second term in the last line can be written as:
\begin{align}
\dd\{\bsigma\wedge(\star\An)\}\qpp&=\dd\bsigma\wedge(\star\An)\qpp+\bsigma\wedge\{\dd(\star\An)\}\qpp \nn \\&=
\dd\bsigma\wedge(\star\An)+\bsigma\wedge\{\dd(\star\An)\}\qpp,
\label{2nsn}
\end{align}
while the third term in \eqref{1nsn} is:
\begin{align}
-\star\{(\dd\bsigma)\qpp\}\wedge\An&=-\star(\dd\bsigma)\qpp\wedge\An \nn \\ 
&=-(\star\dd\bsigma)\wedge\An:
\label{3nsn}
\end{align}
in both the relations \eqref{2nsn} and \eqref{3nsn} the specific property of right $\Asq$-linearity of the $\star$-Hodge map has been used, namely as  $\star(\An)\qpp=\star(\An\qpp)=\star\An$ in \eqref{2nsn} and 
as $\star\{(\dd\bsigma)\qpp\}=\star(\dd\bsigma)\qpp$ in \eqref{3nsn}. Moreover, from the proposition \ref{prost2} one has $\dd\bsigma\wedge(\star\An)=-(\star\dd\bsigma)\wedge\An$, so that
\beq
\star\nabla\star\nabla\bsigma=
\star\dd\{\star(\dd\bsigma)\qpp\}\qpp -2\star\{(\star\dd\bsigma)\wedge\An\}  +\bsigma\wedge\{\star\dd\star\An\}\qpp-\bsigma\wedge\star\{(\star\An)\wedge\An\}
\label{4nsn}
\eeq
The four terms componing the gauged Laplacian can be individually studied.  
\begin{itemize}
\item
Recalling the result of lemma \ref{cala}, one has:  
\begin{align}
\star\An&=q^{n+1}[n]\star\{\ket{\Psin}\A\bra{\Psin}\}\nn \\
&=q^{n+1}[n]\ket{\Psin}(\star\A)\bra{\Psin}.
\label{sAa}
\end{align}
The fourth term in \eqref{4nsn} is, using once more the result of lemma \ref{cala} with $\bra{\Psin}\in\,\cl^{(0)}_{-n}$:
\begin{align}
-\bsigma\wedge\star\{(\star\An)\wedge\An\}&=-\bsigma\wedge q^{2(1+n)}[n]\star\{\ket{\Psin}(\star\A)\wedge\A\bra{\Psin}\}\nn \\
&=-q^{2}[n]\bsigma\wedge\ket{\Psin}\left(\star\{(\star\A)\wedge\A\}\right)\bra{\Psin}. 
\label{la4}
\end{align}
\item
From  \eqref{sAa} the third term in the expression \eqref{4nsn} of the gauged Laplacian is:
\begin{align}
\bsigma\wedge\{\star\dd\star\An\}\qpp&=\bsigma\wedge q^{1+n}[n]\star\{\dd\left(\ket{\Psin}(\star\A)\bra{\Psin}\right)\}\qpp \nn \\ &=\bsigma\wedge q^{1+n}[n]\star\{\qpp\dd\left(\ket{\Psin}(\star\An)\bra{\Psin}\right)\qpp\}.
\label{5nsn}
\end{align} 
The last term in curly bracket is, by the derivation property of $\dd$:
\begin{align}
\qpp&\dd\left(\ket{\Psin}(\star\An)\bra{\Psin}\right)\qpp= \nn \\
&=\ket{\Psin}\left(\hs{\Psin}{\dd\Psin}(\star\A)\right)\bra{\Psin}-\ket{\Psin}\left((\star\A)\hs{\dd\Psin}{\Psin}\right)\bra{\Psin}+\ket{\Psin}\left(\dd(\star\A)\right)\bra{\Psin}\nn \\
&=\ket{\Psin}\{-q^{1+n}[n]\omega_{z}\wedge(\star\A)-q^{1+n}[n](\star\A)\wedge\omega_{z}+\dd(\star\A)\}\bra{\Psin},
\label{6nsn}
\end{align}
where the last equality comes from the identity $\hs{\Psin}{\dd\Psin}=-q^{1+n}[n]\omega_{z}$. Recalling the remark \ref{reoza}, and using the commutation rules \eqref{bi1} as they were used in  \eqref{la4},  the expression \eqref{5nsn} becomes:
\begin{align}
\bsigma\wedge\{\star\dd\star\An\}\qpp&=q^{1+n}[n]\bsigma\wedge\star\{\ket{\Psin}\dd(\star\A)\bra{\Psin}\}
\nn \\ &=q^{1-n}[n]\bsigma\wedge\ket{\Psin}\{\star\dd\star\A\}\bra{\Psin}.
\label{la3bis}
\end{align}
\item
It is now straightforward to analyse the second term in the expression \eqref{4nsn} of the gauged Laplacian.  From the definition \eqref{gapo3} and the Hodge duality \eqref{sHs}, with again $\mathrm{a}=U\omega_{+}+V\omega_{-}$, $U\in\,\cl^{(0)}_{2}$ and $V\in\,\cl^{(0)}_{-2}$ :
\begin{align}
2\star\{\dd\bsigma\wedge(\star\An)\}&=2i\alpha^{\prime\prime}\nu\,q^{n+1}[n]\star\{(X_{+}\lt\bsigma)\omega_{+}\ket{\Psin}\wedge\mathrm{a}\bra{\Psin}-(X_{-}\lt\bsigma)\omega_{-}\ket{\Psin}\wedge\mathrm{a}\bra{\Psin}\} \nn \\
&=-2i\alpha^{\prime\prime}\nu\,q[n]\{(X_{+}\lt\bsigma\ket{\Psin}V\bra{\Psin}+q^{2}(X_{-}\lt\bsigma)\ket{\Psin}U\bra{\Psin}\}\star(\omega_{-}\wedge\omega_{+}) \nn \\
&=-2q[n]\{\nu(X_{+}\lt\bsigma\ket{\Psin}V\bra{\Psin}+\beta(X_{-}\lt\bsigma)\ket{\Psin}U\bra{\Psin}\}\
\label{7nsn}
\end{align}
\item
To analyse the first term in \eqref{4nsn}, which is the only one not depending on the gauge potential $\A$, start with:
\begin{align} 
\star\{(\dd\bsigma)\qpp\}&=\star\left(\{(X_{+}\lt\bsigma)\omega_{+}+(X_{-}\lt\bsigma)\omega_{-}\}\qpp\right) \nn \\
&=\star\{(X_{+}\lt\bsigma)\qpp\omega_{+}+(X_{-}\lt\bsigma)\qpp\omega_{-}\} \nn \\
&=-i\alpha^{\prime\prime}\nu\{(X_{+}\lt\bsigma)\qpp\omega_{+}-(X_{-}\lt\bsigma)\qpp\omega_{-}\} 
\label{8nsn}
\end{align} 
so to have:
\begin{align}
\dd\star\{(\dd\bsigma)\qpp\}&=-i\alpha^{\prime\prime}\nu\left(X_{-}\lt\left[\{X_{+}\lt\bsigma\}\qpp\right]\omega_{-}\wedge\omega_{+}-X_{-}\lt\left[\{X_{-}\lt\bsigma\}\qpp\right]\omega_{+}\wedge\omega_{-}\right) \nn \\
&=-i\alpha^{\prime\prime}\nu\left(X_{-}\lt\left[\{X_{+}\lt\bsigma\}\qpp\right]+q^{2}X_{+}\lt\left[\{X_{-}\lt\bsigma\}\qpp\right]\right)\omega_{-}\wedge\omega_{+} \nn \\
\star\left(\dd\star\{(\dd\bsigma)\qpp\}\right)&=-i\alpha^{\prime\prime}\left(\nu X_{-}\lt\left[\{X_{+}\lt\bsigma\}\qpp\right]+\beta X_{+}\lt\left[\{X_{-}\lt\bsigma\}\qpp\right]\right)\star(\omega_{-}\wedge\omega_{+}) \nn \\
&=-\left(\nu X_{-}\lt\left[\{X_{+}\lt\bsigma\}\qpp\right]+\beta X_{+}\lt\left[\{X_{-}\lt\bsigma\}\qpp\right]\right)
\label{9nsn}
\end{align}
 \end{itemize}
 The gauged Laplacian can be seen as an operator $\Box_{D}:\cl^{(0)}_{n}\to\cl^{(0)}_{n}$ via the equivalence between equivariant maps $\phi\in\cl^{(0)}_{n}$ and section of the associated line bundles $\sigma\in\,\zce_{n}$, represented by the isomorphism in proposition \ref{isoeqsec}:
\beq
\Box_{D}\phi=(\Box_{\nabla}\bsigma)\ket{\Psin}
\label{BD}
\eeq
on any equivariant $\phi=\hs{\sigma}{\Psin}$. The terms $(X_{\pm}\lt\bsigma)\ket{\Psin}$ in \eqref{7nsn} and\eqref{9nsn} need a specific analysis. Given the coproduct $\Delta X_{\pm}=1\otimes X_{\pm}+X_{\pm}\otimes K^{2}$, one has:
\begin{align}
(X_{\pm}\lt\bsigma)\ket{\Psin}&=(X_{\pm}\lt\{\phi\bra{\Psin}\})\ket{\Psin}\nn \\
&=\phi(X_{\pm}\lt\bra{\Psin})\ket{\Psin}+q^{-n}(X_{\pm}\lt \phi)\nn \\
&=q^{-n}(X_{\pm}\lt\phi).
\label{10nsn}
\end{align}
This last equality is clear from \eqref{lact} with $X_{+}$ and $n<0$, and with $X_{-}$ and $n>0$. In the other two cases, it is possible to apply once more the  deformed Leibniz rule to products of elements in $\ASU$, having:
\begin{align} 
q^{n}(X_{\pm}\lt\bra{\Psin})\ket{\Psin}&=X_{\pm}\lt\hs{\Psin}{\Psin}-\bra{\Psin}(X_{\pm}\lt\ket{\Psin}) \nn \\
&=X_{\pm}\lt(1)-\bra{\Psin}(X_{\pm}\lt\ket{\Psin})\nn \\
&=-\bra{\Psin}(X_{\pm}\lt\ket{\Psin})=0;
\label{11nsn}
\end{align} 
since again from \eqref{lact} one has $X_{+}\lt\ket{\Psin}=0$ with $n>0$, and $X_{-}\lt\ket{\Psin}=0$ with $n<0$. 

Recollecting the four terms from \eqref{4nsn} and making use of the relation \eqref{10nsn},  one has:
\begin{align}
-\sigma\wedge\star\{(\star\An)\wedge\An\}\ket{\Psin}&=
-q^{2}[n]\phi\wedge\star\{(\star\A)\wedge\A\}, \nn \\ 
\sigma\wedge\{\star\dd\star\An\}\ket{\Psin}&=q^{1-n}[n]\phi\wedge\{\star\dd\star\A\}, \nn \\
2\star\{\dd\sigma\wedge(\star\A)\}\ket{\Psin}&=
-2q^{1-n}[n]\left(\nu(X_{+}\lt\phi)V+\beta(X_{-}\lt\phi)U\right),
\nn \\
\left[\star\left(\dd\star\{(\dd\sigma)\qpp\}\right)\right]\ket{\Psin}&=-q^{-2n}\left(\nu X_{-}X_{+}+\beta X_{+}X_{-}\right)\lt\phi.
\label{11nsn}
\end{align}
It is clear that the gauged Laplacian operator can be completely diagonalised only if one chooses the gauge potential $\A=0$, that is if one gauges the Laplacian by the monopole connection. Such a gauged Laplacian $\Box_{D_{0}}:\cl^{(0)}_{n}\to\cl^{(0)}_{n}$ can be written as:
\beq
\Box_{D_{0}}\phi=
-q^{-2n}\left(\nu X_{-}X_{+}+\beta X_{+}X_{-}\right)\lt\phi, \qquad\mathrm{for}\,\phi\in\,\cl^{(0)}_{n}.
\label{lagc}
\eeq
The diagonalisation is straightforward, following \eqref{egv}. One has:
\begin{align}
\Box_{D_{0}}\phi_{n,J,l}&=-q^{1-n}\nu \{[J-\frac{n}{2}][J+1+\frac{n}{2}]\}-q^{-1-n}\beta
\{[J-\frac{n}{2}][J+1+\frac{n}{2}]+[n]\}
\phi_{n,J,l} \nn \\
&=-q^{1-n}\nu \{2[J-\frac{n}{2}][J+1+\frac{n}{2}]+[n]\}\phi_{n,J,l}.
\label{spegcl}
\end{align}
Recall the Laplacian operators on $\ASU$ and on $\Asq$ from equations \eqref{la3} and \eqref{la2d}:
\begin{align}
&\Box_{\SU}\phi=-(\nu X_{-}X_{+}+\beta X_{+}X_{-}+\gamma X_{z}X_{z})\lt\phi,\qquad &\phi\in\,\cl^{(0)}_{n}, \nn \\
&\Box_{\sq}f=-(\nu X_{-}X_{+}+\beta X_{+}X_{-})\lt f, \qquad & f\in\,\Asq\simeq\cl^{(0)}_{0}, \nn \\
&\Box_{D_{0}}\phi=-q^{-2n}\left(\nu X_{-}X_{+}+\beta X_{+}X_{-}\right)\lt\phi, \qquad & \phi\in\,\cl^{(0)}_{n}.
\label{glp}
\end{align}
One has that the  restriction of $\Box_{D_{0}}$ to $\phi\in\,\cl^{(0)}_{0}$ coincides with the operator $\Box_{\sq}$. Moreover it is now possible to generalise to the quantum Hopf bundle with the specific differential calculi studied so far, the classical relation \eqref{uno}, from which this analysis started: 
\beq
q^{2n}\Box_{D_{0}}\lt\phi=\left(\Box_{\SU}+\gamma X_{z}X_{z}\right)\lt\phi, \qquad\phi\in\,\cl^{(0)}_{n}.
\label{due}
\eeq
This relation appears as the natural generalisation of the classical relation \eqref{uno} to this specific quantum setting. The quantum Casimir operator \eqref{cas} can not be written as a polynomial in the basis derivations $X_{j}$ \eqref{Xq} of the 3D left covariant calculus from Woronowicz, so its role is played by the Laplacian $\Box_{\SU}$. Its quantum vertical part can still be written as a quadratic operator in the vertical field $X_{z}$ of the quantum Hopf fibration.

\section{An algebraic formulation of the classical Hopf bundle}\label{limcla}

The aim of this section is to apply the formalism developed to study the quantum Hopf bundle to the case  when all the space algebras are commutative, in order to recover the standard formulation of the classical Hopf bundle described at the beginning of the paper, from a dual viewpoint. 

\subsection{An algebraic description of the differential calculus on the group manifold $SU(2)$}

Rephrasing the relations \eqref{desu} which define the matrix Lie group $SU(2)$,  
the coordinate algebra  $\ca(SU(2))$ of the simple Lie group $SU(2)$ is the commutative $*$-algebra generated by 
$u$ and $v$, satisfying the spherical relation $u^{*}u+v^{*}v=1$. The Hopf algebra structure is given by the coproduct:
\beq
\Delta\,\left[
\begin{array}{cc} u & -v^* \\ v & u^*
\end{array}\right]=\left[
\begin{array}{cc} u & -v^* \\ v & u^*
\end{array}\right]\otimes\left[
\begin{array}{cc} u & -v^* \\ v & u^*
\end{array}\right] ,
\eeq
antipode:
\beq
S\,\left[
\begin{array}{cc} u & -v^* \\ v & u^*
\end{array}\right]=\left[
\begin{array}{cc} u^* & v^* \\ -v & u
\end{array}\right],
\eeq
and counit:
\beq
\epsilon \left[
\begin{array}{cc} u & -v^* \\ v & u^*
\end{array}\right]=\left[
\begin{array}{cc} 1 & 0 \\ 0 & 1
\end{array}\right].
\eeq
The universal envelopping algebra $\cu\mathfrak{(su(2))}$ is the Hopf $*$-algebra generated by the three elements 
$\e,\f,\h$ which satisfy  the algebraic relations \eqref{clil} coming from the Lie algebra structure in $\mathfrak{su(2)}$:
\begin{align}
&[\e,\f]=2\h, \nn \\
&[\f,\h]=\f, \nn \\
&[\e,\h]=-\e.
\label{clie}
\end{align}
The $*$-structure is: 
$$
\h^*=\h,\qquad\e^*=\f,\qquad\f^*=\e,
$$
and the Hopf algebra structure is provided by the coproduct:
\begin{align}
&\Delta(\e)=\e\otimes 1+1\otimes \e, \nn\\ 
&\Delta(\f)=\f\otimes 1+1\otimes \f, \nn \\  
&\Delta(\h)=\h\otimes 1+1\otimes \h; \nn
\end{align}
antipode:
$$
S(\e)=-\e, \qquad 
 S(\f)=-\f, \qquad 
S(\h)=-\h;
$$
and a counit which is trivial:
$$
\varepsilon(\e)=\varepsilon(\f)=\varepsilon(\h)=0.
$$
The centre of the algebra $\cu\mathfrak{(su(2))}$ is generated by the Casimir element:
\beq
C=\h^2+\frac{1}{2}(\e\f+\f\e)
\label{cca}
\eeq
The irreducible finite dimensional $*$-representations $\sigma_{j}$ of $\cu\mathfrak{(su(2))}$ are well known and labelled by nonnegative half-integers $j\in\frac{1}{2}\IN$. They are given by:
\begin{align}
&\sigma_{j}(\h)\ket{j,m}=m\ket{j,m}, \nn \\ 
&\sigma_{j}(\e)\ket{j,m}=\sqrt{(j-m)(j+m+1)}\ket{j,m+1}, \nn\\
&\sigma_{j}(\f)\ket{j,m}=\sqrt{(j-m+1)(j+m)}\ket{j,m-1}.
\end{align}
The algebras $\ca(SU(2))$ and $\cu\mathfrak{(su(2))}$ are dually paired. The bilinear \eqref{extc} mapping $\hs{~}{~}:\cu\mathfrak{(su(2))}\times\ca(SU(2))\,\to\,\IC$, compatible with the $*$-structures,  is set by: 
\begin{align}
&\hs{\h}{u}=-1/2, \nn \\ 
&\hs{\h}{u^*}=1/2, \nn \\ 
&\hs{\e}{v}=1, \nn \\ 
&\hs{\f}{v^*}=-1;
\label{cp}
\end{align}  
all other couples of generators pairing to $0$. This pairing is non degenerate:  the condition  $\hs{l}{x}=0\,\, \forall l\in\cu\mathfrak{(su(2))}$  implies  $x=0$, while  $\hs{l}{x}=0\,\, \forall x\in\ca(SU(2))$ implies $h=0$. 

It is possible to prove \cite{KS97} that a finite dimensional vector space $\mathcal{X}\subset\ch^{\prime}$ of linear functionals on a Hopf algebra $\ch$ is a tangent space of a finite dimensional left covariant first order differential calculus $(\Omega^{1}(\ch),\dd)$ if and only if $X(1)=0$ and 
$(\Delta(X)-\varepsilon\otimes X)\,\in\mathcal{X}\otimes\ch^{o}$,
for any $X\in\mathcal{X}$, where $\ch^{o}\subset\ch^{\prime}$ is the dual Hopf  algebra to $\ch$. The ideal $\cq=\{x\in \ker\,\varepsilon_{\ch}:\,X(x)=0\,\forall X\in\mathcal{X}\}$ characterises the calculus,  the bimodule of 1-forms being isomorphic to $\Omega^{1}(\ch)=\Omega^{1}_{un}(\ch)/\cn_{\cq}$ with $\cn_{\cq}=r^{-1}(\ch\otimes\cq)$. This result shows the path to prove the following proposition.

\begin{prop}
\label{clac}
Given the nondegenerate bilinear pairing $\hs{~}{~}:\cu\mathfrak{(su(2))}\times\ca(SU(2))\to\IC$ as in \eqref{cp}, the set $\{\e,\f,\h\}$ of generators in $\cu\mathfrak{(su(2))}$ defines a basis of the tangent space $\mathcal{X}_{SU(2)}$ for a bicovariant differential $*$-calculus on $\ca(SU(2))$. Such a differential calculus is isomorphic to the differential calculus \eqref{decla}, once the algebra $C^{\infty}(S^{3})$ is restricted to the polynomial algebra $\ca(SU(2))$.
\begin{proof}
 The definition of counit in the Hopf algebra $\cu\mathfrak{(su(2))}$ shows that the generators $l_{a}=\{\e,\f,\h\}$, seen as linear functionals on $\ca(SU(2))$ via the pairing, are such that:
\begin{align}
&\e(1)=\hs{\e}{1}=\varepsilon(\e)=0, \nonumber \\
&\f(1)=\hs{\f}{1}=\varepsilon(\f)=0, \nonumber \\
&\h(1)=\hs{\h}{1}=\varepsilon(\h)=0; \nonumber 
\end{align}
while the coproduct relations can be cast in the form:
\begin{align}
&\Delta(\e)-1\otimes\e=\e\otimes 1,\nn\\ 
&\Delta(\f)-1\otimes\f=\f\otimes 1,\nn\\ 
&\Delta(\h)-1\otimes\h=\h\otimes 1;
\label{coc}
\end{align} 
thus proving that the set $\{\e,\f,\h\}$ in $\cu\mathfrak{(su(2))}$ defines a complex vector space basis of a tangent space $\mathcal{X}_{SU(2)}$ for a left covariant differential calculus. The obvious inclusion $\mathcal{X}_{SU(2)}^{*}\subset\mathcal{X}_{SU(2)}$ proves, as described in section \ref{se:qphb}, that such a calculus admits a $*$ structure.

 In order to recover the ideal $\cq_{SU(2)}\subset \ker\varepsilon_{SU(2)}$ for this specific calculus, consider a generic element $x\in\,\ker\varepsilon_{SU(2)}$. It  must necessarily be written as $x=\{(u-1)x_{1},(u^{*}-1)x_{2},vx_{3},v^{*}x_{4}\}$ with $x_{j}\in\ca(SU(2))$. Such an element $x$ will belong to $\cq_{SU(2)}$ if $\hs{l_{a}}{x}=0$ for any of the generators $l_{a}\in\cu\mathfrak{(su(2))}$, since they form a vector space basis for  the tangent space $\mathcal{X}_{SU(2)}$ relative to this calculus. For the element $x=(u-1)x_{1}$ the three conditions are: 
\begin{align}
&\hs{\e}{(u-1)x_{1}}=\hs{\e}{u-1}\hs{1}{x_{1}}+\hs{1}{u-1}\hs{\e}{x_{1}}=0, \nn\\
&\hs{\f}{(u-1)x_{1}}=\hs{\f}{u-1}\hs{1}{x_{1}}+\hs{1}{u-1}\hs{\f}{x_{1}}=0, \nn \\
&\hs{\h}{(u-1)x_{1}}=\hs{\h}{u-1}\hs{1}{x_{1}}+\hs{1}{u-1}\hs{\h}{x_{1}}=-\frac{1}{2}\hs{1}{x_{1}}=-\frac{1}{2}\varepsilon(x_{1}),
\end{align}
where, in each of the three lines, the first equality comes from the general properties  of dual pairing and from the specific coproduct in $\cu\mathfrak{(su(2))}$, while the final result depends on the specific form of the  pairing. This means that $x=(u-1)x_{1}$ belongs to $\cq_{SU(2)}$ if and only if $x_{1}\in \ker\varepsilon_{SU(2)}$. The analysis is similar for the other three elements $x=\{(u^{*}-1)x_{2}, vx_{3},v^{*}x_{4}\}$. It is then proved that  this  left covariant differential calculus on $\ca(SU(2))$ - whose tangent space is 3 dimensional -  can be characterised by the  ideal $\cq_{SU(2)}=\{\ker\varepsilon_{SU(2)}\}^{2}\subset\,\ker\varepsilon_{SU(2)}$, which is generated by the ten elements: $\cq_{SU(2)}=\{(u-1)^2,(u-1)(u^{*}-1),(u-1)v, (u-1)v^{*},(u^{*}-1)^2,(u^{*}-1)v,(u^{*}-1)v^{*},v^2,vv^{*},v^{*2}\}$.  The equation \eqref{ded} allows then to write the exterior derivative for this calculus as:
\beq
\dd x=(\e\lt x)\omega_{\e}+(\f\lt x)\omega_{\f}+(\h\lt x)\omega_{\h}
\label{dcd}
\eeq
The commutation properties between the left invariant forms $\{\omega_{\e},\omega_{\f},\omega_{\h}\}$ and elements of the algebra $\ca(SU(2))$ depend on the  functionals $f_{ab}$ defined as $\Delta(l_{a})=1\otimes l_{a}+l_{b}\otimes f_{ba}$. From \eqref{coc} one has  $f_{ab}=\delta_{ab}$, so  1-forms do commute with elements of the algebra $\ca(SU(2))$, $\omega_{a}x=x\omega_{a}$.

The ideal $\cq_{SU(2)}$ is in addition stable under the right coaction $\Ad$ of the algebra $\ca(SU(2))$ onto itself: $\Ad(\cq_{SU(2)})\subset\cq_{SU(2)}\otimes\ca(SU(2))$. The proof of this result consists of a direct computation. The stability of the ideal $\cq_{SU(2)}$ under the right coaction $\Ad$ means that this differential calculus is bicovariant.

The explicit form of the left action of the generators of $\cu\mathfrak{(su(2))}$ on the generators of the coordinate algebra $\ca(SU(2))$ is: 
\beq
\begin{array}{l}
\h\lt u=-\frac{1}{2} u \\
\h\lt u^{*}=\frac{1}{2} u^{*} \\
\h\lt v=-\frac{1}{2} v \\
\h\lt v^{*}=\frac{1}{2} v^{*} 
\end{array}
\qquad\qquad
\begin{array}{l}
\e\lt u=-v^{*}  \\
\e\lt u^{*}=0   \\
\e\lt v= u^{*}  \\
\e\lt v^{*}=0 
\end{array}
\qquad\qquad
\begin{array}{l}
\f\lt u=0  \\
\f\lt u^{*}=v   \\
\f\lt v=0   \\
\f\lt v^{*}=-u 
\end{array}
\eeq
Starting from these relations it is immediate to see  that the left action of the generators $l_{a}\in\cu\mathfrak{(su(2))}$ is equivalent  to the Lie derivative along the left invariant vector fields $L_{a}$ \eqref{lcf}. This equivalence can now be written as:
\begin{align}
&\e\lt(x)=-iL_{+}(x),\nn \\
&\f\lt(x)=-iL_{-}(x), \nn \\ 
&\h\lt(x)=iL_{z}(x),
\label{efhlie}
\end{align}
and it is valid for any $x\in\ca(SU(2))$, as the Leibniz rule for the action of the derivations $L_{a}$ is encoded in the definition of the left action \eqref{deflr} and the properties of the functionals $f_{ab}=\delta_{ab}$.  From relation \eqref{dcd} it is possible to recover:
\begin{align*}  
&\dd u=-v^{*}\omega_{\e}-\frac{1}{2}u\omega_{\h}, \\ 
&\dd u^{*}=v\omega_{\f}+\frac{1}{2}u^{*}\omega_{\h}, \\
&\dd v=u^{*}\omega_{\e}-\frac{1}{2}v\omega_{\h}, \\
&\dd v^{*}=-u\omega_{\f}+\frac{1}{2}v^{*}\omega_{\h}.
\end{align*}
These relations can be inverted, so that left invariant 1-forms $\{\omega_{\e},\omega_{\f},\omega_{\h}\}$ can be compared to \eqref{linv}:
\begin{align}
&\omega_{\e}=u\dd v-v\dd u=i\tilde{\omega}_{+}, \nn \\
&\omega_{\f}=v^{*}\dd u^{*}-u^{*}\dd v^{*}=i\tilde{\omega}_{-}, \nn \\
&\omega_{\h}=-2(u^{*}\dd u+v^{*}\dd v)=-i\tilde{\omega}_{z}.
\label{isof}
\end{align}
The $*$-structure is given, on the basis of left-invariant generators, as $\omega_{\e}^{*}=-\omega_{\f}, \,\omega_{\h}^{*}=-\omega_{\h}$. 
The equalities \eqref{isof}, which are dual to \eqref{efhlie}, represent the isomorphism between the first order differential calculus   introduced via the action of the exterior derivative in \eqref{dcd}, and the one  analysed in section \ref{se:cdc}.

\end{proof}
\end{prop}

\noindent It is now straightforward to recover this bicovariant calculus as the classical limit of the quantum 3D left covariant calculus $(\Omega(\SU,\dd)$ described in section \ref{se:lcc}. In the classical  limit  $\ASU\to\ca(SU(2))$ as $q\to 1$, with $\phi\rightarrow x$, one has:
$$
\begin{array}{lcl}
\omega_{+}\rightarrow\omega_{\e},&\qquad&(X_{+}\lt\phi)\rightarrow(\e\lt x), \nn \\
\omega_{-}\rightarrow\omega_{\f},&\qquad&(X_{-}\lt\phi)\rightarrow(\f\lt x), \nn \\
\omega_{z}\rightarrow-\frac{1}{2}\,\omega_{\h},&\qquad&(X_{z}\lt\phi)\rightarrow(-2\h\lt x). \nn
\end{array}
$$
The coaction $\Delta_{R}^{(1)}$ of $\ca(SU(2))$  on the basis of left invariant forms defines the matrix $\Delta_{R}^{(1)}(\omega_{a})=\omega_{b}\otimes J_{ba}$:
\begin{align}
&\Delta_{R}^{(1)}(\omega_{\f})=\omega_{\f}\otimes u^{*2}+\omega_{\h}\otimes u^{*}v^{*}-\omega_{\e}\otimes v^{*2}, \nn\\
&\Delta_{R}^{(1)}(\omega_{\h})=-\omega_{\f}\otimes2u^{*}v+\omega_{\h}\otimes(u^{*}u-v^{*}v)-\omega_{\e}\otimes 2uv^{*}, \nn\\
&\Delta_{R}^{(1)}(\omega_{\e})=-\omega_{\f}\otimes v^{2}+\omega_{\h}\otimes{uv}+\omega_{\e}\otimes u^{2},
\label{Jd}
\end{align} 
which is used to define a basis of right invariant one forms $\eta_{a}=\omega_{b}S(J_{ba})$:
\begin{align}
&\eta_{\f}=u^{2}\omega_{\f}-uv^{*}\omega_{\h}-v^{*2}\omega_{\e}=v^{*}\dd u-u\dd v^{*}, \nn\\
&\eta_{\h}=2uv\omega_{\f}+(uu^{*}-vv^{*})\omega_{\h}+2u^{*}v^{*}\omega_{\e}=2(u\dd u^{*}+v^{*}\dd v), \nn\\
&\eta_{\e}=-v^{2}\omega_{\f}-u^{*}v\omega_{\h}+u^{*2}\omega_{\e}=u^{*}\dd v-v\dd u^{*};
\label{deeta}
\end{align}
- note that it has been made explicit use of the commutativity between forms $\omega_{a}$ and elements of the algebra $\ca(SU(2))$. The right acting derivation associated to this basis are given by \eqref{dedr} as 
$$
\dd x=\eta_{a}\rt(-S^{-1}(l_{a}))= \eta_{a}\rt l_{a}
$$ 
for any $x\in\ca(SU(2))$,  since an immediate evaluation gives $S^{-1}(l_{a})=-l_{a}$ for the three vector basis elements of the tangent space $l_{a}\in\mathcal{X}$. Using again the commutativity of the right invariant one forms $\eta_{a}$ with element of $\ca(SU(2))$, the action of the exterior derivation \eqref{dcd} can be written as:
\beq
\dd x=(x\rt\f)\eta_{\f}+(x\rt\h)\eta_{\h}+(x\rt\e)\eta_{\e}.
\label{dde}
\eeq
Comparing \eqref{deeta} to \eqref{rinv} one has:
\begin{align}
&\eta_{\f}=i\tilde{\eta}_{-},\nn \\
&\eta_{\h}=-i\tilde{\eta}_{z},\nn \\ 
&\eta_{\e}=i\tilde{\eta}_{+},
\label{ifrr}
\end{align}
while for the right action of the generators $l_{a}$ on $\ca(SU(2))$ one computes:
\beq
\begin{array}{l}
u\rt\h=-\frac{1}{2} u \\ 
u^{*}\rt\h=\frac{1}{2} u^{*} \\
v\rt\h=\frac{1}{2} v \\
v^{*}\rt\h=-\frac{1}{2} v^{*} 
\end{array}
\qquad
\begin{array}{l}
u\rt\e=0   \\
u^{*}\rt\e=-v^{*}    \\
v\rt\e= u  \\
v^{*}\rt\e=0 
\end{array}
\qquad
\begin{array}{l}
u\rt\f=v  \\
u^{*}\rt\f=0   \\
v\rt\f=0   \\
v^{*}\rt\f=-u^{*}; 
\end{array}
\eeq
so that the identification with the action of the right invariant vector fields \eqref{rcf} can be recovered as:
\begin{align}
&(x)\rt\f=-iR_{-}(x),\nn \\
&(x)\rt\e=-iR_{+}(x),\nn \\
&(x)\rt\h=iR_{z}(x),
\label{efhrie}
\end{align}
being dual to the identification \eqref{ifrr}. It is also evident that relations \eqref{ifrr} and \eqref{efhrie} define a different realisation of the isomorphism between the differential calculus introduced in this section \eqref{dde} and the differential calculus from section \ref{se:cdc}. 

\begin{rema}
The identification \eqref{efhlie} can be read as a Lie algebra isomorphism between the Lie algebra $\{\e,\f,\h\}$ given in \eqref{clie} and the Lie algebra of the left invariant vector fields $\{L_{a}\}$ \eqref{clil}:
\beq
\e=-iL_{+},\qquad\f=-iL_{-},\qquad\h=iL_{z}.
\label{efhcq}
\eeq
The notion of pairing between the algebras $\cu\mathfrak{(su(2))}$ and $\ca(SU(2))$ can be recovered  as the Lie derivative of the coordinate functions along the vector fields $L_{a}$, evaluated at the identity of the group manifold. The terms in \eqref{cp} giving  the nonzero terms of the pairing are:
$$
\begin{array}{lcl}
\left.L_{z}(u)\right|_{\id}=\frac{i}{2}&\qquad\Rightarrow\qquad&\hs{\h}{u}=-\frac{1}{2} \\
\left.L_{z}(u^{*})\right|_{\id}=-\frac{i}{2}&\qquad\Rightarrow\qquad&\hs{\h}{u^{*}}=\frac{1}{2} \\
\left.L_{+}(v)\right|_{\id}=i&\qquad\Rightarrow\qquad&\hs{\e}{v}=1 \\
\left.L_{-}(v^{*})\right|_{\id}=-i&\qquad\Rightarrow\qquad&\hs{\f}{v^{*}}=-1
\end{array}
$$
\end{rema}

The whole exterior algebra $\Omega(SU(2))$ can now be constructed from the differential calculus \eqref{dcd}. Any 1-form $\theta\in\,\Omega^{1}(SU(2))$ can be written on the basis of left invariant forms as $\theta=\sum_{k}x_{k}\omega_{k}=\omega_{k}x_{k}$ with $x_{k}\in\,\ca(SU(2))$. Higher dimensional forms can be defined by requiring their total antisimmetry, and that $\dd^{2}=0$. One has then $\omega_{a}\wedge\omega_{b}+\omega_{b}\wedge\omega_{a}=0$ and:
\begin{align}
&\dd\omega_{\f}=\omega_{\h}\wedge\omega_{\f},\nonumber \\
&\dd\omega_{\e}=\omega_{\e}\wedge\omega_{\h}, \nonumber \\
&\dd\omega_{\h}=2\omega_{\f}\wedge\omega_{\e}.
\label{2fc}
\end{align}
Finally, there is a unique volume top form $\omega_{\f}\wedge\omega_{\e}\wedge\omega_{\h}$.

The algebra $\ca(SU(2))$ can be partitioned into finite dimensional blocks, whose elements are related to the Wigner $D$-functions \cite{Mos} for the group $SU(2)$. Considering all the unitary irreducible representations of $SU(2)$, their matrix elements will give a Peter-Weyl basis for the Hilbert space $\cl^{2}(SU(2),\mu)$ of complex valued functions defined on the group manifold with respect to the Haar invariant measure. The Wigner $D$-function $D^{J}_{ks}(g)$ is defined to be the matrix element ($k,s$ are the matrix indices) representing the element $g\simeq(u,v)$ in $SU(2)$ \eqref{desu} in the representation of weight $J$. They are known:
\beq
D^{J}_{ks}=(-i)^{s+k}[(J+s)!(J-s)!(J+k)!(J-k)!]^{1/2}\sum_{l}(-1)^{k+l}\frac{u^{*l}v^{*J-k-l}v^{J-s-l}u^{*k+s+l}}{l!(J-k-l)!(J-s-l)!(s+k+l)!}
\label{wdf}
\eeq 
with $J=0,1/2,1,\ldots$ and $k=-J,\ldots,+J$, $s=-J,\ldots,+J$. In \eqref{wdf} the index $l$ runs over the set of natural numbers such that all the arguments of the factorial are non negative. To illustrate the meaning of this partition, proceed as in the quantum setting, and  consider the element $u^{*}\in\ca(SU(2))$. Representing the left action $\f\lt$ with a horizontal arrow and the right action $\rt\e$ with a vertical one yields the box:
\beq
\begin{array}{ccc}
u^{*} & \rightarrow & v \\
\downarrow & & \downarrow \\
-v^{*} & \rightarrow & u 
\end{array}
\label{bou}
\eeq
while starting from $u^{*2}\in\ca(SU(2))$ yields the box:
\beq 
\begin{array}{ccccc}
u^{*2} & \rightarrow & 2u^{*}v & \rightarrow & 2v^{2} \\
\downarrow & & \downarrow & & \downarrow \\
-2u^{*}v^{*} & \rightarrow & 2(u^{*}u-v^{*}v) & \rightarrow & 4uv \\
\downarrow & & \downarrow & & \downarrow \\
2v^{*2} & \rightarrow & -4v^{*}u & \rightarrow & 4u^{2}
\end{array}
\label{bod}
\eeq 
A recursive structure emerges now clear. For each positive integer $p$ one has a box $W_{p}$ made up of the $(p+1)\times(p+1)$ elements $w_{p:t,r}=\f^{t}\lt u^{*p}\rt \e^{r}$. An explicit calculation proves that:
\beq
\f^{t}\lt u^{*p}\rt \e^{r}=i^{t+r}j!\left[\frac{t!r!}{(p-t)!(p-r)!}\right]^{1/2}D^{p/2}_{t-p/2,r-p/2}
\eeq
with $t\leq p, r\leq p$. As an element in $\cu(\mathfrak{su(2))}$, the quadratic Casimir $C$ \eqref{cca} of the Lie algebra $\mathfrak{su(2)}$ acts on $x\in\ca(SU(2))$ as $C\lt x=x\rt C$, and its action clearly commutes with the actions $\f\lt$ and $\rt\e$. This means that the decomposition $\ca(SU(2))=\oplus_{j\in \IN}W_{p}$ gives the spectral resolution of the action of $C$:
\beq
C\lt w_{p:t,r}=\frac{p}{2}(\frac{p}{2}+1)w_{p:t,r}.
\label{egC}
\eeq

\subsection{The bundle structure}

\subsubsection{The base algebra of the bundle}

Given the abelian $*$-algebra $\ca(\U(1))=\IC[z,z^{*}]/<zz^{*}-1>$, the map $\check{\pi}:\ca(SU(2))\to(\U(1))$
\beq
\check{\pi}\left[\begin{array}{cc} u & -v^{*} \\ v & u^{*} \end{array}\right]=
\left[\begin{array}{cc} z & 0 \\ 0 & z^{*} \end{array}\right],
\label{chp}
\eeq 
is a surjective Hopf $*$-algebra homomorphism, so that $\ca(\U(1))$ can be seen as a $*$-subalgebra of $\ca(SU(2))$, with a right coaction:
\beq
\check{\Delta}_{R}=(1\otimes \check{\pi})\circ\Delta,\qquad\ca(SU(2)\to\ca(SU(2))\otimes\ca(\U(1)).
\label{rcgd}
\eeq
The coinvariant elements for this coaction, that is elements $b\in\ca(SU(2))$ for which $\check{\Delta}_{R}(b)=b\otimes 1$, form the subalgebra $\ca(S^{2})\subset\ca(SU(2))$, which is the coordinate subalgebra of the sphere $S^{2}$. From:
\begin{align}
&\check{\Delta}_{R}(u)=u\otimes z, \nn \\
&\check{\Delta}_{R}(u^{*})=u^{*}\otimes z^{*} , \nn \\
&\check{\Delta}_{R}(v)=v\otimes z, \nn \\
&\check{\Delta}_{R}(v^{*})=v^{*}\otimes z^{*} ,
\label{rcg}
\end{align}
one has that a set of generators for $\ca(S^{2})$ is given by \eqref{bdef}:
\begin{align}
&b_{z}=uu^{*}-vv^{*}, \nn \\
&b_{y}=uv^{*}+vu^{*}, \nn \\
&b_{x}=-i(vu^{*}-uv^{*})
\end{align}
The comparison with section \ref{se:pbc} shows that $\check{\pi}$ dually describes the choice of the gauge group $\U(1)$ as a subgroup of $SU(2)$, whose right principal pull-back action $\check{\mathrm{r}}_{k}^{*}$  is now replaced by the right $\ca(\U(1))$- coaction $\check{\Delta}_{R}$. The basis of the principal Hopf bundle $S^{2}\simeq SU(2)/\U(1)$ will be given as the algebra $\ca(S^{2})$ of right coinvariant elements $b_{a}\in\ca(SU(2))$, which is a homogeneous space algebra. The  coproduct $\Delta$ of $\ca(SU(2))$ restricts to a left coaction $\Delta:\ca(SU(2))\mapsto\ca(SU(2))\otimes \ca(S^{2})$ as:
\begin{align}
&\Delta(b_{\f})=u^{2}\otimes b_{\f}- v^{*}u\otimes b_{\h}-v^{*2}\otimes b_{\e}, \nn \\
&\Delta(b_{\h})=2uv\otimes b_{\f}+(u^{*}u-v^{*}v)\otimes b_{\h}+2u^{*}v^{*}\otimes b_{\e}, \nn \\
&\Delta(b_{\e})=-v^{2}\otimes b_{\f}-u^{*}v\otimes b_{\h}+u^{*2}\otimes b_{\e}.
\end{align}
with $b_{\f}=1/2(b_{y}-ib_{x})=uv^{*}$, $b_{\e}=1/2(b_{y}+ib_{x})=vu^{*}$, $b_{\h}=b_{z}$. The choice of this specific basis shows that $\Delta(b_{a})=S(J_{ka})\otimes b_{k}$ where the matrix $J$ is exactly the one defined in \eqref{Jd} as $\Delta_{R}^{(1)}(\omega_{a})=\omega_{b}\otimes J_{ba}$. 

The identification \eqref{efhlie} between the left action $\h\lt x$ -- given the generator $\h\in\,\cu(\mathfrak{su(2)})$ on any $x\in\,\ca(SU(2))$ -- and the action $iL_{z}(x)$ -- given the left invariant vector field $L_{z}$ -- as well as the definition of the $\ca(\U(1))$-right coaction $\check{\Delta}_{R}$ on $\ca(SU(2))$ \eqref{rcg}, allow to recover the set of the $\U(1)$-equivariant functions $\mathfrak{L}_{n}^{(0)}\subset\ca(SU(2))$ in 
\eqref{cln0} as:
\beq
\mathfrak{L}_{n}^{(0)}=\{\phi\in\,\ca(SU(2))\,:\,\h\lt\phi=\frac{n}{2}\phi\,\Leftrightarrow\,\check{\Delta}_{R}(\phi)=\phi\otimes z^{-n}\}.
\label{clnc2}
\eeq

\subsubsection{A differential calculus on the gauge group algebra}\label{secu}

The strategy underlining the proof of the  proposition \ref{clac} brings also to the definition of a differential calculus on the gauge group algebra $\ca(\U(1))$. The bilinear pairing $\hs{\cdot}{\cdot}:\cu(\mathfrak{su(2)})\times\ca(SU(2))\to\IC$ \eqref{cp} is restricted via the surjection $\check{\pi}$ \eqref{chp} to  a bilinear pairing   $\hs{\cdot}{\cdot}:\cu\{\h\}\times\ca(\U(1))\to\IC$, which is still compatible with the $*$-structure, given on generators as:
\begin{align}
&\hs{\h}{z}=-\frac{1}{2}, \nn \\
&\hs{\h}{z^{-1}}=\frac{1}{2}. \nn
\end{align}
The set $\mathcal{X}_{\U(1)}=\{\h\}$ is proved to be the basis of the tangent space for a 1-dimensional bicovariant commutative calculus on $\ca(\U(1))$. The ideal $\cq_{\U(1)}\subset\ker\varepsilon_{\U(1)}$ turns out again to be $\ca_{\U(1)}=(\ker\varepsilon_{\U(1)})^{2}$ generated by $\{(z-1)^{2}, (z-1)(z^{-1}-1),(z^{-1}-1)^{2}\}$, which can also be recovered as $\cq_{\U(1)}=\check{\pi}((\ker\varepsilon_{SU(2)})^{2})$. From:
\begin{align}
&\h\lt z=-\frac{1}{2}z, \nn \\
&\h\lt z^{-1}=\frac{1}{2}z^{-1}, \nn 
\end{align}
one has that:
\begin{align}
&\dd z=-\frac{1}{2}z\check{\omega}, \nn \\
&\dd z^{-1}=\frac{1}{2}z^{-1}\check{\omega} 
\end{align}
with $z\dd z=(\dd z)z$. The only left invariant 1-form is 
$$
\check{\omega}=-2z^{-1}\dd z=2z\dd z^{-1},
$$
while the role of the right invariant derivation associated to $\h\in\,\cu\{\h\}$ is played by $-S^{-1}(\h)=\h$, so that the right invariant form generating this calculus is:
$$
\begin{array}{lcl}
\dd z=\check{\eta}(z\rt \h)=\check{\eta}(-\frac{1}{2}z)&\qquad\Rightarrow\qquad& \check{\eta}=-2z^{-1}\dd z, \\
\dd z^{-1}=\check{\eta}(z^{-1}\rt \h)=\check{\eta}(\frac{1}{2}z^{-1})&\qquad \Rightarrow \qquad&\check{\eta}=2z\dd z^{-1} 
\end{array}
$$
so that one obtains $\check{\eta}=\check{\omega}$.

It is possible to characterise the quotient $\ker\varepsilon_{\U(1)}/\cq_{\U(1)}=\ker\varepsilon_{\U(1)}/(\ker\varepsilon_{\U(1)})^{2}$.  The three elements generating the ideal $\cq_{\U(1)}=(\ker\varepsilon_{\U(1)})^{2}$ can be written as:
\begin{align}
&\xi=(z-1)(z^{-1}-1)=(z-1)+(z^{-1}-1), \nn \\
&\xi^{\prime}=(z-1)(z-1)=\xi+\xi(z-1), \nn \\
&\xi^{\prime\prime}=(z^{-1}-1)(z^{-1}-1)=\xi+\xi(z^{-1}-1), \nn 
\end{align}
so that $\cq_{\U(1)}$ can be seen generated by $\xi=(z-1)+(z^{-1}-1)$. Set a map $\lambda:\ker\varepsilon_{\U(1)}\to\IC$ by $\lambda(u(z-1))=\sum_{j\in\,\IZ}u_{j}$, where $u=\sum_{j\in\,\IZ}u_{j}z^{j}$ is generic element in $\ca(\U(1))$. The techniques outlined in lemma \ref{pu} in the quantum setting enable to prove that $\lambda$ can be used to define a complex vector space isomorphism between $\ker\varepsilon_{\U(1)}/(\ker\varepsilon_{\U(1)})^{2}$ and $\IC$, whose inverse is given by $\lambda^{-1}:w\in\,\IC\mapsto\,\lambda^{-1}(w)=w(z-1)\,\in\,\ker\varepsilon_{\U(1)}$. It is evident that such a map $\lambda$ gives the projection $\pi_{\cq_{\U(1)}}:\ker\varepsilon_{\U(1)}\to\ker\varepsilon_{\U(1)}/\cq_{\U(1)}\simeq\IC$, since it chooses a representative in each equivalence class in the quotient $\ker\varepsilon_{\U(1)}/\cq_{\U(1)}$.

\subsubsection{The Hopf bundle structure}
With the 3D bicovariant calculus on the total space algebra $\ca(SU(2))$ and the 1D bicovariant calculus on the gauge group algebra $\ca(\U(1))$, one needs to prove the compatibility conditions that lead to the exact sequence:
\begin{multline*}
0\, \to\, \ca(SU(2))\left(\Omega^{1}(S^{2})\right)\ca(SU(2))\, \to \\ \to\,\Omega^{1}(\ca(SU(2))\,\stackrel{\sim_{\mathcal{N}_{SU(2)}}}  \longrightarrow \,\ca(SU(2))\otimes \ker\varepsilon_{\U(1)}/\mathcal{Q}_{\U(1)}\,\to\,0 ,
\end{multline*}
where  the map $\sim_{\mathcal{N}_{SU(2)}}$ is defined as in the diagram \eqref{qdia} which now acquires the form:
\beq
\begin{array}{lcl}
\Omega^{1}(SU(2))_{un} & \stackrel{ \pi_{\cq_{SU(2)}} }{\longrightarrow} 
& \Omega^{1}(\ca(SU(2)) \\
\downarrow \chi  &  & \downarrow \sim_{\cn_{SU(2))}} \\
\ca(SU(2)) \otimes \ker\varepsilon_{\U(1)} & \stackrel{ \id\otimes\pi_{\cq_{\U(1)}} }{\longrightarrow} 
&\ca(SU(2))\otimes (\ker\varepsilon_{\U(1)}/\mathcal{Q}_{\U(1)}) \,.
\end{array}
\label{qdiacl}
\eeq
The proof of the compatibility conditions is in the following lemmas. The first one analyses the right covariance of the differential structure on $\ca(SU(2))$. 
\begin{lemm}
\label{l1c}
From the 3D bicovariant calculus on $\ca(SU(2))$ generated by the ideal $\cq_{SU(2)}=(\ker\varepsilon_{SU(2)})^{2}\subset\ker\varepsilon_{SU(2)}$ given in proposition \ref{clac}, one has $\check{\Delta}_{R}\cn_{SU(2)}\subset\cn_{SU(2)}\otimes\ca(\U(1))$. 
\begin{proof} 
Using the bijection given in \eqref{3p2}, it is $\Omega^{1}(SU(2))\simeq\Omega^{1}(SU(2))/\cn_{SU(2)}$ with $\cn_{SU(2)}=r^{-1}(\ca(SU(2))\otimes\cq_{SU(2)})$.  For this specific calculus one has that $\cn_{SU(2)}$ is the sub-bimodule generated by $\{\delta\phi\,\delta\psi\}$ for any $\phi, \psi\,\in\,\ca(SU(2))$, where $\delta\phi=(1\otimes\phi-\phi\otimes 1)\in\,\Omega^{1}(SU(2))_{un}$. Choose $\phi\in\,\mathfrak{L}_{n}^{(0)}$ and $\psi\in\,\mathfrak{L}_{m}^{(0)}$ so to have $\check{\Delta}_{R}\phi=\phi\otimes z^{-n}$ and $\check{\Delta}_{R}\psi=\psi\otimes z^{-m}$. Extending the coaction $\check{\Delta}_{R}$ to a coaction $\check{\Delta}_{R}:\ca(SU(2))\otimes\ca(SU(2))\to\ca(SU(2))\otimes\ca(SU(2))\otimes\ca(U(1))$ as $\check{\Delta}_{R}=(\id\otimes \id\otimes m)\circ(\id\otimes \tau\otimes \id)\circ (\check{\Delta}_{R}\otimes\check{\Delta}_{R})$ in terms of the flip operator $\tau$, it becomes an easy calculation to find:
\begin{align}
\check{\Delta}_{R}(\delta\phi\,\delta\psi)&=(1\otimes\phi\psi+\phi\psi\otimes 1-\phi\otimes\psi-\psi\otimes\phi) \nn \\ 
&=(1\otimes\phi\psi+\phi\psi\otimes 1-\phi\otimes\psi-\psi\otimes\phi)\otimes z^{-m-n}=(\delta\phi\,\delta\psi)\otimes z^{-m-n}. \nn
\end{align}

\end{proof}
\end{lemm}

\begin{lemm}
\label{l2c}
The map $\chi:\Omega^{1}(SU(2))_{un}\to\ca(SU(2))\otimes\ca(\U(1))$ defined in \eqref{chimap} as $\chi=(m\otimes \id)\circ (\id\otimes \check{\Delta}_{R})$ is surjecive.
\begin{proof}
The proof of this result closely follows the proof of the proposition \ref{pro:un}. From the spherical relation  $1=(u^{*}u+v^{*}v)^{n}=\sum_{a=0}^{n}\left(\begin{array}{c} n \\ a \end{array}\right)\,u^{*a}v^{*n-a}v^{n-a}u^{a}$ it is possible to set $\ket{\Psin}_{a}\in\mathfrak{L}_{n}^{(0)}$ for  $a=0,\ldots,\mn$ with $\hs{\Psin}{\Psin}=1$ as:
\begin{align}
&n>0:\,\ket{\Psin}_{a}=\sqrt{\left(\begin{array}{c} n \\ a \end{array}\right)}\,v^{*a}u^{*n-a}, \nn \\
&n<0:\,\ket{\Psin}_{a}=\sqrt{\left(\begin{array}{c} \mn \\ a \end{array}\right)}\,v^{*\mn -a}u^{a}. \nn 
\end{align}
Fixed $n\in\,\IZ$, define $\gamma=\hs{\Psi^{(-n)}}{\delta\Psi^{(-n)}}$. Since $\ket{\Psi^{(-n)}}\in\,\mathfrak{L}_{-n}^{(0)}$, one computes that $\chi(\gamma)=1\otimes(z^{n}-1)$, and this 
sufficient to prove the surjectivity of the map $\chi$, being $\chi$ left $\ca(SU(2))$-linear and $\ker\varepsilon_{U(1)}$ is a complex vector space with a basis $(z^{n}-1)$. 

\end{proof}
\end{lemm}

\begin{lemm}
\label{l3c}
Given the map $\chi$ as in the previous lemma, it is $\chi(\cn_{SU(2)})\subset\ca(SU(2))\otimes\cq_{\U(1)}$, where $\cn_{SU(2)}$  is as in lemma \ref{l1c} and $\cq_{\U(1)}=(\ker\varepsilon_{\U(1)})^{2}$. 
\begin{proof}
To be definite, consider $\phi\in\,\mathfrak{L}_{n}^{(0)}$ and $\psi\in\,\mathfrak{L}_{m}^{(0)}$. One has:
\begin{align}
\chi(\delta\phi\,\delta\psi)&=\phi\psi\otimes\{z^{-n-m}+1-z^{-n}-z^{-m}\} \nn \\
&=\phi\psi\otimes\{(1-z^{-n})(1-z^{-m})\}\subset\ca(SU(2))\otimes(\ker\varepsilon_{\U(1)})^{2}.\nn
\end{align}

\end{proof}
\end{lemm} 

The results of these lemmas allow to define the map $\sim_{\cn_{SU(2)}}:\Omega^{1}(SU(2))\to\ca(SU(2))\otimes\ker\varepsilon_{\U(1)}/\cq_{\U(1)}$ from the diagram \eqref{qdiacl}. Using the isomorphism $\lambda:\ker\varepsilon_{\U(1)}/\cq_{\U(1)}\to\IC$ described in section \ref{secu}, one has:
\begin{align}
\sim_{\cn_{SU(2)}}(\omega_{\e})&=0 \nn \\
\sim_{\cn_{SU(2)}}(\omega_{\f})&=0 \nn \\
\sim_{\cn_{SU(2)}}(\omega_{\h})&=-2\otimes\pi_{\cq_{\U(1)}}(z-1)=-2\otimes1.  
\label{ondcl}
\end{align}

The next lemma completes the analysis of the compatibility conditions between the differential structures on $\ca(SU(2))$ and on $\ca(\U(1))$. The horizontal part of the set of $k$-forms  out of $\Omega^{k}(SU(2))$ is defined as $\Omega^{k}_{\mathrm{hor}}(SU(2))=\Omega^{k}(S^{2})\ca(SU(2))=\ca(SU(2))\Omega^{k}(S^{2})$.

\begin{lemm}
\label{l4c}
Given the differential calculus on the basis $\Omega^{1}(S^{2})=\Omega^{1}(S^{2})_{un}/\cn_{S^{2}}$ with $\cn_{S^{2}}=\cn_{SU(2)}\cap\Omega^{1}(S^{2})_{un}$, it is $\ker\sim_{\cn_{SU(2)}}=\Omega^{1}(S^{2})\ca(SU(2))=\ca(SU(2))\Omega^{1}(S^{2})=\Omega^{1}_{\mathrm{hor}}(SU(2))$.
\begin{proof}
Consider a 1-form $[\eta]\in\,\Omega^{1}(SU(2))$ and choose the element $\eta=\psi\,\delta\phi\in\,\Omega^{1}(SU(2))_{un}$ as a representative of $[\eta]$, with $\phi\in\,\mathfrak{L}_{n}^{(0)}$ and $\psi\in\,\mathfrak{L}_{m}^{(0)}$. One finds:
\begin{align}
&\chi(\psi\,\delta\phi)=\psi\phi\otimes(z^{-n}-1), \nn \\
&\sim_{\cn_{SU(2)}}(\eta)=\psi\phi\otimes\pi_{\cq_{\U(1)}}(z^{-n}-1). \nn 
\end{align}
Recalling once more the isomorphism $\lambda:\ker\varepsilon_{\U(1)}/\cq_{\U(1)}\to\IC$, it is $\lambda(z^{-n}-1)=0$ if and only if $n=0$, so to have $\eta=\psi\,\delta\phi$ with $\delta\phi\in\,\Omega^{1}(S^{2})_{un}$ and then $\eta\in\,\Omega^{1}(S^{2})_{un}\ca(\U(1))$.
It is clear that the condition $\chi(\cn_{SU(2)})\subset\ca(SU(2))\otimes\cq_{\U(1)}$ proved in lemma \ref{l3c} ensures that the map $\sim_{\cn_{SU(2)}}$ is well-defined: its image does not depend on the specific choice of the representative $\eta\in\,[\eta]\subset\Omega^{1}(SU(2))$. 

\end{proof}
\end{lemm}

The property of right covariance of the calculus on $\ca(SU(2))$ -- proved in lemma \ref{l1c} -- allows to extend the coaction $\check{\Delta}_{R}$ to a coaction $\check{\Delta}_{R}^{(k)}:\Omega^{k}(SU(2))\to\Omega^{k}(SU(2))\otimes\ca(\U(1))$ via $\check{\Delta}^{(k)}_{R}\circ\dd=(\dd\otimes id)\circ\check{\Delta}_{R}^{(k-1)}$. Via such a coaction it is possible to recover \eqref{Omr} the set $\Omega^{k}(SU(2))_{\rho_{(n)}}$ as the  $\rho_{(n)}(\U(1))$-equivariant k-forms on the Hopf bundle:
$$
\Omega^{k}(SU(2))_{\rho_{(n)}}=\{\phi\in\,\Omega^{k}(SU(2)):\,\check{\Delta}_{R}^{(k)}(\phi)=\phi\otimes z^{-n}\}.
$$ 
as well as the $\ca(S^{2})$-bimodule $\mathfrak{L}_{n}^{(k)}$ of horizontal elements in $\Omega^{k}(SU(2))_{\rho_{(n)}}$.

\subsubsection{Connections and covariant derivative on the classical Hopf bundle}

The compatibility conditions bring the exactness of the sequence:
\beq
0\longrightarrow\Omega^{1}_{\mathrm{hor}}(SU(2))\longrightarrow\Omega^{1}(SU(2))\stackrel{\sim_{\mathcal{N}_{SU(2)}}}{\longrightarrow}\ca(SU(2))\otimes\ker\varepsilon_{\U(1)}/\cq_{\U(1)},
\label{excl}
\eeq
whose every right invariant splitting 
$\sigma:\ca(SU(2))\otimes\ker\varepsilon_{\U(1)}/\cq_{\U(1)}\to\Omega^{1}(SU(2))$  represents a connection \eqref{si}. With $w\in\,\IC\simeq\ker\varepsilon_{\U(1)}/\cq_{\U(1)}$, one has:
\begin{align}
&\sigma(1\otimes w)=-\frac{w}{2}(\omega_{\h}+U\omega_{\e}+V\omega_{\f}), \nn \\
&\sigma(\phi\otimes w)=-\frac{w}{2}\phi(\omega_{\h}+U\omega_{\e}+V\omega_{\f}) \label{si3c}
\end{align}
where $\phi\in\,\ca(SU(2))$, and $U\in\,\mathfrak{L}_{2}^{(0)}$, $V\in\,\mathfrak{L}_{-2}^{(0)}$. The right invariant projection defined in\eqref{Pi} $\Pi:\Omega^{1}(SU(2))\to\Omega^{1}(SU(2))$ associated to this splitting is, from \eqref{ondcl}:
\begin{align}
&\Pi(\omega_{\e})=\Pi(\omega_{\f})=0, \nn \\
&\Pi(\omega_{\h})=\omega_{\h}+U\omega_{\e}+V\omega_{\f}. \label{Pi3c}
\end{align}
The connection one form $\omega:\ca(\U(1))\mapsto\Omega^{1}(SU(2))$ defined in \eqref{ome} is:
\beq
\omega(z^{n})=\sigma(1\otimes [z^{n}-1])=-\frac{n}{2}(\omega_{\h}+U\omega_{\e}+V\omega_{\f}). 
\label{ome3c}
\eeq
The horizontal projector $(1-\Pi):\Omega^{1}(SU(2))\to\Omega^{1}_{\mathrm{hor}}(SU(2))$ can be extended to whole exterior algebra $\Omega(SU(2))$, since it is compatible with the wedge product: one finds that $\{(1-\Pi)\omega_{a}\wedge(1-\Pi)\omega_{b}\}+\{(1-\Pi)\omega_{b}\wedge(1-\Pi)\omega_{a}\}=0$ or any pair of 1-forms. This property, which is \emph{not} valid in the quantum setting for a general connection -- recall the remark \ref{rempi} --, allows to define an operator of covariant derivative $D:\Omega^{k}(SU(2))\mapsto\Omega^{k+1}(SU(2))$ as:
\beq
D\phi=(1-\Pi)\dd\phi, \qquad \forall\,\phi\in\,\Omega^{k}(SU(2)).   
\label{cDg}
\eeq
This definition is the dual counterpart of definition \eqref{2p1}. It is not difficult to prove the main properties of such an operator of covariant derivative $D$:
\begin{itemize}
\item For any $\phi\in\,\Omega^{k}(SU(2))$, $D\phi\in\,\Omega^{k+1}_{\mathrm{hor}}(SU(2))$.
\item The operator $D$ is 'covariant'. One has $\check{\Delta}_{R}^{(k)}\phi=\phi\otimes z^{n}\,\Leftrightarrow\,\check{\Delta}_{R}^{(k+1)}(D\phi)=D\phi\otimes z^{n}$.
\item Given $\phi\in\,\mathfrak{L}_{n}^{(k)}$, that is $\phi\in\,\Omega^{k}_{\mathrm{hor}}(SU(2))$ such that $\check{\Delta}_{R}^{(k)}\phi=\phi\otimes z^{n}$, it is $D\phi=\dd\phi+\omega(z^{n})\wedge\phi$. This last property recovers the relation \eqref{2p2}. 
\end{itemize}

\section{Back on  a covariant derivative  on the exterior algebra $\Omega(\SU)$}\label{bproj}
The analysis in section \ref{se:conn} presents  the formalism of connections on a quantum principal bundle \cite{bm93} and explicitly describes both the set of connections on a quantum Hopf bundle 
and the corresponding set of covariant derivative operators $\nabla:\ce_{n}^{(k)}\to\ce_{n}^{(k+1)}$ acting on $k$-form valued sections  of the associated quantum line bundles. The left $\Asq$-module equivalence between $\ce_{n}^{(k)}$ and horizontal elements $\cl_{n}^{(k)}\subset\Omega_{\mathrm{hor}}^{k}(\SU)$ allows then for the definition of a covariant derivative operator $D:\cl_{n}^{(k)}\to\cl_{n}^{(k+1)}$ with $k=0,1,2$. 

 The equation \eqref{1mpsq} in  remark \ref{rempi} clarifies the reasons why, presenting a  connection via the projector \eqref{piom} $\Pi:\Omega^{1}(\SU)\to\Omega^{1}(\SU)$ given in  \eqref{Pi3b}, the operator $\check{D}=(1-\Pi)\dd:\Omega^{1}(\SU)\to\Omega^{2}_{\mathrm{hor}}(\SU)$ as in \eqref{Dc} defined a consistent covariant derivative on the whole exterior algebra $\Omega(\SU)$  only in the case of the monopole connection: the operator $(1-\Pi):\Omega^{1}(\SU)\to\Omega^{1}_{\mathrm{hor}}(\SU)$ is a covariant projector compatible with the properties of the wedge product \eqref{1mpsq0} in the exterior algebra $\Omega(\SU)$ only if the connection is the monopole connection.     

The problem of defining, for any connection on a principal quantum bundle,  a consistent covariant projection operator on the whole exterior algebra on the total space of the bundle whose range is given by the horizontal exterior forms has been studied in \cite{durI,durII}. The aim of this section is, from one side,  to describe the properties of the horizontal projector arising from that analysis, and then to show that in such a formulation  of the Hopf bundle more than one horizontal covariant projector can be consistently introduced.

As already mentioned, the formulation presented in  \cite{durI,durII}  of the geometrical structures of a quantum principal bundle  slightly differs from that described  in section \ref{QPB} and a comparison between them is in \cite{durcomm}. This formalism will not be explicitly reviewed: the main results concerning how to define an horizontal covariant projector will be translated into the language extensively described in the previous sections.

\bigskip


The differential $*$-calculus $(\Omega(\U(1)),\dd)$ on the gauge group algebra $\U(1)$ is described in  section \ref{se:csg}. It canonically corresponds to the right $\ca(\U(1))$-ideal $\cq_{\U(1)}\subset\ker\varepsilon_{\U(1)}$ generated by the element $\{(z^*-1)+q^{2}(z-1)\}$, so that by lemma \ref{pu} 
it is $\Omega^{1}(\U(1))_{\mathrm{inv}}\simeq\ker\varepsilon_{\U(1)}/\cq_{\U(1)}\simeq\IC.$ Such a calculus is bicovariant: given the left and right coactions \eqref{Dl1} of the $*$-Hopf algebra $\ca(\U(1))$ on $\Omega^{1}(\U(1))$ one has  
that the 1-form 
$\omega_{z}$ is both left and right invariant,
\beq
\begin{array}{lcl}
\Delta_{\ell}^{(1)}:\Omega^{1}(\U(1))\to\ca(\U(1))\otimes\Omega^{1}(\U(1)),&\qquad&\Delta_{\ell}^{(1)}(\omega_{z})=1\otimes\omega_{z};\\
\Delta_{\wp}^{(1)}:\Omega^{1}(\U(1))\to\Omega^{1}(\U(1))\otimes\ca(\U(1)),&\qquad&\Delta_{\wp}^{(1)}(\omega_{z})=\omega_{z}\otimes 1.
\end{array}
\label{lrinv}
\eeq
The exterior algebra on this differential calculus is built following \cite{KS97}, as explained in section \eqref{se:lcc}, where the same procedure has been applied to the analysis of the 3D left-covariant calculus on $\SU$. It results $S_{\cq_{\U(1)}}=(\Omega^{1}(\U(1)))^{\otimes2}$, so that 
\beq
\Omega(\U(1))\,=\,\sum_{k\geq0}\,^{\oplus}\Omega(\U(1))^{\wedge k}\,=\,\ca(\U(1))\oplus\Omega^{1}(\U(1)).
\label{Gsd}
\eeq
The coproduct map in the Hopf $*$-algebra $\ca(\U(1))$ can be extended to a homomorphism $\hat{\Delta}_{\U(1)}:\Omega(\U(1))\to\Omega(\U(1))\otimes\Omega(\U(1))$ given by  
\begin{align}
&\hat{\Delta}_{\U(1)}(\varphi)=\Delta(\varphi)=\varphi\otimes\varphi, \nn \\
&\hat{\Delta}_{\U(1)}(\varphi\,\omega_{z})=\Delta_{\ell}^{(1)}(\varphi\,\omega_{z})+\Delta_{\wp}^{(1)}(\varphi\,\omega_{z})=\varphi(1\otimes\varphi\,\omega_{z}+\omega_{z}\otimes\varphi),
\label{phih}
\end{align}
for any $\varphi\,\in\,\ca(\U(1))$. Given the principal bundle structure, the compatibility conditions among calculi on the total space algebra and the gauge group algebra allow to prove that there exists a unique extension of the coaction \eqref{cancoa} of the gauge group $\U(1)$ on the total space $\SU$ to a 
left $\ASU$-module homomorphism $\mathfrak{F}:\Omega(\SU)\to\Omega(\SU)\otimes\Omega(\U(1))$ implicitly defined by:
\begin{align*}
&(\mathfrak{F}\otimes\id)\mathfrak{F}=(\id\otimes\hat{\Delta}_{\U(1)})\mathfrak{F}, \\
&\mathfrak{F}*_{\SU}=(*_{\SU}\otimes *_{\U(1)})\mathfrak{F}:
\end{align*} 
where the second condition expresses a compatiblity between the map $\mathfrak{F}$ and the $*$-structures on the exterior algebras built over the calculi on $\SU$ and $\U(1)$. One has 
\begin{align}
&\mathfrak{F}(x)=\Delta_{R}(x)=x\otimes z^{-n}, \nn \\
&\mathfrak{F}(x\,\omega_{-})=x\,\omega_{-}\otimes z^{-2-n}, \nn \\
&\mathfrak{F}(x\,\omega_{+})=x\,\omega_{+}\otimes z^{2-n}, \nn \\
&\mathfrak{F}(x\,\omega_{z})=(x\otimes z^{-n}\omega_{z})+(x\,\omega_{z}\otimes z^{-n}), \nn \\
&\mathfrak{F}(x\,\omega_{-}\wedge\omega_{+})=x\,\omega_{-}\wedge\omega_{+}\otimes z^{-n}, \nn \\
&\mathfrak{F}(x\,\omega_{+}\wedge\omega_{z})=(x\,\omega_{+}\otimes z^{2-n}\omega_{z})+(x\,\omega_{+}\wedge\omega_{z}\otimes z^{2-n}), \nn \\
&\mathfrak{F}(x\,\omega_{z}\wedge\omega_{-})=(x\,\omega_{-}\otimes z^{-2-n}\omega_{z})+(x\,\omega_{z}\wedge\omega_{-}\otimes z^{-2-n}), \nn \\
&\mathfrak{F}(x\,\omega_{-}\wedge\omega_{+}\wedge\omega_{z})=(x\,\omega_{-}\wedge\omega_{+}\otimes z^{-n}\omega_{z})+(x\,\omega_{-}\wedge\omega_{+}\wedge\omega_{z}\otimes z^{-n}),
\label{Fha}
\end{align}
with $x\in\,\ASU$, such that $\Delta_{R}(x)=x\otimes z^{-n}\,\Leftrightarrow\,x\in\,\cl_{n}^{(0)}$. The homomorphism $\mathfrak{F}$ can be restricted to the right coaction $\Delta_{R}^{(k)}:\Omega^{k}(\SU)\to\Omega^{k}(\SU)\otimes\ca(\U(1))$ given in \eqref{dre}:
$$
\Delta^{(k)}_{R}(\phi)=(\id\otimes p_{0})\mathfrak{F}(\phi)
$$
with $\phi\in\,\Omega^{k}(\SU)$ and $p_{0}$ the projection $\Omega(\U(1))\to\ca(\U(1))$ coming from \eqref{Gsd}. The horizontal subset of the exterior algebra $\Omega(\SU)$ can be defined via:
\beq
\Omega_{\mathrm{hor}}(\SU)=\{\phi\in\,\Omega(\SU)\quad:\quad\mathfrak{F}(\phi)=(id\otimes p_{0})\mathfrak{F}(\phi)\},
\label{horD}
\eeq
while the exterior algebra $\Omega(\sq)$ described in section \ref{se:cals2} can be recovered as
$$
\Omega(\sq)=\{\phi\in\,\Omega(\SU)\quad:\quad\mathfrak{F}(\phi)=\phi\otimes 1\}.
$$
From the analysis in section \ref{se:conn} one has that a connection 1-form is given via a map $\tilde{\omega}:\Omega^{1}(\U(1))_{\mathrm{inv}}\to\Omega(\SU)$ satisfying the conditions \eqref{ome}.   The equation \eqref{ome3} shows that any connection can be written as:
$$
\tilde{\omega}(\omega_{z})=\omega_{z}+\mathrm{a},
$$
with $\mathrm{a}\in\,\Omega^{1}(\sq)$. Given a connection, one can define a map 
\beq
\mathit{m}_{\omega}:\Omega_{\mathrm{hor}}(\SU)\otimes\Omega(\U(1))_{\mathrm{inv}}\to\Omega(\SU),
\label{mmap}
\eeq
where the relation \eqref{Gsd} enables to recover $\Omega(\U(1))_{\mathrm{inv}}\simeq\{\IC\oplus\Omega^{1}(\U(1))_{\mathrm{inv}}\}$:
given $\psi\in\,\mathfrak{hor}(\SU)$ and $\theta=\lambda+\mu\,\omega_{z}\,\in\Omega(\U(1))_{\mathrm{inv}}$ (with $\lambda,\mu\in\,\IC$) set:
\beq
\mathit{m}_{\omega}(\psi\otimes\theta)=\psi\wedge(\mu+\lambda\tilde{\omega}(\omega_{z})).
\label{defmo}
\eeq
The map $\mathit{m}_{\omega}$ is proved to be bijective, and the operator 
\beq
h_{\omega}=(\id\otimes p_{0})\mathit{m}_{\omega}^{-1}
\label{hproj}
\eeq
 a covariant horizontal projector $h_{\omega}:\Omega(\SU)\to\mathfrak{hor}(\SU)$. Given an element $\phi\in\,\Omega^{k}(\SU)$, define its covariant derivative:
 \beq
\mathfrak{D}\phi=h_{\omega}\dd\phi.
\label{coDeD}
\eeq 
In the formulation developed in \cite{durI,durII} this definition is meant to be the quantum analogue of the classical relation \eqref{cDg}.

\bigskip

The previous analysis allows for a complete study of this quantum horizontal projector. Consider a connection 1-form $\tilde{\omega}(\omega_{z})=\omega_{z}+U\omega_{-}+V\omega_{+}=\omega_{z}+\mathrm{a}$ with $U\in\cl_{2}^{(0)}$ and $V\in\,\cl_{-2}^{(0)}$ as in equation \eqref{Pi3}. The inverse of the multiplicative map $\mathit{m}_{\omega}$ -- the map  $\mathit{m}_{\omega}^{-1}:\Omega(\SU)\to\mathfrak{hor}(\SU)\otimes\Omega(\U(1))_{\mathrm{inv}}$ -- as well as the horizontal projector are given  on 0-forms and 1-forms by:
\beq
\begin{array}{lcl}
\mathit{m}_{\omega}^{-1}(x)=x\otimes1& \Rightarrow& h_{\omega}(x)=x; \\
\mathit{m}_{\omega}^{-1}(x\,\omega_{\pm})=x\,\omega_{\pm}\otimes1 & \Rightarrow & h_{\omega}(x\,\omega_{\pm})=x\,\omega_{\pm}, \\ 
\mathit{m}_{\omega}^{-1}(x\,\omega_{z})=(-x\,\mathrm{a}\otimes1)+(x\otimes\omega_{z}) & \Rightarrow &h_{\omega}(x\,\omega_{z})=-x\,\mathrm{a}
\end{array}
\label{homega01}
\eeq
with $x\in\,\ASU$. This means that one has $\mathfrak{D}\phi=D\phi$ where $\phi\in\,\ASU$ with respect to the covariant derivative defined in \eqref{Dom}. On higher order exterior forms one has: 
\beq
\begin{array}{lcl}
\mathit{m}_{\omega}^{-1}(x\,\omega_{-}\wedge\omega_{+})=x\,\omega_{-}\wedge\omega_{+}\otimes 1 & \Rightarrow & h_{\omega}(x\,\omega_{-}\wedge\omega_{+})=x\,\omega_{-}\wedge\omega_{+}, \\
\mathit{m}_{\omega}^{-1}(x\,\omega_{+}\wedge\omega_{z})=(-x\omega_{+}\wedge\mathrm{a}\otimes 1)+(x\,\omega_{+}\otimes\omega_{z}) & \Rightarrow & h_{\omega}(x\,\omega_{+}\wedge\omega_{z})=-x\,\omega_{+}\wedge\mathrm{a}=x\,U\,\omega_{-}\wedge\omega_{+}, \\
\mathit{m}_{\omega}^{-1}(x\,\omega_{-}\wedge\omega_{z})=(-x\,\omega_{-}\wedge\mathrm{a}\otimes 1)+(x\,\omega_{-}\otimes\omega_{z}) & \Rightarrow & h_{\omega}(x\,\omega_{-}\wedge\omega_{z})=-x\,\omega_{-}\wedge\mathrm{a}=-q^{2}x\,V\,\omega_{-}\wedge\omega_{+}, \\
\mathit{m}_{\omega}^{-1}(x\,\omega_{-}\wedge\omega_{+}\wedge\omega_{z})=x\,\omega_{-}\wedge\omega_{+}\otimes\omega_{z} & \Rightarrow & h_{\omega}(x\,\omega_{-}\wedge\omega_{+}\wedge\omega_{z})=0.
\end{array}
\label{homega23}
\eeq
Recalling the analysis in remark \ref{rempi}, it is important to stress that the projector $h_{\omega}$ from \eqref{hproj} is well defined on the exterior algebra $\Omega(\SU)$ for any choice of the connection, and defines a covariant derivative $\mathfrak{D}:\Omega^{k}(\SU)\to\Omega^{k+1}_{\mathrm{hor}}(\SU)$ which reduces to the operators \eqref{Dom} on 0-forms and \eqref{D12} on 1-forms. The last equation out of \eqref{homega23} shows also that $\mathfrak{D}:\Omega^{2}(\SU)\to0$.

\begin{rema}
\label{refD}
Is the horizontal projector $h_{\omega}$ defined in \eqref{hproj} the only well-defined horizontal covariant projector operator whose domain coincides with $\Omega(\SU)$ and whose range is $\mathfrak{hor}(\SU)\subset\Omega(\SU)$, such that the associated horizontal projection of the exterior derivative \eqref{coDeD} reduces to the well established operator $D:\cl_{n}^{(k)}\to\cl_{n}^{(k+1)}$ given in 
\eqref{Dphi},\eqref{D12}? The answer is no. To be definite, consider the operator  $h^{\prime}_{\omega}:\Omega(\SU)\to\Omega_{\mathrm{hor}}(\SU)$  given by:
\begin{align}
&h^{\prime}_{\omega}(x)=x;\nn \\
&h^{\prime}_{\omega}(x\,\omega_{\pm})=x\,\omega_{\pm}, \nn \\
&h^{\prime}_{\omega}(x\,\omega_{z})=-x\,\mathrm{a},
\label{hpome01}
\end{align}
so to coincide with the projector $h_{\omega}$ \eqref{homega01} on 0-forms and 1-forms, and:
\begin{align}
&h^{\prime}_{\omega}(x\,\omega_{-}\wedge\omega_{+})=x\,\omega_{-}\wedge\omega_{+}, \nn \\
&h^{\prime}_{\omega}(x\,\omega_{+}\wedge\omega_{z})=q^{4}x\,\mathrm{a}\wedge\omega_{+}=q^{4}xU\,\omega_{-}\wedge\omega_{+}, \nn \\
&h^{\prime}_{\omega}(x\,\omega_{-}\wedge\omega_{z})=q^{-4}x\,\mathrm{a}\wedge\omega_{-}=-q^{-2}xV\,\omega_{-}\wedge\omega_{+};\nn \\
&h^{\prime}_{\omega}(x\,\omega_{-}\wedge\omega_{+}\wedge\omega_{z})=0.
\label{hpome23}
\end{align}
It is clear that the operator $\mathfrak{D}^{\prime}=h^{\prime}_{\omega}\dd:\Omega^{k}(\SU)\to\Omega_{\mathrm{hor}}^{k+1}(\SU)$ defines a consistent covariant derivative on the whole exterior algebra on the total space algebra of the quantum Hopf bundle, which  reduces to the operator $\mathfrak{D}$ from \eqref{coDeD} when restricted to  horizontal elements $\cl_{n}^{(k)}\subset\Omega^{k}(\SU)$. Both the operators $\mathfrak{D}, \mathfrak{D}^{\prime}$ coincide in the classical limit with the covariant derivative  on the classical Hopf bundle \eqref{cDg} presented in section \ref{limcla}.

The last step is to understand from where it is possible to trace  the origin of such a projector $h^{\prime}_{\omega}$ back. It is easy to see that the isomorphism $\mathit{m}_{\omega}^{-1}$ coming from \eqref{defmo} can be recovered as the choice of a specific left $\ASU$-module basis for the exterior algebra $\Omega(\SU)$, namely
\begin{align}
&\Omega(\SU)\,\simeq\,\ASU\{ 1\oplus\omega_{-}\oplus\omega_{+}\oplus\tilde{\omega}(\omega_{z})\}\nn\\
&\qquad\qquad\qquad\qquad\oplus\ASU\{ (\omega_{-}\wedge\omega_{+})\oplus(\omega_{-}\wedge\tilde{\omega}(\omega_{z}))\oplus( \omega_{+}\wedge\tilde{\omega}(\omega_{z}))\oplus(\omega_{-}\wedge\omega_{+}\wedge\tilde{\omega}(\omega_{z}))\},
\label{iba1}
\end{align}
while the horizontal projection obviously annihilates all the coefficients associated to exterior forms having the connection 1-form $\tilde{\omega}(\omega_{z})$ as a term. The projector $h^{\prime}_{\omega}$ in \eqref{hpome01},\eqref{hpome23} comes from the choice of a different left $\ASU$-module basis of $\Omega(\SU)$, that is setting -- as analogue of \eqref{iba1} --  the isomorphism
\begin{align}
&\Omega(\SU)\,\simeq\,\ASU\{ 1\oplus\omega_{-}\oplus\omega_{+}\oplus\tilde{\omega}(\omega_{z})\}\nn\\
&\qquad\qquad\qquad\qquad\oplus\ASU\{ (\omega_{-}\wedge\omega_{+})\oplus(\tilde{\omega}(\omega_{z})\wedge\omega_{-})\oplus (\tilde{\omega}(\omega_{z})\wedge\omega_{+})\oplus(\omega_{-}\wedge\omega_{+}\wedge\tilde{\omega}(\omega_{z}))\}.
\label{iba2}
\end{align}
and then defining $h_{\omega}^{\prime}$ as the projector whose nucleus is given as the  left $\ASU$-module spanned by $\{\tilde{\omega}(\omega_{z}),\tilde{\omega}(\omega_{z})\wedge\omega_{\pm},\omega_{-}\wedge\omega_{+}\wedge\tilde{\omega}(\omega_{z})\}$. An explicit computation shows that
$$
\omega_{-}\wedge\tilde{\omega}(\omega_{z})=(q^{2}-q^{-2})V\omega_{-}\wedge\omega_{+}-q^{-4}\tilde{\omega}(\omega_{z})\wedge\omega_{-}\quad\Rightarrow\quad \ker\,h_{\omega}\neq\ker\,h^{\prime}_{\omega}:
$$
the two projectors are not equivalent, being equivalent if and only if the connection is the monopole connection.
\end{rema}





\subsection*{Acknowledgements}
Some days ago, reading once more this manuscript, I felt that for almost any sentence of it I am indebted to one and to all of the travelmates I had during the last year. 
This paper has been originated and developed as a part of a more general research project with G.Landi: I want to thank him for  his support, guidance  and feedback. I should like to thank M.Marcolli, who suggested me to write it, and S.Albeverio: with them I often discussed about many of the themes here described. I should like to thank G.Marmo, L.Cirio and C.Pagani, who read a draft of the paper and  made my understanding of many important details clearer, and G.Dell'Antonio, who encouraged me to write it.  And I should like to thank the referee: his report was precious and helped me to improve it.

It is a pleasure to thank the Max-Planck-Institut f\"ur Mathematik in Bonn and the Hausdorff Center for Mathematics at the University Bonn for their invitation, the Foundation Blanceflor Boncompagni-Ludovisi n\'ee Bildt (Stockholm) for the support,  the Joint Institut Research for Nuclear Physics in Dubna-Moscow: A.Motovilov was a wonderful host.

\end{document}